\documentclass[12pt,a4paper,reqno]{amsart}
\usepackage{amsmath}
\usepackage{amsfonts}
\usepackage{amssymb}
\numberwithin{equation}{section}

\theoremstyle{plain}

\newtheorem{theorem}[subsection]{Theorem}
\newtheorem{proposition}[subsection]{Proposition}
\newtheorem{lemma}[subsection]{Lemma}

\newtheorem{conjecture}[subsection]{Conjecture}

\theoremstyle{definition}

\newtheorem{definition}[subsection]{Definition}

\renewcommand{\leq}{\leqslant}
\renewcommand{\geq}{\geqslant}

\newsavebox{\proofbox}
\savebox{\proofbox}{\begin{picture}(7,7)%
  \put(0,0){\framebox(7,7){}}\end{picture}}

\newcommand{\md}[1]{\ensuremath{(\mbox{mod}\, #1)}}
\newcommand{\mdsub}[1]{\ensuremath{(\mbox{\scriptsize mod}\, #1)}}

\def\vconst{\ensuremath \nu_{\mbox{\scriptsize const}}}
\def\vconstlem{\ensuremath \nu_{\mbox{\scriptsize \textup{const}}}}

\def\funif{\ensuremath f_{U}}
\def\fanti{\ensuremath f_{U^\perp}}

\def\B{{\mathcal B}}
\def\D{{\mathcal D}}
\def\Z{{\mathbb Z}}
\def\E{{\mathbb E}}
\def\C{{\mathbb C}}
\def\R{{\mathbb R}}

\def\Q{{\mathbb Q}}

\def\x{{\bf x}}
\def\eps{{\varepsilon}}

\def\w{w}

\def\vs{\vspace{11pt}}
\def\ni{\noindent}

\def\proof{\noindent\textit{Proof. }}
\def\endproof{\hfill{\usebox{\proofbox}}}

\def\emph#1{{\it #1}}
\def\textbf#1{{\bf #1}}

\begin{document}

\title{The primes contain arbitrarily long arithmetic progressions}

\author{Ben Green}
\address{School of Mathematics\\
University Walk\\
Bristol\\
BS8 1TW
}
\email{b.j.green@bristol.ac.uk}

\author{Terence Tao}
\address{Department of Mathematics\\University of California at Los Angeles\\ Los Angeles CA 90095}

\email{tao@math.ucla.edu}

\thanks{While this work was carried out the first author was a PIMS postdoctoral fellow at the University of British Columbia, Vancouver, Canada.  The second author was
a Clay Prize Fellow and was supported by a grant from the Packard Foundation.}

\subjclass{11N13, 11B25, 374A5}

\begin{abstract}
We prove that there are arbitrarily long arithmetic progressions of primes.\vs  

\ni There are three
major ingredients.  The first is Szemer\'edi's theorem, which asserts that any subset of the integers of positive density 
contains progressions of arbitrary length.  The second, which is the main new ingredient of this paper,
is a certain transference principle. This allows us to deduce from Szemer\'edi's theorem that any subset of a sufficiently 
pseudorandom set (or measure) of positive \emph{relative} density 
contains progressions of arbitrary length.  The third ingredient is a recent result of Goldston
and Y{\i}ld{\i}r{\i}m, which we reproduce here. Using this, one may place (a large fraction of) the primes 
inside a pseudorandom set of ``almost primes'' (or more precisely, a pseudorandom measure concentrated on almost primes)
with positive relative density. 
\end{abstract}

\maketitle

\section{Introduction}
     
\ni It is a well-known conjecture that there are arbitrarily long arithmetic progressions of prime numbers. The conjecture is best described as ``classical'', or maybe even ``folklore''. In Dickson's \textit{History} it is stated that around 1770 Lagrange and Waring investigated how large the common difference of an arithmetic progression of $L$ primes must be, and it is hard to imagine that they did not at least wonder whether their results were sharp for all $L$.\vs

\ni It is not surprising that the conjecture should have been made, since a simple heuristic based on the prime number theorem would suggest that there are $\gg N^2/\log^k N$ $k$-tuples of primes $p_1,\dots,p_k$ in arithmetic progression, each $p_i$ being at most $N$. Hardy and Littlewood \cite{hardy-littlewood}, in their famous paper of 1923, advanced a very general conjecture which, as a special case, contains the hypothesis that the number of such $k$-term progressions is asymptotically $C_kN^2/\log^k N$ for a certain explicit numerical factor $C_k > 0$ (we do not come close to establishing this conjecture here, obtaining instead
a lower bound $(\gamma(k) + o(1)) N^2 / \log^k N$ for some very small $\gamma(k) > 0$).\vs

\ni The first theoretical progress on these conjectures was made by van der Corput \cite{van-der-corput} (see also \cite{chowla}) who, in 1939, used Vinogradov's method of prime number sums to establish the case $k = 3$, that is to say that there are infinitely many triples of primes in arithmetic progression.  However, the question of longer arithmetic progressions seems to have remained completely open (except for upper bounds), even for $k=4$.  On the other hand, it has been known for some time that better results can be obtained if one replaces the primes with
a slightly larger set of \emph{almost primes}. The most impressive such result is due to Heath-Brown \cite{heath-brown1}. He showed that there are infinitely many 4-term progressions consisting of three primes and a number which is either prime or a product of two primes. In a somewhat different direction, let us mention the beautiful results of Balog \cite{balog1,balog2}. Among other things he shows that for any $m$ there are $m$ distinct primes $p_1,\dots,p_m$ such that all of the averages $\frac{1}{2}(p_i + p_j)$ are prime.\vs

\ni The problem of finding long arithmetic progressions in the primes has also attracted the interest of computational mathematicians. At the time of writing the longest known arithmetic progression of primes is of length 23, and was found in 2004 by Markus Frind, Paul Underwood, and Paul Jobling:
\[ 56211383760397 + 44546738095860 k; \quad k = 0, 1, \ldots, 22.\]
An earlier arithmetic progression of primes of length 22 was found by Moran, Pritchard and Thyssen \cite{moran-pritchard-thyssen}:
\[ 11410337850553 + 4609098694200k; \quad k = 0, 1, \ldots, 21.\]

\ni Our main theorem resolves the above conjecture.

\begin{theorem}\label{mainthm} The prime numbers contain infinitely many arithmetic progressions of length $k$ for all $k$.  
\end{theorem}

\ni In fact, we can say something a little stronger:

\begin{theorem}[Szemer\'edi's theorem in the primes]\label{sz-primes} Let $A$ be any subset of the prime numbers of positive relative upper density, thus $\limsup_{N \to \infty} \pi(N)^{-1}|A \cap [1,N]| > 0$, where $\pi(N)$ denotes the number of primes less than or equal to $N$.  Then $A$ contains infinitely many arithmetic progressions of length $k$ for all $k$.
\end{theorem}

\ni If one replaces ``primes'' in the statement of Theorem \ref{sz-primes} by the set of all positive integers $\Z^+$, then this is a famous theorem of Szemer\'edi \cite{szemeredi}. The special case $k = 3$ of Theorem \ref{sz-primes} was recently established by the first author \cite{green} using methods of Fourier analysis.  In contrast, our
methods here have a more ergodic theory flavour and do not involve much Fourier analysis (though the argument does rely on Szemer\'edi's theorem which can be proven by either combinatorial, ergodic theory, or Fourier analysis arguments).  We also remark that if the primes were replaced
 by a random subset of the integers, with density at least $N^{-1/2+\eps}$ on each interval $[1,N]$, then the $k=3$ case of the above theorem was established in \cite{klr}.  
\vs

\ni\textit{Acknowledgements} The authors would like to thank Jean Bourgain, Enrico Bombieri, Tim Gowers, Bryna Kra, Elon Lindenstrauss, Imre Ruzsa, Roman Sasyk, Peter Sarnak and Kannan Soundararajan for helpful conversations. We are particularly indebted to Andrew Granville for drawing our attention to the work of Goldston and Y{\i}ld{\i}r{\i}m, and to Dan Goldston for making the preprint \cite{goldston-yildirim} available. We are also indebted to Yong-Gao Chen and his students, Bryna Kra, Jamie Radcliffe, Lior Silberman and Mark Watkins for corrections to earlier versions of the manuscript.  We are 
particularly indebted to the anonymous referees for a very thorough reading and many helpful corrections and suggestions, which have
been incorporated into this version of the paper.
Portions of this work were completed while the first author was visiting UCLA and Universit\'e de Montr\'eal, and he would like to thank these institutions for their hospitality. He would also like to thank Trinity College, Cambridge for support over several years.

\section{An outline of the proof}

\ni Let us start by stating Szemer\'edi's theorem properly. In the introduction we claimed that it was a statement about sets of integers with positive upper density, but there are other equivalent formulations.  A ``finitary'' version of the theorem is as follows.

\begin{proposition}[Szemer\'edi's theorem]\label{sz-original}\cite{szemeredi-4,szemeredi} Let $N$ be a positive integer and let $\mathbb{Z}_N := \mathbb{Z}/N\mathbb{Z}$.\footnote{We will retain this notation throughout the paper, thus $\Z_N$ will never refer to the $N$-adics. We always assume for convenience that $N$ is prime. It is very convenient to work in $\mathbb{Z}_N$, rather than the more traditional $[-N,N]$, since we are free to divide by $2,3,\dots,k$ and it is possible to make linear changes of variables without worrying about the ranges of summation.  There is a slight price to pay for this, in that one must now address some ``wraparound'' issues when identifying $\Z_N$ with a subset of the integers, but these will be easily dealt with.}  
 Let $\delta > 0$ be a fixed positive real number, and let $k \geq 3$ be an integer. Then there is a minimal $N_0(\delta,k) < \infty$ with the following property. If $N \geq N_0(\delta,k)$ and $A \subseteq \mathbb{Z}_N$ is any set of cardinality at least $\delta N$, then $A$ contains an arithmetic progression of length $k$.
\end{proposition}

\ni Finding the correct dependence of $N_0$ on $\delta$ and $k$ (particularly $\delta$) is a famous open problem. It was a great breakthrough when Gowers \cite{gowers-4,gowers} showed that 
\[ N_0(\delta,k) \leq 2^{2^{\delta^{-c_k}}},\] where $c_k$ is an explicit constant (Gowers obtains $c_k = 2^{2^{k+9}}$). It is possible that a new proof of Szemer\'edi's theorem could be found, with sufficiently good bounds that Theorem \ref{mainthm} would follow immediately. To do this one would need something just a little weaker than\begin{equation}\label{too-hard} N_0(\delta,k) \leq 2^{c_k\delta^{-1}}\end{equation} (there is a trick, namely passing to a subprogression of common difference $2\times 3 \times 5 \times \dots \times \w(N)$ for appropriate $\w(N)$, which allows one to consider the primes as a set of density essentially $\log \log N/\log N$ rather than $1/\log N$; we will use a variant of this ``$W$-trick'' later in this paper to eliminate local irregularities arising from small divisors).  In our proof of Theorem \ref{sz-primes}, we will need to use Szemer\'edi's theorem, but we will not need any quantitative estimimates on $N_0(\delta,k)$.\vs

\ni Let us state, for contrast, the best known lower bound which is due to Rankin \cite{rankin} (see also Lacey-{\L}aba \cite{laba-lacey}):
\[ N_0(\delta,k) \geq \exp(C(\log 1/\delta)^{1 + \lfloor \log_2 (k-1)\rfloor}).\]

\ni At the moment it is clear that a substantial new idea would be required to obtain a result of the strength \eqref{too-hard}. In fact, even for $k = 3$ the best bound is $N_0(\delta,3) \leq 2^{C\delta^{-2}\log(1/\delta)}$, a result of Bourgain \cite{bourgain-triples}. The hypothetical bound \eqref{too-hard} is closely related to the following very open conjecture of Erd\H{o}s:
\begin{conjecture}[Erd\H{o}s conjecture on arithmetic progressions] Suppose that $A = \{a_1 < a_2 < \dots\}$ is an infinite sequence of integers such that $\sum 1/a_i = \infty$. Then $A$ contains arbitrarily long arithmetic progressions.
\end{conjecture}
\ni This would imply Theorem \ref{mainthm}.\vs

\ni We do not make progress on any of these issues here. In one sentence, our argument can be described instead as a \textit{transference principle} which allows us to deduce Theorems \ref{mainthm} and \ref{sz-primes} from Szemer\'edi's theorem, regardless of what bound we know for $N_0(\delta,k)$; in fact we prove a more general statement in Theorem \ref{main} below. Thus, in this paper, we must assume Szemer\'edi's theorem. However with this one (rather large!) caveat\footnote{We will also require some standard facts from analytic number theory such as the prime number theorem, Dirichlet's theorem on primes in arithmetic progressions, and the classical zero-free region for the Riemann $\zeta$-function (see Lemma \ref{basic-zeta-facts}).} our paper is self-contained.\vs

\ni Szemer\'edi's theorem can now be proved in several ways. The original proof of Szemer\'edi \cite{szemeredi-4,szemeredi} was combinatorial. In 1977, Furstenberg made a very important breakthrough by providing an ergodic-theoretic proof \cite{furst}. Perhaps surprisingly for a result about primes, our paper has at least as much in common with the ergodic-theoretic approach as it does with the harmonic analysis approach of Gowers. We will use a language which suggests this close connection, without actually relying explicitly on any ergodic-theoretical concepts\footnote{It has become clear that there is a deep connection between harmonic analysis (as applied to solving linear equations in sets of integers) and certain parts of ergodic theory. Particularly exciting is the suspicion that the notion of a $k$-step nilsystem, explored in many ergodic-theoretical works (see e.g. \cite{host-kra1,host-kra2,host-kra3,ziegler}), might be analogous to a kind of ``higher order Fourier analysis'' which could be used to deal with systems of linear equations that cannot be handled by conventional Fourier analysis (a simple example being the equations $x_1 + x_3 = 2x_2$, $x_2 + x_4 = 2x_3$, which define an arithmetic progression of length 4). We will not discuss such speculations any further here, but suffice it to say that much is left to be understood.}.
In particular we shall always remain in the finitary setting of $\Z_N$, in contrast to the standard ergodic theory framework in which one takes weak limits (invoking the axiom of choice) to pass to an infinite measure-preserving system.
As will become clear in our argument, 
in the finitary setting one can still access many
tools and concepts from ergodic theory, but often one must incur error terms of the form $o(1)$ when one does so.\vs

\ni Here is another form of Szemer\'edi's theorem which suggests the ergodic theory analogy more closely. We use the conditional expectation notation $\mathbb{E}(f | x_i \in B)$ to denote the average of $f$ as certain variables $x_i$ range over the set $B$,
and $o(1)$ for a quantity which tends to zero as $N \rightarrow \infty$ (we will give more precise definitions later).

\begin{proposition}[Szemer\'edi's theorem, again]\label{sz} Write $\vconstlem : \mathbb{Z}_N \rightarrow \R^+$ for the constant function $\vconstlem \equiv 1$. Let $0 < \delta \leq 1$ and $k \geq 1$ be fixed.  Let $N$ be a large integer parameter, and let
$f: \mathbb{Z}_N \to \R^+$ be a non-negative function obeying the bounds
\begin{equation}\label{f-bound} 0 \leq f(x) \leq \vconstlem(x) \hbox{ for all } x \in \mathbb{Z}_N
\end{equation}
and
\begin{equation}\label{pos-dens} \E(f(x) | x \in \mathbb{Z}_N) \geq \delta.\end{equation}
Then we have
$$ \E( f(x) f(x+r) \ldots f(x+(k-1)r) | x,r \in \mathbb{Z}_N) \geq c(k,\delta) - o_{k,\delta}(1)$$
for some constant $c(k,\delta) > 0$ which does not depend on $f$ or $N$.
\end{proposition}
\noindent\textit{Remark.} Ignoring for a moment the curious notation for the constant function $\vconst$, there are two main differences between this and Proposition \ref{sz-original}. One is the fact that we are dealing with functions rather than sets: however, it is easy to pass from sets to functions, for instance by probabilistic arguments. Another difference, if one unravels the $\mathbb{E}$ notation, is that we are now asserting the existence of $\gg N^2$ arithmetic progressions, and not just one. Once again, such a statement can be deduced from Proposition \ref{sz-original} with some combinatorial trickery (of a less trivial nature this time -- the argument was first worked out by Varnavides \cite{varnavides}).  A direct proof of Proposition \ref{sz} can be found in \cite{tao:ergodic}.  A formulation of Szemer\'edi's theorem similar to this one was also used by Furstenberg \cite{furst}.  Combining this argument with the one in Gowers gives an explicit bound on $c(k,\delta)$ of the form $c(k,\delta) \geq \exp(-\exp(\delta^{-c_k}))$ for some $c_k > 0$.\vs

\ni Now let us abandon the notion that $\nu$ is the constant function. We say that $\nu : \mathbb{Z}_N \rightarrow \R^+$ is a \textit{measure}\footnote{The term \emph{normalized probability density} might be more accurate here, but \emph{measure} has the advantage of brevity.
One may think of $\vconstlem$ as the uniform probability distribution on $\Z_N$, and $\nu$ as some other probability distribution which can concentrate on a subset of $\Z_N$ of very small density (e.g. it may concentrate on the ``almost primes'' in $[1,N]$).} if 
\begin{equation}\label{numean}
\mathbb{E}(\nu) = 1 + o(1). 
\end{equation}
We are going to exhibit a class of measures, more general than the constant function $\vconst$, for which Proposition \ref{sz} still holds. These measures, which we will call \textit{pseudorandom}, will be ones satisfying two conditions called the \textit{linear forms condition} and the \textit{correlation condition}. These are, of course, defined formally below, but let us remark that they are very closely related to the ergodic-theory notion of weak-mixing. It is perfectly possible for a ``singular'' measure - for instance, a measure for which $\E(\nu^2)$ grows like a power of $\log N$ - to be pseudorandom. Singular measures are the ones that will be of interest to us, since they generally support rather sparse sets. This generalisation of Proposition \ref{sz} is Proposition \ref{main} below.\vs

\ni Once Proposition \ref{main} is proved, we turn to the issue of finding primes in AP. A possible choice for $\nu$ would be $\Lambda$, the von Mangoldt function (this is defined to equal $\log p$ at $p^m$, $m = 1,2,\dots,$ and $0$ otherwise). Unfortunately, verifying the linear forms condition and the correlation condition for the von Mangoldt function (or minor variants thereof) is strictly harder than proving that the primes contain long arithmetic progressions; indeed, this task is comparable in difficulty to the notorious Hardy-Littlewood prime tuples conjecture, for which our methods here yield no progress.\vs

\ni However, all we need is a measure $\nu$ which (after rescaling by at most a constant factor) majorises $\Lambda$ pointwise. Then, \eqref{pos-dens} will be satisfied with $f = \Lambda$. Such a measure is provided to us\footnote{Actually, there is an extra technicality which is caused by the very irregular distribution of primes in arithmetic progressions to small moduli (there are no primes congruent to $4 \md{6}$, for example). We get around this using something which we refer to as the $W$-trick, which basically consists of restricting the primes to the arithmetic progression $n \equiv 1 \md{W}$, where $W = \prod_{p < \w(N)} p$ and $w(N)$ tends slowly to infinity with $N$. Although this looks like a trick, it is actually an extremely important feature of that part of our argument which concerns primes.} by recent work of Goldston and Y{\i}ld{\i}r{\i}m \cite{goldston-yildirim} concerning the size of gaps between primes. The proof that the linear forms condition and the correlation condition are satisfied is heavily based on their work, so much so that parts of the argument are placed in an appendix.\vs

\ni The idea of using a majorant to study the primes is by no means new -- indeed in some sense sieve theory is precisely the study of such objects. For another use of a majorant in an additive-combinatorial setting, see \cite{ramare,ramare-ruzsa}. \vs

\ni It is now timely to make a few remarks concerning the proof of Proposition \ref{main}. It is in the first step of the proof that our original investigations began, when we made a close examination of Gowers' arguments. If $f : \mathbb{Z}_N \rightarrow \R^+$ is a function then the normalised count of $k$-term arithmetic progressions
\begin{equation}\label{apcount}  \E( f(x) f(x+r) \ldots f(x+(k-1)r) | x,r \in \mathbb{Z}_N)\end{equation} 
is closely controlled by certain norms $\Vert \cdot \Vert_{U^d}$, which we would like to call the \textit{Gowers uniformity norms}\footnote{Analogous objects have recently surfaced in the genuinely ergodic-theoretical work of Host and Kra \cite{host-kra1,host-kra2,host-kra3} concerning non-conventional ergodic averages, thus enhancing the connection between ergodic theory and additive number theory.}. They are defined in \S \ref{sec5}. The formal statement of this fact can be called a generalised von Neumann theorem. Such a theorem, in the case $\nu = \vconst$, was proved by Gowers \cite{gowers} as a first step in his proof of Szemer\'edi's theorem, using $k-2$ applications of the Cauchy-Schwarz inequality. In Proposition \ref{vn} we will prove a generalised von Neumann theorem relative to an arbitrary pseudorandom measure $\nu$. Our main tool is again the Cauchy-Schwarz inequality. We will use the term \textit{Gowers uniform} loosely to describe a function which is small in some $U^d$ norm. This should not be confused with the term pseudorandom, which will be reserved for measures on $\mathbb{Z}_N$.\vs

\ni Sections \ref{sec6}-\ref{sec7} are devoted to concluding the proof of Proposition \ref{main}. Very roughly the strategy will be to decompose the function $f$ under consideration into a Gowers uniform component plus a bounded ``Gowers anti-uniform'' object (plus a negligible error). The notion\footnote{We note that Gowers uniformity, which is a measure of ``randomness'', ``uniform distribution'', or ``unbiasedness'' in a function should not be confused with the very different notion of uniform boundedness.  Indeed, in our arguments, the Gowers uniform functions will be highly unbounded, whereas the Gowers anti-uniform functions will be uniformly bounded.  Anti-uniformity can in fact be viewed as a measure of ``smoothness'', ``predictability'', ``structure'', or ``almost periodicity''.} of Gowers anti-uniformity is captured using the dual norms $(U^d)^*$, whose properties are laid out in \S \ref{sec6}.\vs

\ni The contribution of the Gowers-uniform part to the count \eqref{apcount} will be negligible\footnote{Using the language of ergodic theory, we are essentially claiming that the Gowers anti-uniform functions form a characteristic factor for the expression \eqref{apcount}.  The point is that even though $f$ is not necessarily bounded uniformly, the fact that it is bounded pointwise by a pseudorandom measure $\nu$ allows us to conclude that the \emph{projection} of $f$ to the Gowers anti-uniform component is bounded, at which point we can invoke the standard Szemer\'edi theorem.} by the generalised von Neumann theorem. The contribution from the Gowers anti-uniform component will be bounded from below by Szemer\'edi's theorem in its traditional form, Proposition \ref{sz}. 

\section{Pseudorandom measures}\label{sec3}

\ni In this section we specify exactly what we mean by a pseudorandom measure on $\mathbb{Z}_N$. First, however, we set up some notation. We fix the length $k$ of the arithmetic progressions we are seeking.  $N = |\mathbb{Z}_N|$ will always be assumed to be prime and large (in particular, we can invert any of the numbers $1,\ldots, k$ in $\Z_N$), and we will write $o(1)$ for a quantity that tends to zero as $N \rightarrow \infty$. We will write $O(1)$ for a bounded quantity. Sometimes quantities of this type will tend to zero (resp. be bounded) in a way that depends on some other, typically fixed, parameters. If there is any danger of confusion as to what is being proved, we will indicate such dependence using subscripts, thus for instance $O_{j,\eps}(1)$ denotes a quantity whose magnitude is bounded by $C(j,\eps)$ for some quantity $C(j,\eps) > 0$ depending only $j,\eps$.  Since every quantity in this paper will depend on $k$, however, we will not bother indicating the $k$ dependence throughout this paper.  As is customary
we often abbreviate $O(1) X$ and $o(1) X$ as $O(X)$ and $o(X)$ respectively for various non-negative quantities $X$.\vs

\ni If $A$ is a finite non-empty set (for us $A$ is usually just $\mathbb{Z}_N$) and $f: A \to \R$ is a function,
we write $\E( f ) := \E( f(x) | x \in A )$ for the average value of $f$, that is to say
\[ \E(f) := \frac{1}{|A|} \sum_{x \in A} f(x).\] Here, as is usual, we write $|A|$ for the cardinality of the set $A$. More generally, if $P(x)$ is any statement concerning an element of $A$ which is true for at least one $x \in A$,
we define
\[ \E(f(x) | P(x)) := \frac{ \sum_{x \in A: P(x)} f(x)}{|\{ x \in A: P(x) \}|}.\]
This notation extends to functions of several variables in the obvious manner.  We now define two notions of randomness
for a measure, which we term the linear forms condition and the correlation condition.

\begin{definition}[Linear forms condition]\label{linear-forms-condition}
Let $\nu : \mathbb{Z}_N \rightarrow \mathbb{R}^{+}$ be a measure. Let $m_0,t_0$ and $L_0$ be small positive integer parameters. Then we say that $\nu$ satisfies the $(m_0,t_0,L_0)$-linear forms condition if the following holds. Let $m \leq m_0$ and $t \leq t_0$ be arbitrary, and suppose that $( L_{ij} )_{1 \leq i \leq m, 1 \leq j \leq t}$ are arbitrary rational numbers with numerator and denominator at most $L_0$ in absolute value, and that $b_i$, $1 \leq i \leq m$, are arbitrary elements of $\mathbb{Z}_N$.  For $1 \leq i \leq m$, let $\psi_i: \Z_N^t \to \Z_N$ be the linear forms $\psi_i(\mathbf{x}) = \sum_{j = 1}^t L_{ij}x_j + b_i$, where $\mathbf{x} = (x_1,\dots,x_t) \in \mathbb{Z}_N^t$, and where the rational numbers $L_{ij}$ are interpreted as elements of $\Z_N$ in the usual manner (assuming $N$ is prime and larger than $L_0$). Suppose that as $i$ ranges over $1,\ldots,m$, the $t$-tuples $(L_{ij})_{1 \leq j \leq t} \in \Q^t$ are non-zero, and no $t$-tuple is a rational multiple of any other.
Then we have 
\begin{equation}\label{lfc}
 \mathbb{E}\left( \nu(\psi_1(\mathbf{x})) \dots \nu(\psi_m(\mathbf{x})) \;|\; \mathbf{x} \in \mathbb{Z}_N^t \right) = 1 + o_{L_0,m_0,t_0}(1).
\end{equation}
Note that the rate of decay in the $o(1)$ term is assumed to be uniform in the choice of $b_1, \ldots, b_m$.
\end{definition}

\ni\textit{Remarks.} It is the parameter $m_0$, which controls the number of linear forms, that is by far the most important, and will be kept relatively small. It will eventually be set equal to $k \cdot 2^{k-1}$.  Note that the $m=1$ case of the linear forms condition recovers the measure condition \eqref{numean}.  Other simple examples of the linear forms condition which we will encounter later are
\begin{equation}\label{example-1}
\E( \nu(x) \nu(x+h_1) \nu(x+h_2) \nu(x+h_1+h_2) \; | \; x,h_1,h_2 \in \Z_N) = 1 + o(1) \end{equation}
(here $(m_0, t_0, L_0) = (4,3,1)$);
\begin{equation}\label{example-2} 
\E \big( \nu(x + h_1) \nu(x + h_2) \nu(x + h_1 + h_2) \; | \; h_1 , h_2 \in \Z_N \big) = 1 + o(1)
\end{equation}
for all $x \in \Z_N$ (here $(m_0,t_0, L_0) = (3,2,1)$) and
\begin{eqnarray}\nonumber& & 
\E\bigg( \nu((x-y)/2) \nu((x-y+h_2)/2) \nu(-y) \nu(-y-h_1) \times \\ & & \nonumber
\qquad \times \; \nu((x-y')/2) \nu((x-y'+h_2)/2)  \nu(-y') \nu(-y'-h_1) \times \\ & &
\qquad \qquad \times \; \nu(x) \nu(x+h_1) \nu(x+h_2) \nu(x+h_1+h_2) \; 
\bigg| \; x,h_1,h_2,y,y' \in \Z_N\bigg) \nonumber \\ & & \qquad \qquad \qquad \qquad = 1 + o(1)\label{example-3}
\end{eqnarray} (here $(m_0, t_0, L_0) = (12,5,2)$).
For those readers familiar with the Gowers uniformity norms $U^{k-1}$ 
(which we shall discuss in detail later), the example \eqref{example-1} demonstrates that $\nu$ is close to 1 in the $U^2$ norm (see Lemma \ref{mua}). Similarly, the linear forms condition with appropriately many parameters implies that $\nu$ is close to $1$ in the $U^d$ norm, for any fixed $d \geq 2$. However, the linear forms condition is much stronger than simply asserting that $\Vert \nu - 1 \Vert_{U^d}$ is small for various $d$.\vs

\ni For the application to the primes, the measure $\nu$ will be constructed using truncated divisor sums, and the linear forms condition will be deduced from some arguments of Goldston and Y{\i}ld{\i}r{\i}m.  From a probabilistic point of view, the linear
forms condition is asserting a type of joint independence between the ``random variables'' $\nu(\psi_j(\mathbf{x}))$; in the application to the primes,
$\nu$ will be concentrated on the ``almost primes'', and the linear forms condition is then saying that the events ``$\psi_j(\mathbf{x})$ is almost prime'' are essentially independent of each other as $j$ varies\footnote{This will only be true after first eliminating some local correlations in the almost primes arising from small divisors.  This will be achieved by a simple ``$W$-trick'' which we will come to later in this paper.}.  

\begin{definition}[Correlation condition]\label{correlation-condition}
Let $\nu : \mathbb{Z}_N \rightarrow \mathbb{R}^{+}$ be a measure, and let $m_0$ be a positive integer parameter. We say that $\nu$ satisfies the $m_0$-correlation condition if for every $1 < m \leq m_0$
there exists a weight 
function $\tau = \tau_{m}: \mathbb{Z}_N \to \R^+$ which obeys the moment conditions
\begin{equation}\label{eq3.1}
 \E( \tau^q ) = O_{m,q}(1)
\end{equation}
for all $1 \leq q < \infty$
and such that
\begin{equation}\label{eq3.2}
\E( \nu(x+h_1) \nu(x+h_2) \ldots \nu(x+h_m) \; |\;  x \in \mathbb{Z}_N)
\leq \sum_{1 \leq i < j \leq m} \tau(h_i-h_j)
\end{equation}
for all $h_1, \ldots, h_m \in \Z_N$ (not necessarily distinct).
\end{definition}

\ni\emph{Remarks.}  The condition \eqref{eq3.2} may look a little strange, since if $\nu$ were to be chosen randomly then we would expect such a condition to hold with $1 + o(1)$ on the right-hand side, at least when $h_1,\ldots,h_m$ are distinct.  Note that one cannot use the linear
forms condition to control the left-hand side of \eqref{eq3.2} because the linear components of the forms $x + h_j$ are all the same.
 The correlation condition has been designed with the primes in mind\footnote{A simpler, but perhaps less interesting, model case occurs when one is trying to prove Szemer\'edi's theorem relative to a random subset of $\{1,\ldots,N\}$ of density $1/\log N$ (cf. \cite{klr}).  The pseudorandom weight $\nu$ would then be a Bernoulli random variable, with each $\nu(x)$ equal to $\log N$ with independent probability $1/\log N$ and equal to $0$ otherwise.  In such a case,
 we can (with high probability)
 bound the left-hand side of \eqref{eq3.2} more cleanly by $O(1)$ (and even obtain the asymptotic $1+o(1)$) 
 when the $h_j$ are distinct, and by $O(\log^m N)$ otherwise.}, 
 because in that case we must tolerate slight ``arithmetic'' nonuniformities. Observe, for example, that the number of $p \leq N$ for which $p - h$ is also prime is not bounded above by a constant times $N/\log^2 N$ if $h$ contains a very large number of prime factors, although such exceptions will of course be very rare and one still expects to have moment conditions such as \eqref{eq3.1}. It is phenomena like this which prevent us from assuming an $L^{\infty}$ bound for $\tau$.  While $m_0$ will be restricted to be small (in fact, equal to $2^{k-1}$), it will be important for us that there is no upper bound required on $q$ (which we will eventually need to be a very large function of $k$, but still independent of $N$ of course).  Since the correlation condition is an upper bound rather than an asymptotic, it is
 fairly easy to obtain; we shall prove it using the arguments of Goldston and Y{\i}ld{\i}r{\i}m (since we are using those methods in any case to prove
 the linear forms condition), but these upper bounds could also be obtained by more standard sieve theory methods.
 
\begin{definition}[Pseudorandom measures]\label{mu-pseudo-def}
Let $\nu : \mathbb{Z}_N \rightarrow \mathbb{R}^{+}$ be a measure. We say that $\nu$ is $k$-pseudorandom if it satisfies the $(k\cdot 2^{k-1},3k - 4,k)$-linear forms condition and also the $2^{k-1}$-correlation condition.
\end{definition}

\ni\emph{Remarks.} The exact values $k \cdot 2^{k-1}$, $3k-4$, $k$, $2^{k-1}$ of the parameters chosen here
are not too important; in our application to the primes, any quantities
which depend only on $k$ would suffice.  It can be shown that if $C = C_k > 1$ is any constant independent of $N$ and if $S \subseteq \mathbb{Z}_N$ is chosen at random, each $x \in \mathbb{Z}_N$ being selected to lie in $S$ independently at random with probability $1/\log^C N$, then (with high probability) the measure $\nu = \log^C N {\bf 1}_S$ is $k$-pseudorandom, and the Hardy-Littlewood prime tuples conjecture can be viewed as an assertion that the Von Mangoldt function is essentially of this form (once one eliminates the obvious obstructions to pseudorandomness coming from small prime divisors). While we will of course not attempt to establish this conjecture here, in \S \ref{sec8} we will construct pseudorandom measures which are concentrated on the \emph{almost primes} instead of the primes; this is of course consistent with the so-called ``fundamental lemma of sieve theory'',
but we will need a rather precise variant of this lemma due to Goldston and Y{\i}ld{\i}r{\i}m. \vs

\ni The function $\vconst \equiv 1$ is clearly $k$-pseudorandom for any $k$.  In fact the pseudorandom measures are star-shaped
around the constant measure:

\begin{lemma}\label{halfway}  Let $\nu$ be a $k$-pseudorandom measure.  Then $\nu_{1/2} := (\nu + \vconstlem)/2 = (\nu + 1)/2$ is also a $k$-pseudorandom measure \textup{(}though possibly with slightly different bounds in the $O()$ and $o()$ terms\textup{)}.
\end{lemma}

\proof It is clear that $\nu_{1/2}$ is non-negative and has expectation $1 + o(1)$.  To verify the linear forms condition \eqref{lfc}, we simply
replace $\nu$ by $(\nu+1)/2$ in the definition and expand as a sum of $2^m$ terms, divided by $2^m$.  Since each term can be verified to be
$1 + o(1)$ by the linear forms condition \eqref{lfc}, the claim follows.  The correlation condition is verified
in a similar manner.  (A similar result holds for $(1-\theta) \nu + \theta \vconst$ for any $0 \leq \theta \leq 1$, but we will not need to use
this generalization.)\endproof\vs

\ni The following result is one of the main theorems of the paper.  It asserts that for the purposes of Szemer\'edi's theorem (and ignoring $o(1)$ errors), there is no distinction
between a $k$-pseudorandom measure $\nu$ and the constant measure $\vconst$.

\begin{theorem}[Szemer\'edi's theorem relative to a pseudorandom measure]\label{main}
Let $k \geq 3$ and $0 < \delta \leq 1$  be fixed parameters. Suppose that $\nu: \mathbb{Z}_N \to \R^+$ is $k$-pseudorandom. Let $f: \mathbb{Z}_N \to \R^+$ be any non-negative function obeying the bound
\begin{equation}\label{f-bound-2}
0 \leq f(x) \leq \nu(x) \hbox{ for all } x \in \mathbb{Z}_N
\end{equation}
and
\begin{equation}\label{f-density}
 \E(f) \geq \delta.
\end{equation}
Then we have
\begin{equation}\label{recurrence}
 \E( f(x) f(x+r) \ldots f(x+(k-1)r) | x,r \in \mathbb{Z}_N ) \geq c(k,\delta) - o_{k,\delta}(1)
\end{equation}
where $c(k,\delta) > 0$ is the same constant which appears in Proposition \ref{sz}.  \textup{(}The decay rate $o_{k,\delta}(1)$, on the other hand, decays significantly slower than that in Proposition \ref{sz}, and depends of course on the decay rates in the linear forms and correlation conditions\textup{)}.
\end{theorem}

\ni We remark that while we do not explicitly assume that $N$ is large in Theorem \ref{main}, we are free to do so since the conclusion
\eqref{recurrence} is trivial when $N = O_{k,\delta}(1)$.  We certainly encourage the reader to think of $N$ has being extremely large compared
to other quantities such as $k$ and $\delta$, and to think of $o(1)$ errors as being negligible.\vs

\ni The proof of this theorem will occupy the next few sections, \S\ref{sec4}--\ref{sec7}.  Interestingly, the proof requires no Fourier analysis, additive combinatorics, or number theory; the argument is instead a blend of quantitative ergodic theory arguments with some combinatorial estimates related to Gowers uniformity and sparse hypergraph regularity.  From \S \ref{sec8} onwards we will apply this theorem to
the specific case of the primes, by establishing a pseudorandom majorant for
(a modified version of) the von Mangoldt function.  

\section{Notation}\label{sec4}

\ni We now begin the proof of Theorem \ref{main}. Thoughout this proof we fix the parameter $k \geq 3$
and the probability density $\nu$ appearing in 
Theorem \ref{main}.  All our constants in the $O()$ and $o()$ notation are allowed to depend on $k$ (with all future dependence on this parameter being suppressed), and are also allowed to depend on
the bounds implicit in the right-hand sides of \eqref{lfc} and \eqref{eq3.1}.  We may take $N$ to be 
sufficiently large with respect to $k$ and $\delta$ since \eqref{recurrence} is trivial otherwise.\vs

\ni We need some standard $L^q$ spaces.

\begin{definition}  For every $1 \leq q \leq \infty$ and $f: \mathbb{Z}_N \to \R$, we define the $L^q$ norms as
$$\|f\|_{L^q} := \E(|f|^q)^{1/q}$$
with the usual convention that $\|f\|_{L^\infty} := \sup_{x \in \mathbb{Z}_N} |f(x)|$.  We let $L^q(\Z_N)$ be the Banach space of all
functions from $\Z_N$ to $\R$ equipped with the $L^q$ norm; of course since $\Z_N$ is finite these spaces are all equal to each other
as vector spaces, but the norms are only equivalent up to powers of $N$.  We also observe that $L^2(\Z_N)$ is a real Hilbert space
with the usual inner product
$$ \langle f, g \rangle := \E( f g ).$$
If $\Omega$ is a subset of $\mathbb{Z}_N$, we use ${\bf 1}_\Omega: \mathbb{Z}_N \to \R$ to denote the indicator function of
$\Omega$, thus ${\bf 1}_\Omega(x) = 1$ if $x \in \Omega$ and ${\bf 1}_\Omega(x) = 0$ otherwise.  Similarly if $P(x)$ is a statement
concerning an element $x \in \Z_N$, we write ${\bf 1}_{P(x)}$ for ${\bf 1}_{\{ x \in \Z_N: P(x) \}}(x)$.
\end{definition}

\ni In our arguments we shall frequently be performing linear changes of variables and then taking expectations.  To facilitate this we 
adopt the following definition.
Suppose that $A$ and $B$ are finite non-empty sets and that $\Phi: A \to B$ is a map. Then we say that $\Phi$ is a \textit{uniform cover of $B$ by $A$} if $\Phi$ is surjective and all the fibers $\{ \Phi^{-1}(b): b \in B\}$ have the same cardinality (i.e. they have
cardinality $|A|/|B|$).  Observe that if $\Phi$ is a uniform cover of $B$ by $A$, then for any function $f: B \to \R$ 
we have 
\begin{equation}\label{uniform-cover}
\E( f(\Phi(a)) | a \in A ) = \E( f(b) | b \in B ).
\end{equation}

\section{Gowers uniformity norms, and a generalized von Neumann theorem}\label{sec5}

\ni As mentioned in earlier sections, the proof of Theorem \ref{main} relies on splitting the given function $f$ into
a Gowers uniform component and a Gowers anti-uniform component.  We will come to this splitting in later sections, but for this section
we focus on defining the notion of Gowers uniformity, introduced in \cite{gowers-4,gowers}.  The main result of
this section will be a generalized von Neumann theorem (Proposition \ref{vn}), which basically asserts that Gowers uniform functions
are negligible for the purposes of computing sums such as \eqref{recurrence}.

\begin{definition}  Let $d \geq 0$ be a dimension\footnote{In practice, we will have $d = k-1$, where $k$ is the 
length of the arithmetic progressions under consideration.}.  We let $\{0,1\}^d$ be the standard discrete $d$-dimensional
cube, consisting of $d$-tuples $\omega = (\omega_1,\ldots,\omega_d)$ where $\omega_j \in \{0,1\}$ for $j=1,\ldots,d$.
If $h = (h_1,\ldots,h_d)
\in \mathbb{Z}_N^d$ we define $\omega \cdot h := \omega_1 h_1 + \ldots + \omega_d h_d$.  If $(f_\omega)_{\omega \in \{0,1\}^d}$
is a $\{0,1\}^d$-tuple of functions in $L^\infty(\mathbb{Z}_N)$, we define the \emph{$d$-dimensional Gowers inner product}
$\langle (f_\omega)_{\omega \in \{0,1\}^d} \rangle_{U^d}$ by the formula
\begin{equation}\label{inner-def}
 \langle (f_\omega)_{\omega \in \{0,1\}^d} \rangle_{U^d}
:= \E\bigg( \prod_{\omega \in \{0,1\}^d} f_\omega(x + \omega \cdot h) \; \bigg| \; x \in \mathbb{Z}_N, h \in \mathbb{Z}_N^d \bigg).
\end{equation}
\end{definition}

\ni Henceforth we shall refer to a configuration $\{x + \omega \cdot h : \omega \in \{0,1\}^d\}$ as a \emph{cube of dimension $d$}.\vs

\ni \emph{Example.}  When $d=2$, we have
$$ \langle f_{00}, f_{10}, f_{01}, f_{11} \rangle_{U^d} =
\E( f_{00}(x) f_{10}(x+h_1) f_{01}(x+h_2) f_{11}(x+h_1+h_2) \; | \; x, h_1, h_2 \in \Z_N ).$$

\ni We recall from \cite{gowers} the positivity properties of the Gowers inner product \eqref{inner-def} when $d \geq 1$ (the $d=0$
case being trivial).  First suppose that $f_\omega$ does not depend on the final digit $\omega_d$ of $\omega$, thus
$f_\omega = f_{\omega_1, \ldots, \omega_{d-1}}$.  Then we may rewrite \eqref{inner-def} as
\begin{align*}
\langle (f_\omega)_{\omega \in \{0,1\}^d} \rangle_{U^d}
= \E\bigg(& \prod_{\omega' \in \{0,1\}^{d-1}}  
f_{\omega'}(x + \omega' \cdot h')
f_{\omega'}(x + h_d + \omega' \cdot h')\\
& \qquad\qquad\qquad\qquad\qquad \bigg| \; x \in \mathbb{Z}_N, h' \in \mathbb{Z}_N^{d-1}, h_d \in \mathbb{Z}_N \bigg),
\end{align*}
where we write $\omega' := (\omega_1, \ldots, \omega_{d-1})$ and $h' := (h_1, \ldots, h_{d-1})$.  This can be
rewritten further as
\begin{equation}\label{f-cube}
\langle (f_\omega)_{\omega \in \{0,1\}^d} \rangle_{U^d}
=
\E\bigg(
\big|\E\big( \prod_{\omega' \in \{0,1\}^{d-1}} f_{\omega'}(y + \omega' \cdot h') | y \in \mathbb{Z}_N\big)\big|^2
 \;\; \bigg| \; h' \in \mathbb{Z}_N^{d-1} \bigg),
\end{equation}
so in particular we have the positivity property
$\langle (f_\omega)_{\omega \in \{0,1\}^d} \rangle_{U^d} \geq 0$
when $f_\omega$ is independent of $\omega_d$.  This proves the positivity property
\begin{equation}\label{positivity}
 \langle (f)_{\omega \in \{0,1\}^d} \rangle_{U^d} \geq 0
\end{equation}
when $d \geq 1$.  We can thus define the \emph{Gowers uniformity norm} $\|f\|_{U^d}$ of a function 
$f: \mathbb{Z}_N \to \R$ by the formula
\begin{equation}\label{ud-def}
 \|f \|_{U^d} := \langle (f)_{\omega \in \{0,1\}^d} \rangle_{U^d}^{1/2^d}
= \E\bigg( \prod_{\omega \in \{0,1\}^d} f(x + \omega \cdot h) \bigg| x \in \mathbb{Z}_N, h \in \mathbb{Z}_N^d \bigg)^{1/2^d}.
\end{equation}

\ni When $f_\omega$ does depend on $\omega_d$, \eqref{f-cube} must be rewritten as
\begin{align*}
\langle (f_\omega)_{\omega \in \{0,1\}^d} \rangle_{U^d}
=
\E\bigg(&
\E\big( \prod_{\omega' \in \{0,1\}^{d-1}} f_{\omega',0}(y + \omega' \cdot h') \big| y \in \mathbb{Z}_N\big) \times \\
& \times \E\big( \prod_{\omega' \in \{0,1\}^{d-1}} f_{\omega',1}(y + \omega' \cdot h') \big| y \in \mathbb{Z}_N\big )
\bigg| h' \in \mathbb{Z}_N^{d-1} \bigg).
\end{align*}
From the Cauchy-Schwarz inequality in the $h'$ variables, we thus see that
$$ |\langle (f_\omega)_{\omega \in \{0,1\}^d} \rangle_{U^d}|
\leq \langle (f_{\omega',0})_{\omega \in \{0,1\}^d} \rangle_{U^d}^{1/2}
\langle (f_{\omega',1})_{\omega \in \{0,1\}^d} \rangle_{U^d}^{1/2}.$$
Similarly if we replace the role of the $\omega_d$ digit by any of the other digits.  Applying this Cauchy-Schwarz inequality once in each digit, we obtain the \emph{Gowers Cauchy-Schwarz inequality}
\begin{equation}\label{gcz}
 |\langle (f_\omega)_{\omega \in \{0,1\}^d} \rangle_{U^d}| \leq \prod_{\omega \in \{0,1\}^d} \|f_\omega \|_{U^d}.
\end{equation}
From the multilinearity of the inner product, and the binomial formula, we then obtain the inequality
$$ |\langle (f+g)_{\omega \in \{0,1\}^d} \rangle_{U^d}| \leq (\|f\|_{U^d} + \|g\|_{U^d})^{2^d}$$
whence we obtain the \emph{Gowers triangle inequality}
$$ \| f + g \|_{U^d} \leq \|f\|_{U^d} + \|g\|_{U^d}.$$
(cf. \cite{gowers} Lemmas 3.8 and 3.9).\vs

\ni\textit{Example.} Continuing the $d=2$ example, we have
$$ \|f\|_{U^2} := \E( f(x) f(x+h_1) f(x+h_2) f(x+h_1+h_2) \; | \; x,h_1,h_2 \in \Z_N )^{1/4}$$
and the Gowers Cauchy-Schwarz inequality then states
\begin{eqnarray*} & &  |\E( f_{00}(x) f_{10}(x+h_1) f_{01}(x+h_2) f_{11}(x+h_1+h_2) | x,h_1,h_2 \in \Z_N )| \\ & & \qquad\qquad\qquad\qquad\qquad \leq
\|f_{00} \|_{U^2} \|f_{10} \|_{U^2} \|f_{01} \|_{U^2} \|f_{11} \|_{U^2}.\end{eqnarray*}
Applying this with $f_{10}$, $f_{01}$, $f_{11}$ set equal to Kronecker delta functions, one can easily verify that
$$ f_{00} \equiv 0 \hbox{ whenever } \|f_{00} \|_{U^2} = 0.$$
This, combined with the preceding discussion, shows that the $U^2$ norm is indeed a genuine norm.  This can also be seen by
the easily verified identity
$$ \| f\|_{U^2} = \big(\sum_{\xi \in \mathbb{Z}_N} |\widehat{f}(\xi)|^4\big)^{1/4}$$
(cf. \cite{gowers}, Lemma 2.2), where the Fourier transform $\hat f: \Z_N \to \C$ of $f$ is defined by the formula\footnote{The Fourier
transform of course plays a hugely important r\^ ole in the $k=3$ theory, and provides some
very useful intuition to then think about the higher $k$ theory, but will not be used in this paper except as motivation.}
\[ \widehat{f}(\xi) := \E( f(x) e^{-2\pi i x \xi/N}| x \in \mathbb{Z}_N ).\]
for any $\xi \in \mathbb{Z}_N$. 

\ 

\ni We return to the study of general $U^d$ norms.  Since 
\begin{equation}\label{1-norm}
\|\vconst\|_{U^d} = \|1\|_{U^d} = 1, 
\end{equation}
we see from \eqref{gcz} that
$$  |\langle (f_\omega)_{\omega \in \{0,1\}^d} \rangle_{U^d}| \leq \|f\|_{U^d}^{2^{d-1}}$$
where $f_\omega := 1$ when $\omega_d = 1$ and $f_\omega := f$ when $\omega_d = 0$.  But the left-hand side can
easily be computed to be $\|f\|_{U^{d-1}}^{2^{d-1}}$, and thus we have the monotonicity relation
\begin{equation}\label{ud-monotone} \| f\|_{U^{d-1}} \leq \| f \|_{U^d}
\end{equation}
for all $d \geq 2$. Since the $U^2$ norm was already shown to be strictly positive, we see that the higher norms $U^d$, $d \geq 2$ are also.
Thus the $U^d$ norms are genuinely norms for all $d \geq 2$.  On the other hand, the $U^1$ norm is not actually a norm, since one can compute from \eqref{ud-def} that $\|f\|_{U^1} = |\E(f)|$ and thus $\|f\|_{U^1}$
may vanish without $f$ itself vanishing.  

\ 

\ni From the linear forms condition one can easily verify that $\| \nu\|_{U^d} = 1 + o(1)$ (cf. \eqref{example-1}).
In fact more is true, namely that pseudorandom measures $\nu$ are close to the
constant measure $\vconst$ in the $U^d$ norms; this is of course
consistent with our philosophy of deducing Theorem \ref{main} from
Theorem \ref{sz}.

\begin{lemma}\label{mua} Suppose that $\nu$ is $k$-pseudorandom \textup{(}as defined in Definition \ref{mu-pseudo-def}\textup{)}. Then we have
\begin{equation}\label{mu-uniform}
\| \nu - \vconstlem \|_{U^d} = \| \nu - 1 \|_{U^d} = o(1)
\end{equation}
for all $1 \leq d \leq k-1$.
\end{lemma}

\proof 
By \eqref{ud-monotone} it suffices to prove the claim for $d=k-1$.  Raising to the power $2^{k-1}$, it suffices
from \eqref{ud-def} to show that
\[ \E\bigg( \prod_{\omega \in \{0,1\}^{k-1}} (\nu(x + \omega \cdot h)-1) \;
\bigg| \; x \in \mathbb{Z}_N, h \in \mathbb{Z}_N^{k-1}\bigg) = o(1).\]
The left-hand side can be expanded as
\begin{equation}\label{eq5.1} \sum_{A \subseteq \{0,1\}^{k-1}} (-1)^{|A|}
\E\bigg( \prod_{\omega \in A} \nu(x + \omega \cdot h) \; \bigg| \; x \in \mathbb{Z}_N, h \in \mathbb{Z}_N^{k-1} \bigg).\end{equation}
Let us look at the expression
\begin{equation}\label{linear-forms-1}
 \E\bigg( \prod_{\omega \in A} \nu(x + \omega \cdot h) \; \bigg| \; x \in \mathbb{Z}_N, h \in \mathbb{Z}_N^{k-1} \bigg)
\end{equation}
for some fixed $A \subseteq \{0,1\}^{k-1}$.
This is of the form
\[ \E\left(\nu(\psi_1(\mathbf{x}))\dots \nu(\psi_{|A|}(\mathbf{x})) \; |\; \mathbf{x} \in \mathbb{Z}_N^k\right),\]
where $\mathbf{x} := (x,h_1,\dots,h_{k-1})$ and the $\psi_1, \ldots, \psi_{|A|}$ are some ordering of the $|A|$ linear forms $x + \omega \cdot h$, $\omega \in A$. It is clear that none of these forms is a rational multiple of any other. Thus we may invoke the $(2^{k-1},k,1)$-linear forms condition, which is a consequence of the fact that $\nu$ is $k$-pseudorandom, to conclude that the expression \eqref{linear-forms-1} is $1 + o(1)$.\vs

\ni Referring back to \eqref{eq5.1}, one sees that the claim now follows from the binomial theorem $\sum_{A \subseteq \{0,1\}^{k-1}} (-1)^{|A|} = (1-1)^{2^{k-1}} = 0$.
\endproof\vspace{11pt}

\ni It is now time to state and prove our ``generalised von Neumann theorem'', which explains how the expression \eqref{recurrence}, which counts $k$-term arithmetic progressions, is governed by the Gowers uniformity norms. All of this, of course, is relative to a pseudorandom measure $\nu$.

\begin{proposition}[Generalised von Neumann]\label{vn} Suppose that $\nu$ is $k$-pseudorandom. Let $f_0, \ldots, f_{k-1} \in L^1(\Z_N)$ be functions which
are pointwise bounded by $\nu+\vconstlem$, or in other words
\begin{equation}\label{fj-bounds}
 |f_j(x)| \leq \nu(x) + 1 \hbox{ for all } x \in \mathbb{Z}_N, 0 \leq j \leq k-1.
\end{equation}
Let $c_0, \ldots, c_{k-1}$ be a permutation of some $k$ consecutive elements of $\{-k+1,\dots,-1,0,1,\ldots,k-1\}$ \textup{(}in practice we will take $c_j := j$\textup{)}. Then 
\[ \E\bigg( \prod_{j=0}^{k-1} f_j(x + c_j r) \;  \bigg| \; x,r \in \mathbb{Z}_N \bigg) = O\big(\inf_{0 \leq j \leq k-1} \| f_j \|_{U^{k-1}}\big) + o(1).\]
\end{proposition}

\ni {\it Remark.}  This proposition is standard when $\nu = \vconst$ (see for instance \cite[Theorem 3.2]{gowers} or, for an analogous result in the ergodic setting, \cite[Theorem 3.1]{FKO}). The 
novelty is thus the extension to the pseudorandom $\nu$ studied in Theorem \ref{main}.  The reason we have an upper bound of $\nu(x)+1$ instead
of $\nu(x)$ is because we shall be applying this lemma to functions $f_j$ which roughly have the form $f_j = f - \E(f|\B)$, where $f$ is some function bounded pointwise by $\nu$, and $\B$ is a $\sigma$-algebra such that $\E(\nu|\B)$ is essentially bounded (up to $o(1)$ errors)
by $1$, so that we can essentially bound $|f_j|$ by $\nu(x)+1$; see Definition \ref{sigma-def} for the notations we are using here.  
The techniques here are inspired by similar Cauchy-Schwarz arguments relative to pseudorandom hypergraphs in \cite{gowers-reg}.
Indeed, the estimate here can be viewed as a kind of ``sparse counting lemma'' that utilises a regularity hypothesis (in the guise of $U^{k-1}$ control
on one of the $f_j$) to obtain control on an expression which can be viewed as a weighted count of arithmetic progressions concentrated in
a sparse set (the support of $\nu$).  See \cite{gowers-reg,klr} for some further examples of such lemmas.
\vs

\proof  By replacing $\nu$ with $(\nu+1)/2$ (and by dividing $f_j$ by 2), and using Lemma \ref{halfway}, we see that we may in fact assume
without loss of generality that we can improve \eqref{fj-bounds} to
\begin{equation}\label{fj-bounds-better}
 |f_j(x)| \leq \nu(x) \hbox{ for all } x \in \mathbb{Z}_N, 0 \leq j \leq k-1.
\end{equation}
For similar reasons we may assume that $\nu$ is strictly positive everywhere.\vs

\ni By permuting the $f_j$ and $c_j$ if necessary, we may assume that the infimum \[ \inf_{0 \leq j \leq k-1}
\| f_j \|_{U^{k-1}}\] is attained when $j=0$.  By shifting $x$ by $c_0 r$ if necessary we may assume that $c_0 = 0$.
Our task is thus to show
\begin{equation}\label{gvn-task}
 \E\bigg( \prod_{j=0}^{k-1} f_j(x + c_j r) \; \bigg| \; x,r \in \mathbb{Z}_N \bigg) = O\big(\| f_0 \|_{U^{k-1}}\big) + o(1).
\end{equation}
The proof of this will fall into two parts. First of all we will use the Cauchy-Schwarz inequality $k - 1$ times (as is standard in the proof of theorems of this general type). In this way we will bound the left hand side of \eqref{gvn-task} by a \textit{weighted} sum of $f_0$ over $(k-1)$-dimensional cubes. After that, we will show using the linear forms condition that these weights are roughly 1 on average, which will enable us to deduce \eqref{gvn-task}.\vs

\ni Before we give the full proof, let us first give the argument in the case
$k=3$, with $c_j := j$. This is conceptually no easier than the general case, but the notation is substantially less fearsome.
Our task is to show that  
$$\E\big( f_0(x) f_1(x+r) f_2(x+2r) \;  \big| \; x,r \in \mathbb{Z}_N \big) = O\big( \| f_0 \|_{U^2} \big) + o(1).$$
It shall be convenient to reparameterise the progression $(x,x+r,x+2r)$ as $(y_1+y_2, y_2/2, -y_1)$. The fact that the first term does not depend on $y_1$ and the second term does not depend on $y_2$ will allow us to perform
Cauchy-Schwarz in the arguments below without further changes of variable.  
Since $N$ is a large prime, we are now faced with estimating the quantity
\begin{equation}\label{j0-def}
J_0 := \E\big( f_0(y_1+y_2) f_1(y_2/2) f_2(-y_1) \;  \big| \; y_1,y_2 \in \mathbb{Z}_N \big).
\end{equation}
We estimate $f_2$ in absolute value by $\nu$ and bound this by
$$ |J_0| \leq \E\big( |\E(f_0(y_1+y_2) f_1(y_2/2)\; | \; y_2 \in \Z_N)| \nu(-y_1) \;  \big| \; y_1 \in \mathbb{Z}_N \big).$$
Using Cauchy-Schwarz and \eqref{numean}, we can bound this by
$$ (1 + o(1)) \E\big( |\E(f_0(y_1+y_2) f_1(y_2/2)\; | \; y_2 \in \Z_N)|^2 \nu(-y_1) \;  \big| \; y_1 \in \mathbb{Z}_N \big)^{1/2}$$
which we rewrite as $(1 + o(1)) J_1^{1/2}$, where
$$ J_1 := \E\big( f_0(y_1+y_2) f_0(y_1+y'_2) f_1(y_2/2) f_1(y'_2/2) \nu(-y_1) \;  \big| \; y_1,y_2,y'_2 \in \mathbb{Z}_N \big).$$
We now estimate $f_1$ in absolute value by $\nu$, and thus bound \[
J_1 \leq  \E\big( |\E( f_0(y_1+y_2) f_0(y_1+y'_2) \nu(-y_1) | y_1 \in \Z_N )| \nu(y_2/2) \nu(y'_2/2) \;  \big| \; y_2,y'_2 \in \mathbb{Z}_N \big).
\]
Using Cauchy-Schwarz and \eqref{numean} again, we bound this by $1 + o(1)$ times
$$ \E\big( |\E( f_0(y_1+y_2) f_0(y_1+y'_2) \nu(-y_1) | y_1 \in \Z_N )|^2 \nu(y_2/2) \nu(y'_2/2) \;  \big| \; y_2,y'_2 \in \mathbb{Z}_N \big)^{1/2}.$$
Putting all this together, we conclude the inequality
\begin{equation}\label{jojo}
 |J_0| \leq (1 + o(1)) J_2^{1/4},
\end{equation}
where
\begin{eqnarray*}  J_2  & := & \E\big( f_0(y_1+y_2) f_0(y_1+y'_2) f_0(y'_1+y_2) f_0(y'_1+y'_2) \nu(-y_1) \nu(-y'_1) \nu(y_2/2) \nu(y'_2/2) \; \\ & & \qquad\qquad\qquad\qquad \big| \; 
y_1,y'_1,y_2,y'_2 \in \mathbb{Z}_N \big).\end{eqnarray*}
If it were not for the weights involving $\nu$, $J_2$ would be the $U^2$ norm of $f_0$, and we would be done.  If we reparameterise
the cube $(y_1+y_2,y'_1+y_2,y_1+y'_2,y'_1+y'_2)$ by $(x,x+h_1,x+h_2,x+h_1+h_2)$, the above expression becomes
$$ J_2 = \E\big( f_0(x) f_0(x+h_1) f_0(x+h_2) f_0(x+h_1+h_2) W(x,h_1,h_2) \;  \big| \; 
x,h_1,h_2 \in \mathbb{Z}_N \big)$$
where $W(x,h_1,h_2)$ is the quantity
\begin{equation}\label{WEIGHT-def}
 W(x,h_1,h_2) := \E\big( \nu(-y) \nu(-y-h_1) \nu((x-y)/2) \nu((x-y-h_2)/2) \; | \; y \in \Z_N \big).
 \end{equation}
 In order to compare $J_2$ to $\Vert f_0 \Vert_{U^2}^4$, we must compare $W$ to $1$. To that end it suffices to show that the error
$$ \E\big( f_0(x) f_0(x+h_1) f_0(x+h_2) f_0(x+h_1+h_2) (W(x,h_1,h_2)-1) \;  \big| \; 
x,h_1,h_2 \in \mathbb{Z}_N \big)$$
is suitably small (in fact it will be $o(1)$).  To achieve this we estimate $f_0$ in absolute value by $\nu$ and use Cauchy-Schwarz one last time
to reduce to showing that
$$ \E\big( \nu(x) \nu(x+h_1) \nu(x+h_2) \nu(x+h_1+h_2) (W(x,h_1,h_2)-1)^n \;  \big| \; 
x,h_1,h_2 \in \mathbb{Z}_N \big) = 0^n + o(1)$$
for $n=0,2$. Expanding out the $W-1$ term, it suffices to show that
$$ \E\big( \nu(x) \nu(x+h_1) \nu(x+h_2) \nu(x+h_1+h_2) W(x,h_1,h_2)^q \;  \big| \; 
x,h_1,h_2 \in \mathbb{Z}_N \big) = 1+o(1)$$
for $q=0,1,2$.  But this follows from the linear forms condition (for instance, the case $q=2$ is just \eqref{example-3}).\vs

\ni We turn now to the proof of \eqref{gvn-task} in general. As one might expect in view of the above discussion, this shall consist of a large number of applications
of Cauchy-Schwarz to replace all the functions $f_j$ with $\nu$, and then applications of the linear forms condition.
In order to expedite these applications of Cauchy-Schwarz we shall need some notation.
Suppose that $0 \leq d \leq k-1$, and that we have
two vectors $y = (y_1,\ldots,y_{k-1}) \in \Z_N^{k-1}$ and
$y' = (y'_{k-d},\ldots,y'_{k-1}) \in \Z_N^d$ of length $k-1$ and $d$
respectively.  For any set $S \subseteq \{k-d,\ldots,k-1\}$, we define the vector $y^{(S)} = (y^{(S)}_1, \ldots, y^{(S)}_{k-1})
\in \Z_N^{k-1}$ as 
$$ y^{(S)}_i := \left\{ \begin{array}{ll}
y_i &\hbox{ if } i \not \in S \\
y'_i &\hbox{ if } i \in S.
\end{array}\right.$$
The set $S$ thus indicates which components of $y^{(S)}$ come from $y'$ rather
than $y$.

\begin{lemma}[Cauchy-Schwarz]\label{weighted-cubes}
Let $\nu : \mathbb{Z}_N \rightarrow \mathbb{R}^+$ be any measure. Let $\phi_0,\phi_1,\dots,\phi_{k-1} : \mathbb{Z}_N^{k-1} \rightarrow \mathbb{Z}_N$ be functions of $k-1$ variables $y_1,\ldots,y_{k-1}$, such that $\phi_i$ does not depend on $y_i$ for $1 \leq i \leq k-1$. Suppose that $f_0,f_1,\dots,f_{k-1} \in L^1(\Z_N)$ are functions satisfying $|f_i(x)| \leq \nu(x)$ for all $x \in \mathbb{Z}_N$ and for each $i$, $0 \leq i \leq k-1$. For each $0 \leq d \leq k-1$, define the quantities
\begin{equation}\label{j-def}
 J_{d} :=   \mathbb{E}\bigg(\prod_{S \subseteq \{k-d,\dots,k-1\}} \big(\prod_{i = 0}^{k-d-1} f_i(\phi_i(y^{(S)}))\big) \big(\prod_{i = k-d}^{k-1} \nu^{1/2}(\phi_i (y^{(S)}))\big)\; \bigg| \; y \in \mathbb{Z}_N^{k-1}, y' \in \Z_N^d\bigg)
\end{equation}
and
\begin{equation}\label{alt-p-exp} P_d := \mathbb{E}\bigg( \prod_{S \subseteq \{k-d,\dots,k-1\}} \nu(\phi_{k-d-1}(y^{(S)})) \; \bigg| \; y \in \Z_N^{k-1}, y' \in \Z_N^d\bigg).\end{equation}
Then for any $0 \leq d \leq k-2$, we have the inequality
\begin{equation}\label{to-prove-j}
|J_{d}|^2 \leq P_d J_{d+1}.
\end{equation}
\end{lemma}

\noindent\textit{Remarks.}
The appearance of $\nu^{1/2}$ in \eqref{j-def} may seem odd. Note, however, that since $\phi_i$ does not depend on the $i^{th}$ variable, each factor of $\nu^{1/2}$ in \eqref{j-def} occurs twice. If one takes $k=3$ and 
\begin{equation}\label{phi-example}
\phi_0(y_1,y_2) = y_1+y_2, \quad \phi_1(y_1,y_2) = y_2/2, \quad \phi_2(y_1,y_2) = -y_1,
\end{equation}
then the above notation is consistent with the quantities $J_0,J_1,J_2$ defined in the preceding discussion.
\vs 

\noindent\textit{Proof of Lemma \ref{weighted-cubes}.} 
Consider the quantity $J_d$. Since $\phi_{k-d-1}$ does not depend on $y_{k-d-1}$, we may take all quantities depending on $\phi_{k-d-1}$ outside of the
$y_{k-d-1}$ average.  This allows us to write
\[ J_d = \mathbb{E}\big( G(y,y')H(y,y')\; \big|\; y_1,\dots,y_{k-d-2},y_{k-d},\dots,y_{k-1},y'_{k-d},\dots,y'_{k-1} \in \mathbb{Z}_N\big),\]
where
\[G(y,y') := \prod_{S \subseteq \{k-d,\dots,k-1\}}f_{k-d-1}(\phi_{k - d - 1}(y^{(S)})) \nu^{-1/2}(\phi_{k - d - 1}(y^{(S)})) \]
 and
\[
H(y,y') := \mathbb{E}\bigg(\prod_{S \subseteq \{k-d,\dots,k-1\}} \prod_{i=0}^{k-d-2} f_i(\phi_i(y^{(S)})) \prod_{i=k - d - 1}^{k-1} \nu^{1/2}(\phi_i (y^{(S)})) \; \bigg| \; y_{k - d - 1} \in \mathbb{Z}_N\bigg) \]
(note we have multiplied and divided by several factors of the form
$\nu^{1/2}(\phi_{k-d-1}(y^{(S)})$).  Now apply Cauchy-Schwarz to give
\[ |J_d|^2 \leq \mathbb{E} \big(|G(y,y')|^2 \; \big| \; y_1,\dots,y_{k-d-2},y_{k-d},\dots,y_{k-1},y'_{k-d},\dots,y'_{k-1} \in \mathbb{Z}_N\big) \times \] \[ \qquad \qquad\times\mathbb{E} \big(|H(y,y')|^2 \; \big| \; y_1,\dots,y_{k-d-2},y_{k-d},\dots,y_{k-1},y'_{k-d},\dots,y'_{k-1} \in \mathbb{Z}_N\big).\]
Since $|f_{k-d-1}(x)| \leq \nu(x)$ for all $x$, one sees from \eqref{alt-p-exp}
that
\[\mathbb{E} \big(|G(y,y')|^2 \; \big| \; y_1,\dots,y_{k-d-2},y_{k-d},\dots,y_{k-1},y'_{k-d},\dots,y'_{k-1} \in \mathbb{Z}_N\big) \leq P_{d}\]
(note that the $y_{k-d-1}$ averaging in \eqref{alt-p-exp} is redundant
since $\phi_{k-d-1}$ does not depend on this variable). 
Moreover, by writing in the definition of $H(y,y')$ and expanding out the square, replacing the averaging variable $y_{k - d - 1}$ with the new variables $y_{k - d - 1}, y'_{k - d - 1}$, one sees from \eqref{j-def} that
\[  \mathbb{E} \big(|H(y,y')|^2 \; \big| \; y_1,\dots,y_{k-d-2},y_{k-d},\dots,y_{k-1},y'_{k-d},\dots,y'_{k-1} \in \mathbb{Z}_N\big) = J_{d+1}.\]
The claim follows.
\endproof\vs

\ni Applying the above lemma $k-1$ times, we obtain in particular that
\begin{equation}\label{ind-hyp} |J_0|^{2^{k-1}} \leq J_{k-1} \prod_{d = 0}^{k-2} P_d^{2^{k-2-d}} .\end{equation}
Observe from \eqref{j-def} that
\begin{equation}\label{e-def} J_0 = \mathbb{E}\bigg( \prod_{i = 0}^{k-1} f_i(\phi_i(y)) \; \bigg| \;  y \in \mathbb{Z}_N^{k-1}\bigg).\end{equation}


\ni\noindent\textit{Proof of Proposition \ref{vn}.}
We will apply \eqref{ind-hyp}, observing that \eqref{e-def} can be used to count configurations $(x,x+c_1r,\dots,x + c_{k-1}r)$ by making a judicious choice of the functions $\phi_i$.  For $y = (y_1,\ldots,y_{k-1})$, take 
\[ \phi_i(y) := \sum_{j = 1}^{k-1} \left(1 - \frac{c_i}{c_j}\right)y_j\] for $i = 0,\dots,k-1$. Then $\phi_0(y) = y_1 + \dots + y_{k-1}$, $\phi_i(y)$ does not depend on $y_i$ and, as one can easily check, for any $y$ we have $\phi_i(y) = x + c_i r$ where
\[ r = -\sum_{i=1}^{k-1} \frac{y_i}{c_i}.\]
Note that \eqref{phi-example} is simply the case $k=3$, $c_j = j$ of this more general construction.
Now the map $\Phi : \mathbb{Z}_N^{k-1} \rightarrow \mathbb{Z}_N^2$ defined by
\[ \Phi(y) := (y_1+\dots+y_{k-1},\frac{y_1}{c_1} + \frac{y_2}{c_2} + \dots + \frac{y_{k-1}}{c_{k-1}})\] is a uniform cover, and so 
\begin{equation}\label{eqn5.11}
\mathbb{E}\bigg(\prod_{j=0}^{k-1} f_j(x + c_jr) \; \bigg| \; x,r \in \mathbb{Z}_N\bigg) = \mathbb{E}\bigg(\prod_{i=0}^{k-1} f_i(\phi_i(y)) \; \bigg| \; y \in \mathbb{Z}_N^{k-1}\bigg) = J_0
\end{equation} 
thanks to \eqref{e-def} (this generalises \eqref{j0-def}).  On the other hand
we have $P_d = 1 + o(1)$ for each $0 \leq d \leq k - 2$, since
the $k$-pseudorandom hypothesis on $\nu$ implies the
$(2^d, k-1+d,k)$-linear forms condition. Applying \eqref{ind-hyp} 
we thus obtain
\begin{equation}\label{eqn5.12} J_0^{2^{k-1}} \leq  (1 + o(1))J_{k-1}\end{equation}
(this generalises \eqref{jojo}).
Fix $y \in \mathbb{Z}_N^{k-1}$. As $S$ ranges over all subsets of $\{1,\dots,k-1\}$, $\phi_0(y^{(S)})$ ranges over a $(k-1)$-dimensional cube $\{ x + \omega \cdot h : \omega \in \{0,1\}^{k-1}\}$ where $x = y_1 + \dots + y_{k-1}$ and $h_i = y'_i - y_i$, $i = 1,\dots,k-1$. Thus we may write
\begin{equation}\label{eq21} J_{k-1} = \E\bigg( W(x,h) 
\prod_{\omega \in \{0,1\}^{k-1}} 
f_0(x + \omega \cdot h)\; \bigg| \; 
x \in \mathbb{Z}_N, h \in \mathbb{Z}_N^{k-1} \bigg)\end{equation}
where the weight function $W(x,h)$ is given by
\begin{eqnarray*} W(x,h) & = & 
\E \bigg(\prod_{\omega \in \{0,1\}^{k-1}} \prod_{i=1}^{k-1} \nu^{1/2}(\phi_i (y + \omega h)) \;  \bigg| \; y_1,\dots,y_{k-2}\in \mathbb{Z}_N\bigg) \\ 
& = & \E \bigg( \prod_{i=1}^{k-1}\prod_{\substack{ \omega \in \{0,1\}^{k-1} \\ \omega_i = 0}} \nu(\phi_i(y+\omega h))\; 
\bigg| \; y_1,\dots,y_{k-2} \in \mathbb{Z}_N\bigg)\end{eqnarray*}
(this generalises \eqref{WEIGHT-def}).
Here, $\omega h \in \Z_N^{k-1}$ is the vector with components $(\omega h)_j := \omega_j h_j$ for $1 \leq j \leq k-1$,
and $y \in \Z_N^{k-1}$ is the vector with components $y_j$ for $1 \leq j \leq k-2$ and $y_{k-1} := x - y_1 - \ldots - y_{k-2}$.
Now by the definition of the $U^{k-1}$ norm we have
\[
\E\bigg( \prod_{\omega \in \{0,1\}^{k-1}} f_0(x + \omega \cdot h) \; \bigg| \; 
x \in \mathbb{Z}_N, h \in \mathbb{Z}_N^{k-1} \bigg) = \| f_0 \|_{U^{k-1}}^{2^{k-1}}.\]
To prove \eqref{gvn-task} it therefore suffices, by \eqref{eqn5.11}, \eqref{eqn5.12} and \eqref{eq21}, to prove that
\[ \E\bigg( (W(x,h)-1) \prod_{\omega \in \{0,1\}^{k-1}} f_0(x + \omega \cdot h) \; \bigg| \; 
x \in \mathbb{Z}_N, h \in \mathbb{Z}_N^{k-1} \bigg) = o(1).\]
Using \eqref{fj-bounds-better}, it suffices to show that
\[ \E\bigg( |W(x,h)-1| \prod_{\omega \in \{0,1\}^{k-1}} \nu(x + \omega \cdot h) \; \bigg| \; 
x \in \mathbb{Z}_N, h \in \mathbb{Z}_N^{k-1} \bigg) = o(1).\]
Thus by Cauchy-Schwarz it will be enough to prove

\begin{lemma}[$\nu$ covers its own cubes uniformly]\label{mu-cubes}
For $n=0,2$, we 
have \[\E\bigg( |W(x,h)-1|^n \prod_{\omega \in \{0,1\}^{k-1}} \nu(x + \omega \cdot h)\; \bigg| \; 
x \in \mathbb{Z}_N, h \in \mathbb{Z}_N^{k-1} \bigg) = 0^n + o(1).\]
\end{lemma}

\proof Expanding out the square, it then suffices to show that
\[  \E\bigg( W(x,h)^q \prod_{\omega \in \{0,1\}^{k-1}} \nu(x + \omega \cdot h) \; \bigg| \;
x \in \mathbb{Z}_N, h \in \mathbb{Z}_N^{k-1} \bigg) = 1 + o(1)\]
for $q=0,1,2$. This can be achieved by three applications of the linear forms condition, as follows:\\[11pt]
$q = 0$. Use the $(2^{k-1},k,1)$-linear forms property with variables $x,h_1,\dots,h_{k-1}$ and forms
\[ x + \omega \cdot h \qquad \omega \in \{0,1\}^{k-1}.\]
$q = 1$. Use the $(2^{k-2}(k+1), 2k - 2, k)$-linear forms property with variables $x$, $h_1,\dots,h_{k-1}$, $y_1,\dots,y_{k-2}$ and forms
\begin{eqnarray*} \phi_i(y + \omega h) & & \omega \in \{0,1\}^{k-1}, \;\; \omega_i = 0, 1 \leq i \leq k - 1;\\
 x + \omega \cdot h, && \omega \in \{0,1\}^{k-1}.\end{eqnarray*}
$q = 2$. Use the $(k\cdot 2^{k-1}, 3k - 4, k)$-linear forms property with variables $x$, $h_1,\dots,h_{k-1}$, $y_1,\dots,y_{k-2}$, $y'_1,\dots,y'_{k-2}$ and forms
\begin{eqnarray*} \phi_i(y + \omega h) && \omega \in \{0,1\}^{k-1},\;\; \omega_i = 0, 1 \leq i \leq k - 1;\\
 \phi_i(y' + \omega h) && \omega \in \{0,1\}^{k-1},\;\; \omega_i = 0, 1 \leq i \leq k - 1;\\
 x + \omega \cdot h, && \omega \in \{0,1\}^{k-1}.\end{eqnarray*}
Here of course we adopt the convention that 
$y_{k-1} = x - y_1 - \ldots - y_{k-2}$ and
$y'_{k-1} = x - y'_1 - \ldots - y'_{k-2}$.
This completes the proof of the lemma, and hence of Proposition \ref{vn}.\endproof

\section{Gowers anti-uniformity}\label{sec6}

\ni Having studied the $U^{k-1}$ norm, we now introduce the dual $(U^{k-1})^*$ norm, defined in the usual
manner as
\begin{equation}\label{dual-def}
 \| g \|_{(U^{k-1})^*} := \sup\{ |\langle f, g \rangle|: f \in U^{k-1}(\mathbb{Z}_N), \|f\|_{U^{k-1}} \leq 1 \}.
\end{equation}
We say that $g$ is \emph{Gowers anti-uniform} if $\|g\|_{(U^{k-1})^*} = O(1)$ and $\|g\|_{L^\infty} = O(1)$. If $g$ is Gowers anti-uniform, and if $|\langle f,g \rangle|$ is large, then $f$ cannot be Gowers uniform (have small Gowers norm) since
\[ |\langle f,g \rangle| \leq \Vert f \Vert_{U^{k-1}}\Vert g \Vert_{(U^{k-1})^*}.\]
Thus Gowers anti-uniform functions can be thought of as ``obstructions to Gowers uniformity''.
\ni The $(U^{k-1})^*$ are well-defined norms for $k \geq 3$ since $U^{k-1}$ is then a genuine norm (not just a seminorm).
In this section we show how to generate a large class of Gowers anti-uniform functions, in order that we can decompose an arbitrary
function $f$ into a Gowers uniform part and a bounded Gowers anti-uniform part in the next section.\vs

\ni\textit{Remark.} In the $k=3$ case we have the explicit formula \begin{equation}\label{star} \|g\|_{(U^2)^*} = \big(\sum_{\xi \in \Z_N} |\widehat{g}(\xi)|^{4/3}\big)^{3/4} = \Vert \widehat{g} \Vert_{4/3}.\end{equation} We will not, however, require this fact except for motivational purposes.\vs

\ni A basic way to generate Gowers anti-uniform functions is the following.  For each function $F \in L^1(\Z_N)$, define
the \emph{dual function} ${\mathcal D} F$ of $F$ by
\begin{equation}\label{cald}
 {\mathcal D} F(x) := 
\E\bigg( \prod_{\omega \in \{0,1\}^{k-1}: \omega \neq 0^{k-1}} F(x + \omega \cdot h) \; \bigg| \; h \in \mathbb{Z}_N^{k-1} \bigg)
\end{equation}
where $0^{k-1}$ denotes the element of $\{0,1\}^{k-1}$ consisting entirely of zeroes.\vs

\ni \emph{Remark.} Such functions have arisen recently in work of Host and Kra \cite{host-kra2} in the ergodic theory setting (see also \cite{assani}).\vs

\ni The next lemma, while simple, is fundamental to our entire approach;
it asserts that if a function majorised by a pseudorandom measure $\nu$ is not Gowers uniform, then it correlates\footnote{This idea was inspired by the proof of the Furstenberg structure theorem \cite{furst,FKO}; a key point in that proof being that 
if a system is not (relatively) weakly mixing, then it must contain a non-trivial (relatively) almost periodic function, which can then be projected out via conditional expectation.  A similar idea also occurs in the proof of the Szemer\'edi regularity lemma \cite{szemeredi}.} with a bounded Gowers anti-uniform function. Boundedness is the key feature here. The idea in proving Theorem \ref{main} will then be to project out the influence of these bounded Gowers anti-uniform functions (through the machinery of conditional expectation) until one is only left with a Gowers uniform remainder, which can be discarded by the generalised von Neumann theorem (Proposition \ref{vn}). 

\begin{lemma}[Lack of Gowers uniformity implies correlation]\label{antiuniform}  Let $\nu$ be a $k$-pseudorandom measure, and let $F \in L^1(\Z_N)$ be any function.  Then we have the identities
\begin{equation}\label{dual-first}
 \langle F, {\mathcal D} F \rangle = \| F \|_{U^{k-1}}^{2^{k-1}}
\end{equation}
and
\begin{equation}\label{dual}
 \| {\mathcal D} F \|_{(U^{k-1})^*} = \| F \|_{U^{k-1}}^{2^{k-1}-1}.
\end{equation}
If furthermore we assume the bounds
\[ |F(x)| \leq \nu(x) + 1 \hbox{ for all } x \in \mathbb{Z}_N\]
then we have the estimate
\begin{equation}\label{dual-bounded}
 \| {\mathcal D} F \|_{L^\infty} \leq 2^{2^{k-1}-1} + o(1).
\end{equation}
\end{lemma}

\proof
The identity \eqref{dual-first} is clear just by expanding out both sides using \eqref{cald}, \eqref{ud-def}.  
To prove \eqref{dual} we may of course assume $F$ is not identically zero.  By
\eqref{dual-def} and \eqref{dual-first} it suffices to show that
$$ |\langle f, {\mathcal D} F \rangle| \leq \|f\|_{U^{k-1}} \| F \|_{U^{k-1}}^{2^{k-1}-1}$$
for arbitrary functions $f$.  But by \eqref{cald} the left-hand side is simply the Gowers inner product
$\langle (f_\omega)_{\omega \in \{0,1\}^{k-1}} \rangle_{U^{k-1}}$, where $f_\omega := f$ when $\omega = 0$
and $f_\omega := F$ otherwise.  The claim then follows from the Gowers Cauchy-Schwarz inequality \eqref{gcz}.\vs

\ni Finally, we observe that \eqref{dual-bounded} is a consequence of the linear forms condition.  Bounding $F$ by $2 (\nu+1)/2 = 2 \nu_{1/2}$, it suffices to show
that
$$ {\mathcal D} \nu_{1/2}(x) \leq 1 + o(1) $$
uniformly in the choice of $x \in \mathbb{Z}_N$.  The left-hand side can be expanded using \eqref{cald} as
$$ \E\bigg( \prod_{\omega \in \{0,1\}^{k-1}: \omega \neq 0^{k-1}} \nu_{1/2}(x + \omega \cdot h) \; \bigg| \;  h \in \mathbb{Z}_N^{k-1} \bigg).$$
By the linear forms condition \eqref{lfc} (and Lemma \ref{halfway}) this expression is $1 + o(1)$ (this is the only place in the paper that we appeal to the linear forms condition in the non-homogeneous case where some $b_i \neq 0$; here, all the $b_i$ are equal to $x$).  Note that
\eqref{example-2} corresponds to the $k=3$ case of this application of the linear forms condition.\endproof\vspace{11pt}

\noindent\textit{Remarks.} Observe that
if $P: \mathbb{Z}_N \to \mathbb{Z}_N$ is any polynomial on $\mathbb{Z}_N$ of degree
at most $k-2$, and $F(x) = e^{2\pi i P(x)/N}$, then\footnote{To make this assertion precise, one has to generalise the notion of dual function
to complex-valued functions by inserting an alternating sequence of conjugation signs; see \cite{gowers}.}
 ${\mathcal D} F = F$; this is basically
a reflection of the fact that taking $k-1$ successive differences of $P$ yields the zero function. Hence
by the above lemma $\|F\|_{(U^{k-1})^*} \leq 1$, and thus $F$ is Gowers anti-uniform.  One should keep these 
``polynomially quasiperiodic'' functions $e^{2\pi i P(x)/N}$ in mind as model examples of functions 
of the form ${\mathcal D} F$, whilst bearing in mind that they are not the only examples\footnote{The situation again has an intriguing parallel with ergodic theory, in which the r\^ ole of the Gowers anti-uniform functions of order $k-2$ appear to be played by $k-2$-step nilfactors (see \cite{host-kra2,host-kra3,ziegler}), which may contain polynomial eigenfunctions of order $k-2$, but can also exhibit slightly more general behaviour; see \cite{furst-weiss} for further discussion.}. For some further discussion on the role of such polynomials of degree $k-2$ in determining Gowers uniformity especially in the $k=4$ case, 
see \cite{gowers-4,gowers}.
Very roughly speaking, Gowers uniform
functions are analogous to the notion of ``weakly mixing'' functions that appear in ergodic theory proofs of
Szemer\'edi's theorem, whereas Gowers anti-uniform functions are somewhat analogous to the notion of 
``almost periodic'' functions. When $k =3$ there is a more precise relation with linear exponentials (which are the same thing as characters on $\mathbb{Z}_N$). When $\nu = 1$, for example, one has the explicit formula
\begin{equation}\label{dagger} {\mathcal D} F(x) = \sum_{\xi \in \mathbb{Z}_N} |\widehat{F}(\xi)|^2\widehat{F}(\xi) e^{2\pi i x \xi/N}.\end{equation}
Suppose for sake of argument that $F$ is bounded pointwise in magnitude by 1.
By splitting the set of frequencies $\mathbb{Z}_N$ into the sets $S := \{\xi : |\widehat{F}(\xi)| \geq \epsilon\}$ and $\Z_N \backslash S$ one sees that it is possible to write 
\[ \mathcal{D}F(x) = \sum_{\xi \in S} a_{\xi} e^{2\pi i x\xi/N} + E(x),\]
where $|a_{\xi}| \leq 1$ and $\Vert E \Vert_{L^\infty} \leq \epsilon$. Also, we have $|S| \leq \epsilon^{-2}$. Thus $\mathcal{D}F$ is equal to a linear combination of a few characters plus a small error.\vs

\ni Once again, these remarks concerning the relation with harmonic analysis are included only for motivational purposes.\vs

\ni Let us refer to a functions of the form ${\mathcal D} F$, where $F$ is pointwise bounded by $\nu+1$, as a \emph{basic Gowers anti-uniform function}.
Observe from \eqref{dual-bounded} that if $N$ is sufficiently large, then all basic Gowers anti-uniform functions take values
in the interval $I := [-2^{2^{k-1}}, 2^{2^{k-1}}]$.\vs

\ni The following is a statement to the effect that the measure $\nu$ is uniformly distributed with respect not just to each basic Gowers anti-uniform function
(which is a special case of \eqref{dual}), but also to the \emph{algebra} generated by such functions.

\begin{proposition}[Uniform distribution wrt basic Gowers anti-uniform functions]\label{weier} Suppose that $\nu$ is $k$-pseudorandom. Let $K \geq 1$ be a fixed integer, let $\Phi: I^K \to \R$ be a fixed continuous function, let $\D F_1, \ldots, \D F_K$ be basic Gowers anti-uniform functions,
and define the function $\psi: \Z_N \to \R$ by
\[ \psi(x) := \Phi({\mathcal D} F_1(x), \ldots, {\mathcal D} F_K(x)).\]
Then we have the estimate
$$ \langle \nu - 1, \psi \rangle = o_{K,\Phi}(1).$$
Furthermore if $\Phi$ ranges over a compact set $E \subset C^0(I^K)$ of the space $C^0(I^K)$ of continuous functions on $I^K$ \textup{(}in the uniform topology\textup{)} then the bounds here are uniform in $\Phi$ \textup{(}i.e. one can replace $o_{K,\Phi}(1)$ with
$o_{K,E}(1)$ in this case\textup{)}.
\end{proposition}

\noindent\textit{Remark.} In light of the previous remarks, we see in particular that $\nu$ is uniformly distributed with respect to any continuous
function of polynomial phase functions such as $e^{2\pi i P(x)/N}$, where $P$ has degree at most $k-2$.  \vs

\proof We will prove this result in two stages, first establishing the result for $\Phi$ polynomial and then using a Weierstrass approximation argument to deduce the general case.  Fix $K \geq 1$, and let $F_1, \ldots, F_K \in L^1(\Z_N)$ be fixed functions obeying the bounds 
$$ F_j(x) \leq \nu(x) + 1 \hbox{ for all } x \in \Z_N, 1 \leq j \leq K.$$
By replacing
$\nu$ by $(\nu+1)/2$, dividing the $F_j$ by two, and using Lemma \ref{halfway} as before, we may strengthen this bound without loss of generality to
\begin{equation}\label{fjb-better}
 |F_j(x)| \leq \nu(x) \hbox{ for all } x \in \mathbb{Z}_N, 1 \leq j \leq K.
\end{equation}

\begin{lemma}\label{poly-lemma} Let $d \geq 1$.  For any polynomial $P$ of $K$ variables and degree $d$ with real coefficients \textup{(}independent of $N$\textup{)}, we have
\[ \| P({\mathcal D} F_1, \ldots, {\mathcal D} F_K) \|_{(U^{k-1})^*} = O_{K,d,P}(1).\]
\end{lemma}

\noindent\textit{Remark.} It may seem surprising that that there is no size restriction on $K$ or $d$, since we are presumably going to use the linear forms or correlation conditions and we are only assuming those conditions with bounded parameters.  However whilst we do indeed restrict the size of $m$ in \eqref{eq3.2}, we do not need to restrict the size of $q$ in \eqref{eq3.1}.\vs

\proof  By linearity it suffices to prove this when $P$ is a monomial.
By enlarging $K$ to at most $dK$ and repeating the functions $F_j$ as necessary, it in fact suffices to prove this for the monomial $P(x_1,\ldots,x_K) = x_1 \ldots x_K$.
Recalling the definition of $(U^{k-1})^*$, we are thus required to show that
\[ \big\langle f, \prod_{j=1}^K {\mathcal D} F_j \big\rangle = O_K(1) \]
for all $f: \mathbb{Z}_N \to \R$ satisfying $\| f \|_{U^{k-1}} \leq 1$.  
By \eqref{cald} the left-hand side can be expanded as
\[
\E\bigg( f(x)
\prod_{j=1}^K \E\big( \prod_{\omega \in \{0,1\}^{k-1}: \omega \neq 0^{k-1}} F_j(x + \omega \cdot h^{(j)}) \;  \big| \;  h^{(j)} \in \mathbb{Z}_N^{k-1} \big) \; \bigg|\; x \in \mathbb{Z}_N\bigg).\]
We can make the change of variables $h^{(j)} = h + H^{(j)}$ for any $h \in \mathbb{Z}_N^{k-1}$, and then average over $h$,
to rewrite this as
\[
\E\bigg( f(x)
\prod_{j=1}^K \E\big( \prod_{\omega \in \{0,1\}^{k-1}: \omega \neq 0^{k-1}} F_j(x + \omega \cdot H^{(j)} + \omega \cdot h) \; \big| \; H^{(j)} \in \mathbb{Z}_N^{k-1} \big) \; \bigg| \; x \in \mathbb{Z}_N; h \in \mathbb{Z}_N^{k-1}\bigg).\]
Expanding the $j$ product and interchanging the expectations, we can rewrite this in terms of the Gowers inner product as
\[ \E\big( \langle (f_{\omega, H})_{\omega \in \{0,1\}^{k-1}} \rangle_{U^{k-1}} \; \big| \; H \in (\mathbb{Z}_N^{k-1})^K \big)\]
where $H := (H^{(1)}, \ldots, H^{(K)})$, $f_{0,H} := f$, and $f_{\omega,H} := g_{\omega \cdot H}$ for $\omega \neq 0^{k-1}$,
where $\omega \cdot H := (\omega \cdot H^{(1)}, \ldots, \omega \cdot H^{(K)})$ and
\begin{equation}\label{foh-def}
 g_{u^{(1)}, \ldots, u^{(K)}}(x) := \prod_{j=1}^K F_j(x + u^{(j)}) \hbox{ for all } u^{(1)}, \ldots, u^{(K)} \in \mathbb{Z}_N.
\end{equation}
By the Gowers-Cauchy-Schwarz inequality \eqref{gcz} we can bound this as
\[ \E\bigg( \| f\|_{U^{k-1}} \prod_{\omega \in \{0,1\}^{k-1}: \omega \neq 0^{k-1}} \| g_{\omega \cdot H} \|_{U^{k-1}}
\; \bigg| \; H \in (\mathbb{Z}_N^{k-1})^K \bigg)\]
so to prove the claim it will suffice to show that
\[ \E\bigg( \prod_{\omega \in \{0,1\}^{k-1}: \omega \neq 0^{k-1}} \| g_{\omega \cdot H} \|_{U^{k-1}} \; \bigg| \; H \in (\mathbb{Z}_N^{k-1})^K \bigg) = O_K(1).\]
By H\"older's inequality it will suffice to show that
$$\E\big( \| g_{\omega \cdot H} \|_{U^{k-1}}^{2^{k-1}-1} \; \big| \; H \in (\mathbb{Z}_N^{k-1})^K \big) = O_K(1)$$
for each $\omega \in \{0,1\}^{k-1} \backslash 0^{k-1}$.  \vs

\ni Fix $\omega$.  Since $2^{k-1}-1 \leq 2^{k-1}$, another application of H\"older's inequality shows that it in fact suffices
to show that
$$\E\big( \| g_{\omega \cdot H} \|_{U^{k-1}}^{2^{k-1}} \; \big| \; H \in (\mathbb{Z}_N^{k-1})^K \big) = O_K(1).$$
Since $\omega \neq 0^{k-1}$, the map $H \mapsto \omega \cdot H$ is a uniform covering of $\mathbb{Z}_N^K$ 
by $(\mathbb{Z}_N^{k-1})^K$.  Thus by \eqref{uniform-cover} we can rewrite the left-hand side as
$$\E\big( \| g_{u^{(1)}, \ldots, u^{(K)}} \|_{U^{k-1}}^{2^{k-1}} \; \big| \; u^{(1)}, \ldots, u^{(K)} \in \mathbb{Z}_N \big).$$
Expanding this out using \eqref{ud-def} and \eqref{foh-def}, we can rewrite the left-hand side as
$$\E\bigg( \prod_{\tilde \omega \in \{0,1\}^{k-1}} \prod_{j=1}^K
F_j(x + u^{(j)} + h \cdot \tilde \omega) \; \bigg| \; x \in \mathbb{Z}_N, h \in \mathbb{Z}_N^{k-1}, 
u^{(1)}, \ldots, u^{(K)} \in \mathbb{Z}_N \bigg).$$
This factorises as
$$ \E\bigg( \prod_{j=1}^K \E\big( \prod_{\tilde \omega \in \{0,1\}^{k-1}} 
F_j(x + u^{(j)} + h \cdot \tilde \omega) \; \big| \; u^{(j)} \in \mathbb{Z}_N \big) \; \bigg| \;
x \in \mathbb{Z}_N, h \in \mathbb{Z}_N^{k-1} \bigg).$$
Applying \eqref{fjb-better}, we reduce to showing that
$$ \E\bigg( \E\big( \prod_{\tilde \omega \in \{0,1\}^{k-1}} 
\nu(x + u + h \cdot \tilde \omega)\; \big| \; u \in \mathbb{Z}_N \big)^K \; \bigg| \;
x \in \mathbb{Z}_N, h \in \mathbb{Z}_N^{k-1} \bigg) = O_K(1).$$
We can make the change of variables $y := x+u$, and then discard the redundant $x$ averaging, to reduce to showing that
$$ \E\bigg( \E\big( \prod_{\tilde \omega \in \{0,1\}^{k-1}} 
\nu(y + h \cdot \tilde \omega) \; \big| \; y \in \mathbb{Z}_N \big)^K \; \bigg| \;
h \in \mathbb{Z}_N^{k-1} \bigg) = O_K(1).$$
Now we are ready to apply the correlation condition (Definition \ref{correlation-condition}). This is, in fact, the only time we will use that condition. It gives
\[ \E\bigg( \prod_{\tilde \omega \in \{0,1\}^{k-1}} \nu(y + h \cdot \tilde \omega) \; \bigg| \; y \in \mathbb{Z}_N\bigg) \leq \sum_{\tilde \omega, \tilde \omega' \in \{0,1\}^{k-1}: \tilde \omega \neq \tilde \omega'}
\tau(h \cdot (\tilde \omega - \tilde \omega'))\] where, recall, $\tau$ is a weight function satisfying $\mathbb{E}(\tau^q) = O_q(1)$ for all $q$. Applying the triangle inequality
in $L^K(\Z_N^{k-1})$, it thus suffices to show that
\[ \E\big( \tau(h \cdot (\tilde \omega - \tilde \omega'))^K \; \big| \; h \in \mathbb{Z}_N^{k-1} \big) = O_K(1)\]
for all distinct $\tilde \omega, \tilde \omega' \in \{0,1\}^{k-1}$.  But the map $h \mapsto h \cdot (\tilde \omega - \tilde \omega')$
is a uniform covering of $\mathbb{Z}_N$ by $(\mathbb{Z}_N)^{k-1}$, so by \eqref{uniform-cover} the left-hand side is
just $\E(\tau^K) $, which is $O_K(1)$.
\endproof\vspace{11pt}

\noindent\textit{Proof of Proposition \ref{weier}}.
Let $\Phi$, $\psi$ be as in the Proposition, and let $\eps > 0$ be arbitrary.  From \eqref{dual-bounded} we know that the basic
Gowers anti-uniform functions ${\mathcal D} F_1, \ldots, {\mathcal D} F_K$ take values in the compact interval
$I := [-2^{2^{k-1}}, 2^{2^{k-1}}]$ introduced earlier.  By the
Weierstrass approximation theorem, we can thus find a polynomial $P$ (depending only on $K$ and $\eps$) 
such that
$$ \| \Phi({\mathcal D} F_1, \ldots, {\mathcal D} F_K) - P({\mathcal D} F_1, \ldots, {\mathcal D} F_K) \|_{L^\infty} \leq \eps$$
and thus by \eqref{numean} and taking absolute values inside the inner product, we have
$$ |\langle \nu - 1, \Phi({\mathcal D} F_1, \ldots, {\mathcal D} F_K) - P({\mathcal D} F_1, \ldots, {\mathcal D} F_K)
\rangle| \leq (2 + o(1)) \eps.$$
On the other hand, from Lemma \ref{poly-lemma}, Lemma \ref{mua} and \eqref{dual-def} we have
$$ \langle \nu - 1, P({\mathcal D} F_1, \ldots, {\mathcal D} F_K) \rangle = o_{K,\eps}(1)$$
since $P$ depends on $K$ and $\eps$.
Combining the two estimates we thus see that
for $N$ sufficiently large (depending on $K$ and $\eps$) we have
$$ |\langle \nu - 1, \Phi({\mathcal D} F_1, \ldots, {\mathcal D} F_K) \rangle| \leq 4\eps$$
(for instance).
Since $\eps > 0$ was arbitrary, the claim follows.  It is clear that this argument also gives the uniform bounds when
$\Phi$ ranges over a compact set (by covering this compact set by finitely many balls of radius $\eps$ in the uniform
topology).
\endproof\vspace{11pt}

\ni\textit{Remarks.} The philosophy behind Proposition \ref{weier} and its proof is that $(U^{k-1})^{*}$ respects, to some degree, the algebra structure of the space of functions on $\Z_N$. However, $\Vert \cdot \Vert_{(U^{k-1})^*}$ is not itself an algebra norm even in the model case 
$k=3$, as can be seen from \eqref{star} by recalling that $\Vert g \Vert_{(U^2)^*} = \Vert \widehat{g} \Vert_{4/3}$.  Note, however, that
$\Vert g \Vert_{(U^2)^*} \leq \Vert g \Vert_A$, where $\Vert g \Vert_A:= \Vert \widehat{g} \Vert_1$ is the Wiener norm.  From Young's inequality
we see that the Wiener norm \textit{is} an algebra norm, that is to say $\Vert gh \Vert_A \leq \Vert g \Vert_A \Vert h \Vert_A$. Thus while the $(U^2)^*$ norm is not an algebra norm, it is at least majorised by an algebra norm.\vs

\ni Now \eqref{dagger} easily implies that if $0 \leq F(x) \leq \vconst(x)$ then $\Vert \widehat{\mathcal{D}F} \Vert_1 \leq 1$, and so in this case we really have identified an algebra norm (the Weiner norm) such that if $\Vert f \Vert_{U^2}$ is large then $f$ correlates with a bounded function with small algebra norm. The $(U^{k-1})^*$ norms can thus be thought of as combinatorial variants of the Wiener algebra norm which apply
to more general values of $k$ than the Fourier case $k=3$.  (See also \cite{tao:ergodic} for a slightly different generalization of the Wiener
algebra to the case $k > 3$.)\vs

\ni For the majorant $\nu$ that we will use to majorise the primes, it is quite likely that $\Vert \mathcal{D}F \Vert_A = O(1)$ whenever $0 \leq F(x) \leq \nu(x)$, which would allow us to use the Wiener algebra $A$ in place of $(U^2)^*$ in the $k=3$ case of the arguments here. To obtain this estimate, however, requires some serious harmonic analysis related to the restriction phenomenon (the paper \cite{green-tao} may be consulted for further information). Such a property does not seem to follow simply from the pseudorandomness of $\nu$, and generalisation to $U^{k-1}$, $k > 3$, seems very difficult (it is not even clear what the form of such a generalisation would be). \vs

\ni For these reasons, our proof of Proposition \ref{weier} does not mention any algebra norms explicitly.

\section{Generalised Bohr sets and $\sigma$-algebras}

\ni To use Proposition \ref{weier}, we shall associate a $\sigma$-algebra to each basic Gowers anti-uniform function, such that the measurable
functions in each such algebra can be approximated by a function of the type considered in Proposition \ref{weier}.  We begin by setting out
our notation for $\sigma$-algebras.

\begin{definition}\label{sigma-def}  A \emph{$\sigma$-algebra} $\B$ in $\mathbb{Z}_N$ is any collection
of subsets of $\mathbb{Z}_N$ which contains the empty set $\emptyset$ and the full set $\Z_N$, and is closed under
complementation, unions and intersections.  As $\mathbb{Z}_N$ is a finite set, we will not need to distinguish between countable and uncountable unions or intersections.  We define the \emph{atoms} of a $\sigma$-algebra to be the minimal non-empty elements of $\B$ (with respect to set inclusion); it is clear that the atoms in $\B$ form a partition of $\mathbb{Z}_N$, and $\B$ consists precisely of arbitrary unions of its atoms (including the empty union $\emptyset$).
A function $f \in L^q(\mathbb{Z}_N)$ is said to be \emph{measurable} with respect to a $\sigma$-algebra $\B$ if all the level 
sets $\{f^{-1}(\{x\}): x \in \R \}$ of $f$ lie in $\B$, or equivalently if $f$ is constant on each of the atoms of $\B$.  \vs

\ni We define $L^q(\B) \subseteq L^q(\mathbb{Z}_N)$ to be the subspace of $L^q(\mathbb{Z}_N)$ consisting of $\B$-measurable functions,
equipped with the same $L^q$ norm.  We can then define the conditional expectation operator $f \mapsto \E(f|\B)$ to
be the orthogonal projection of $L^2(\mathbb{Z}_N)$ to $L^2(\B)$; this is of course also defined on all the other $L^q(\mathbb{Z}_N)$
spaces since they are all the same vector space.  An equivalent definition of conditional expectation is
$$ \E(f|\B)(x) := \E( f(y) | y \in \B(x) )$$
for all $x \in \mathbb{Z}_N$, where $\B(x)$ is the unique atom in $\B$ which contains $x$.  It is clear that conditional
expectation is a linear self-adjoint orthogonal projection on $L^2(\mathbb{Z}_N)$, is a contraction on $L^q(\mathbb{Z}_N)$ for
every $1 \leq q \leq \infty$, preserves non-negativity, and also preserves constant functions.  Also, if $\B'$ is a
subalgebra of $\B$ then $\E(\E(f|\B) | \B') = \E(f|\B')$.\vs

\ni If $\B_1, \ldots, \B_K$ are $\sigma$-algebras, we use $\bigvee_{j=1}^K \B_j = \B_1 \vee \ldots \vee \B_K$ to denote the $\sigma$-algebra generated by these algebras, or in other words the algebra whose atoms are the intersections of atoms in $\B_1, \ldots, \B_K$.  We adopt the usual convention that when $K=0$, the join $\bigvee_{j=1}^K \B_j$ is just the trivial $\sigma$-algebra $\{ \emptyset, \Z_N\}$.
\end{definition}

\ni We now construct the basic $\sigma$-algebras that we shall use.  We view the basic Gowers anti-uniform functions as generalizations of complex exponentials,
and the atoms of the $\sigma$-algebras we use can be thought of as ``generalised Bohr sets''.

\begin{proposition}[Each function generates a $\sigma$-algebra]\label{sigma-gen}  Let $\nu$ be a $k$-pseudorandom measure, let $0 < \eps < 1$ and $0 < \eta < 1/2$ be parameters, and let $G \in L^\infty(\Z_N)$ be function taking values in the interval $I := [-2^{2^{k-1}}, 2^{2^{k-1}} ]$.
Then there exists a $\sigma$-algebra $\B_{\eps,\eta}(G)$ with the following properties:
\begin{itemize}

\item \textup{($G$ lies in its own $\sigma$-algebra)} For any $\sigma$-algebra $\B$, we have
\begin{equation}\label{trivial}
\| G - \E( G | \B \vee \B_{\eps,\eta}(G)) \|_{L^\infty(\Z_N)} \leq \eps.
\end{equation}

\item \textup{(Bounded complexity)} $\B_{\eps,\eta}(G)$ is generated by at most $O(1/\eps)$ atoms.

\item \textup{(Approximation by continuous functions of $G$)}  If $A$ is any atom in $\B_{\eps,\eta}(G)$, then there exists a continuous function
$\Psi_A: I \to [0,1]$ such that 
\begin{equation}\label{nup-bound}
\| ({\bf 1}_A - \Psi_A(G))(\nu + 1) \|_{L^1(\Z_N)} = O(\eta).
\end{equation}
Furthermore, $\Psi_A$ lies in a fixed compact set $E = E_{\eps,\eta}$ of $C^0(I)$ \textup{(}which is independent of $F$, $\nu$, $N$, or $A$\textup{)}.

\end{itemize}

\end{proposition}

\proof 
Observe from Fubini's theorem and \eqref{numean} that
$$ 
\int_0^1
\sum_{n \in \Z} \E\big( {\bf 1}_{G(x) \in [\eps(n-\eta+\alpha), \eps(n+\eta+\alpha)]} (\nu(x)+1) \; \big| \; x \in \Z_N \big)\ d\alpha
= 2\eta \E( \nu(x)+1 | x \in \Z_N ) = O(\eta) $$
and hence by the pigeonhole principle there exists $0 \leq \alpha \leq 1$ such that
\begin{equation}\label{dual-shift}
\sum_{n \in \Z} \E\big( {\bf 1}_{G(x) \in [\eps(n-\eta+\alpha), \eps(n+\eta+\alpha)]} (\nu(x)+1) \; \big| \; x \in \Z_N \big) = O( \eta ).
\end{equation}
We now set $\B_{\eps,\eta}(G)$ to be the $\sigma$-algebra whose atoms are the sets $G^{-1}([\eps (n+\alpha), \eps(n+1+\alpha)))$ for $n \in \Z$.  This is well-defined since the intervals $[\eps (n+\alpha), \eps(n+1+\alpha))$ tile the real line.\vs

\ni It is clear that if $\B$ is an arbitrary $\sigma$-algebra, then on any atom of $\B \vee \B_{\eps,\eta}(G)$, the function $G$ takes values
in an interval of diameter $\eps$, which yields \eqref{trivial}.  
Now we verify the approximation by continuous functions property.  Let $A := G^{-1}([\eps (n+\alpha), \eps(n+1+\alpha)))$ be an atom.
Since $G$ takes values in $I$, we may assume that $n = O(1/\eps)$, since $A$ is empty otherwise; note that this already establishes
the bounded complexity property.  Let $\psi_\eta: \R \to [0,1]$ be a fixed continuous cutoff function which equals 1 on $[\eta, 1-\eta]$ and
vanishes outside of $[-\eta,1+\eta]$, and define $\Psi_A(x) := \psi_\eta( \frac{x}{\eps} - n - \alpha )$.  Then it is clear that $\Psi_A$
ranges over a compact subset $E_{\eps,\eta}$ of $C^0(I)$ (because $n$ and $\alpha$ are bounded).  Furthermore from \eqref{dual-shift} it is clear
that we have \eqref{nup-bound}.  The claim follows.
\endproof\vs

\ni We now specialise to the case when the functions $G$ are basic Gowers anti-uniform functions.

\begin{proposition}\label{bohr-sigma} Let $\nu$ be a $k$-pseudorandom measure.  Let $K \geq 1$ be an fixed integer and let $\D F_1, \ldots, \D F_K \in L^{\infty}(\Z_N)$ be 
basic Gowers anti-uniform functions.  Let $0 < \eps < 1$ and $0 < \eta < 1/2$ be parameters, and let $\B_{\eps,\eta}(\D F_j)$, $j=1,\ldots,K$,
be constructed as in Proposition \ref{sigma-gen}.  Let $\B := \B_{\eps,\eta}(\D F_1) \vee \ldots \vee \B_{\eps,\eta}(\D F_K)$.  Then if $\eta < \eta_0(\epsilon, K)$ is sufficiently small and $N > N_0(\epsilon, K, \eta)$ is sufficiently large we have
\begin{equation}\label{trivial-specific}
\| \D F_j - \E( \D F_j | \B) \|_{L^\infty(\Z_N)} \leq \eps \hbox{ for all } 1 \leq j \leq K.
\end{equation}
Furthermore there exists a set $\Omega$ which lies in $\B$ such that
\begin{equation}\label{mok-small} \E( (\nu + 1) {\bf 1}_{\Omega} ) = O_{K,\eps}(\eta^{1/2})
\end{equation}
and such that 
\begin{equation}\label{unif-mu}
 \| (1 - {\bf 1}_{\Omega}) \E(\nu - 1 | \B) \|_{L^\infty(\Z_N)} = O_{K,\eps}(\eta^{1/2}).
\end{equation}
\end{proposition}

\noindent\textit{Remark.} We strongly recommend that here and in subsequent arguments the reader pretend that the exceptional set $\Omega$ is empty; in practice we shall be able to set $\eta$ small enough
that the contribution of $\Omega$ to our calculations will be negligible.\vs

\proof
The claim \eqref{trivial-specific} follows immediately from \eqref{trivial}.  Now we prove \eqref{mok-small} and \eqref{unif-mu}.
Since each of the $\B_{\eps,\eta}(\D F_j)$ are generated by $O(1/\eps)$ atoms, we see that $\B$ is generated by $O_{K,\eps}(1)$ atoms.
Call an atom $A$ of $\B$ \emph{small} if $\E( (\nu + 1) {\bf 1}_A) \leq \eta^{1/2}$, and let $\Omega$ be the union of all the small atoms.  Then clearly $\Omega$ lies in $\B$ and obeys \eqref{mok-small}.  To prove the remaining claim \eqref{unif-mu}, it suffices to show that
\begin{equation}\label{sigma-smooth}
\frac{ \E( (\nu - 1) {\bf 1}_A ) }{ \E( {\bf 1}_A ) } = \E( \nu - 1 | A ) = o_{K,\eps,\eta}(1) + O_{K,\eps}(\eta^{1/2})
\end{equation}
for all atoms $A$ in $\B$ which are not small.  However, by definition of ``small'' we have
$$ \E( (\nu - 1) {\bf 1}_A ) + 2 \E( {\bf 1}_A ) = \E( (\nu+1) {\bf 1}_A ) \geq \eta^{1/2}.$$
Thus to complete the proof of \eqref{sigma-smooth} it will suffice (since $\eta$ is small and $N$ is large) to show that
\begin{equation}\label{nubf}
\E( (\nu - 1) {\bf 1}_A ) = o_{K,\eps,\eta}(1) + O_{K,\eps}(\eta).
\end{equation}
On the other hand, since $A$ is the intersection of $K$ atoms $A_1, \ldots, A_K$ from $\B_{\eps,\eta}(\D F_1), \ldots$, $\B_{\eps,\eta}(\D F_K)$
respectively, we see from Proposition \ref{sigma-gen} and an easy induction argument (involving H\"older's inequality and the triangle inequality) that we can find a continuous 
function $\Psi_A: I^K \to [0,1]$ such that
$$\| (\nu + 1) ({\bf 1}_A - \Psi_A(\D F_1, \ldots, \D F_K)) \|_{L^1(\Z_N)} = O_K(\eta),$$
so in particular
$$ \| (\nu - 1) ({\bf 1}_A - \Psi_A(\D F_1, \ldots, \D F_K)) \|_{L^1(\Z_N)} = O_K(\eta).$$
Furthermore one can easily ensure that $\Psi_A$ lives in a compact set $E_{\eps,\eta,K}$ of $C^0(I^K)$.  From this and Proposition \ref{weier} we have
$$ \E( (\nu - 1) \Psi_A(\D F_1, \ldots, \D F_K)) ) = o_{K,\eps,\eta}(1)$$
since $N$ is assumed large depending in $K,\eps,\eta$, and the claim \eqref{nubf} now follows from the triangle inequality.
\endproof\vs

\ni \emph{Remarks.} This $\sigma$-algebra $\B$ is closely related to the (relatively) compact $\sigma$-algebras studied in the ergodic theory proof of Szemer\'edi's theorem, see for instance \cite{furst,FKO}.  In the case $k=3$ they are closely connected to the Kronecker factor of an ergodic system, and for higher $k$ they are related to $(k-2)$-step nilsystems, see e.g. \cite{host-kra2,ziegler}.

\section{A Furstenberg tower, and the proof of Theorem \ref{main}}\label{sec7}

\ni We now have enough machinery to deduce Theorem \ref{main} from Proposition \ref{sz}.  The key proposition is the following decomposition, which splits an arbitrary function into Gowers uniform and Gowers anti-uniform components (plus a negligible error).

\begin{proposition}[Generalised Koopman-von Neumann structure theorem]\label{kvn-decomp}  Let $\nu$ be a $k$-pseudorandom measure,
and let $f \in L^1(\Z_N)$ be a non-negative function satisfying $0 \leq f(x) \leq \nu(x)$ for all $x \in \Z_N$. Let $0 < \eps \ll 1$ be a small parameter, and assume $N > N_0(\eps)$ is sufficiently large.
Then there exists a $\sigma$-algebra $\B$ and an exceptional set $\Omega \in \B$ such that
\begin{itemize}
\item \textup{(}smallness condition\textup{)} 
\begin{equation}\label{small-exceptions-final}
 \E( \nu {\bf 1}_{\Omega} ) = o_\eps(1);
\end{equation}
\item \textup{(}$\nu$ is uniformly distributed outside of $\Omega$\textup{)}
\begin{equation}\label{mu-regular-final}
 \| (1 - {\bf 1}_{\Omega}) \E( \nu - 1 | \B) \|_{L^\infty} = o_{\eps}(1)
\end{equation} and
\item \textup{(}Gowers uniformity estimate\textup{)}
\begin{equation}\label{unif-ok}
 \| (1 - {\bf 1}_{\Omega}) (f - \E(f|\B)) \|_{U^{k-1}} \leq \eps^{1/2^k}.
\end{equation}
\end{itemize}
\end{proposition}

\ni\emph{Remarks.}  As in the previous section, the exceptional set $\Omega$ should be ignored on a first reading. The ordinary Koopman-von Neumann theory in ergodic theory asserts, among other things, that any function $f$
on a measure-preserving system $(X,\B, T, \mu)$ can be orthogonally decomposed into a ``weakly mixing part'' $f-\E(f|\B)$ (in which $f-\E(f|\B)$
is asymptotically orthogonal to its shifts $T^n (f-\E(f|\B))$ on the average) and an ``almost periodic part'' $\E(f|\B)$ (whose shifts form a precompact
set); here $\B$ is the \emph{Kronecker factor}, i.e. the $\sigma$-algebra generated by the almost periodic functions (or equivalently, by the eigenfunctions of $T$).  
This is somewhat related to the $k=3$ case of the above Proposition, hence our labeling of that proposition as a generalised Koopman-von Neumann
theorem.  A slightly more quantitative analogy for the $k=3$ case would be the assertion that any function bounded by a pseudorandom measure can
be decomposed into a Gowers uniform component with small Fourier coefficients, and a Gowers anti-uniform component which consists of only a few Fourier coefficients
(and in particular is bounded). For related ideas see \cite{bourg2,green,green-tao}.\vs

\ni\emph{Proof of Theorem \ref{main} assuming Proposition \ref{kvn-decomp}.}  Let $f$, $\delta$ be as in Theorem \ref{main}, and let $0 < \eps \ll \delta$ be a parameter to be chosen later.  Let $\B$ be as in the above decomposition, and write $\funif := (1-{\bf 1}_\Omega)(f-\E(f|\B))$ and
$\fanti := (1-{\bf 1}_\Omega)\E(f|\B)$ (the subscript $U$ stands for Gowers uniform, and $U^\perp$ for Gowers anti-uniform).  Observe from \eqref{small-exceptions-final}, \eqref{f-bound-2},
\eqref{f-density} and the measurability of $\Omega$ that
$$ \E( \fanti ) = \E( (1-{\bf 1}_\Omega) f) \geq \E(f) - \E( \nu {\bf 1}_\Omega ) \geq \delta - o_\eps(1).$$
Also, by \eqref{mu-regular-final} we see that $\fanti$ is bounded above by $1+o_{\eps}(1)$.  Since $f$ is non-negative, $\fanti$ is also.
We may thus\footnote{There is an utterly trivial issue which we have ignored here, which is that $\fanti$ is not bounded above by $1$ but by $1 + o_\eps(1)$, and that the density is bounded below by $\delta - o_\eps(1)$ rather than $\delta$. One can easily get around this by modifying $\fanti$ by $o_\eps(1)$ before applying Proposition \ref{sz}, incurring a net error of $o_\eps(1)$ at the end since $\fanti$ is bounded.} apply Proposition \ref{sz}  to obtain
\[
 \E\big( \fanti(x) \fanti(x+r) \ldots \fanti(x+(k-1)r) \; \big| \; x,r \in \mathbb{Z}_N \big) \geq c(k,\delta) -  o_{\eps}(1) - o_{k,\delta}(1).\]
On the other hand, from \eqref{unif-ok} we have $\|\funif \|_{U^{k-1}} \leq \eps^{1/2^k}$; since $(1-{\bf 1}_\Omega)f$ is bounded by $\nu$
and $\fanti$ is bounded by $1+o_{\eps}(1)$, we thus see that $\funif$ is pointwise bounded by $\nu+1+o_{\eps}(1)$.
Applying the generalised von Neumann theorem (Proposition \ref{vn}) we thus see that
 \[
 \E\big( f_0(x) f_1(x+r) \ldots f_{k-1}(x+(k-1)r) \; \big| \; x,r \in \mathbb{Z}_N \big) = O(\eps^{1/2^k}) + o_{\eps}(1)\]
whenever each $f_j$ is equal to $\funif$ or $\fanti$, with at least one $f_j$ equal to $\funif$.  Adding these two estimates together
we obtain
\[
 \E( \tilde f(x) \tilde f(x+r) \ldots \tilde f(x+(k-1)r) | x,r \in \mathbb{Z}_N ) \geq c(k,\delta) - O(\eps^{1/2^k}) - o_{\eps}(1) - o_{k,\delta}(1),\]
where $\tilde f := \funif + \fanti = (1-{\bf 1}_\Omega) f$.  But since $0 \leq (1-{\bf 1}_\Omega) f \leq f$ we obtain
$$
 \E( f(x) f(x+r) \ldots f(x+(k-1)r) | x,r \in \mathbb{Z}_N ) \geq c(k,\delta) - O(\eps^{1/2^k}) - o_{\eps}(1)- o_{k,\delta}(1).$$
Since $\eps$ can be made arbitrarily small (as long as $N$ is taken sufficiently large), the error terms on the right-hand side can be taken to
be arbitrarily small by choosing $N$ sufficiently large depending on $k$ and $\delta$.  The claim follows.
\endproof\vspace{11pt}

\ni To complete the proof of Theorem \ref{main}, it suffices to prove Proposition \ref{kvn-decomp}.  To construct the $\sigma$-algebra $\B$
required in the Proposition, we will use the philosophy laid out by Furstenberg in his ergodic structure theorem (see \cite{furst,FKO}), which decomposes any measure-preserving system into a weakly-mixing extension of a tower of compact extensions.  
In our setting, the idea is roughly speaking as follows.  We initialise $\B$ to be the trivial $\sigma$-algebra $\B = \{ \emptyset, \Z_N \}$.
If the function $f - \E(f|\B)$ is already Gowers uniform (in the sense of \eqref{unif-ok}), then we can terminate the algorithm.  Otherwise, we use the
machinery of dual functions, developed in \S \ref{sec6}, to locate a Gowers anti-uniform function $\D F_1$ which has some non-trivial correlation with $f$, and add the level sets of $\D F_1$ to the $\sigma$-algebra $\B$; the non-trivial correlation property will ensure that the $L^2$ norm of $\E(f|\B)$ increases by
a non-trivial amount during this procedure, while the pseudorandomness of $\nu$ will ensure that $\E(f|\B)$ remains uniformly
bounded.  We then repeat the above algorithm until $f - \E(f|\B)$ becomes sufficiently Gowers uniform,
at which point we terminate the algorithm.  In the original ergodic theory arguments of Furstenberg this algorithm was not guaranteed to terminate,
and indeed one required the axiom of choice (in the guise of Zorn's lemma) in order to conclude\footnote{For the specific purpose of $k$-term recurrence, i.e. finding progressions of length $k$, one only needs to run Furstenberg's algorithm for a finite number of steps depending on $k$, and so Zorn's lemma is not needed in this application.  We thank Bryna Kra for pointing out this subtlety.} the structure theorem.  However, in our setting
we can 
terminate in a bounded number of steps (in fact in at most $2^{2^k}/\eps + 2$ steps), because there is a quantitative
 $L^2$-increment to the bounded function $\E(f|\B)$ at each stage.\vs

\ni Such a strategy will be familiar to any reader acquainted with the proof of Szemer\'edi's regularity lemma \cite{szem-reg}. This is no coincidence: there is in fact a close connection between regularity lemmas such as those in  \cite{gowers-reg,green-reg,szem-reg} and ergodic theory of the type we have brushed up against in this paper.  Indeed there are strong analogies between all of
the known proofs of Szemer\'edi's theorem, despite the fact that they superficially appear to use
very different techniques.\vs

\ni We turn to the details.  To prove Proposition \ref{kvn-decomp}, we will iterate the following somewhat technical 
Proposition, which can be thought of as a $\sigma$-algebra variant of Lemma \ref{antiuniform}.

\begin{proposition}[Iterative Step]\label{relative} 
Let $\nu$ be a $k$-pseudorandom measure,
and let $f \in L^1(\Z_N)$ be a non-negative function satisfying $0 \leq f(x) \leq \nu(x)$ for all $x \in \Z_N$.
Let $0 <\eta \ll \eps \ll 1$ be small numbers, and let $K \geq 0$
be an integer.  Suppose that $\eta < \eta_0(\eps, K)$ is sufficiently small and that $N > N_0(\eps, K,\eta)$ is sufficiently large.
Let $F_1,\ldots,F_K \in L^1(\Z_N)$ be a collection of functions
obeying the pointwise bounds
\begin{equation}\label{f-disc-ok}
 |F_{j}(x)| \leq (1 + O_{K,\eps}(\eta^{1/2})) (\nu(x) + 1)
\end{equation}
for all $1 \leq j \leq K$ and $x \in \Z_N$.  Let $\B_K$ be the $\sigma$-algebra
\begin{equation}\label{BK-def}
 \B_K :=  \B_{\eps,\eta}(\D F_1) \vee \ldots \vee \B_{\eps,\eta}(\D F_K)
 \end{equation}
where $\B_{\eps,\eta}(\D F_j)$ is as in Proposition \ref{sigma-gen}, and suppose that there exists a set $\Omega_K$ in $\Z_N$
obeying 
\begin{itemize}
\item \textup{(}smallness bound\textup{)}
\begin{equation}\label{om-small} \E( (\nu + 1) {\bf 1}_{\Omega_K} ) = O_{K,\eps}(\eta^{1/2})
\end{equation}
and \item \textup{(}uniform distribution bound\textup{)}
\begin{equation}\label{uniform-dist}
\| (1 - {\bf 1}_{\Omega_{K}}) \E( \nu - 1 | \B_{K} )\|_{L^\infty(\Z_N)} = O_{K,\eps}(\eta^{1/2}).
\end{equation}
\end{itemize}
Set
\begin{equation}\label{FK1-def}
F_{K+1} := (1 - {\bf 1}_{\Omega_K}) (f - \E(f|\B_K))
\end{equation}
and suppose that $F_{K+1}$ obeys the non-Gowers-uniformity estimate
\begin{equation}\label{non-terminate}
 \| F_{K+1} \|_{U^{k-1}} > \eps^{1/2^k}.
\end{equation}
Then we have the estimates
\begin{equation}\label{fbk-bound}
\| (1 - {\bf 1}_{\Omega_K}) \E(f|\B_K) \|_{L^\infty(\Z_N)} \leq 1 + O_{K,\eps}(\eta^{1/2})
\end{equation}
\textup{(}$\E (f | \B_K)$ is bounded outside $\Omega_K$\textup{)}
and
\begin{equation}\label{f-disc-ok-2}
 |F_{K+1}(x)| \leq (1 + O_{K,\eps}(\eta^{1/2})) (\nu(x) + 1).
\end{equation}
Furthermore, if we let $\B_{K+1}$ be the $\sigma$-algebra
\begin{equation}\label{BK1-def}
 \B_{K+1} := \B_K \vee \B_{\eps,\eta}(\D F_{K+1}) = \B_{\eps,\eta}(\D F_1) \vee \ldots \vee \B_{\eps,\eta}(\D F_{K+1})
\end{equation}
then there exists a set $\Omega_{K+1} \supseteq \Omega_K$ obeying 
\begin{itemize}
\item \textup{(}smallness bound\textup{)}
\begin{equation}\label{om-small-2} \E( (\nu + 1) {\bf 1}_{\Omega_{K+1}} ) = O_{K,\eps}(\eta^{1/2});
\end{equation}
\item \textup{(}uniform distribution bound\textup{)}
\begin{equation}\label{uniform-dist-2}
\| (1 - {\bf 1}_{\Omega_{K+1}}) \E( \nu - 1 | \B_{K+1} )\|_{L^\infty(\Z_N)} = O_{K,\eps}(\eta^{1/2}).
\end{equation}
and such that we have 
\item \textup{(}energy increment property\textup{)}
\begin{equation}\label{fj-increment}
 \| (1 - {\bf 1}_{\Omega_{K+1}}) \E(f|\B_{K+1}) \|_{L^2(\Z_N)}^2 \geq \| (1 - {\bf 1}_{\Omega_K}) \E(f|\B_K) \|_{L^2(\Z_N)}^2 + 2^{-2^k + 1} \eps.
\end{equation}
\end{itemize}
\end{proposition}

\ni\emph{Remark.}  If we ignore the exceptional sets $\Omega_K$, $\Omega_{K+1}$, this proposition is asserting the following: if
$f$ is not ``relatively weakly mixing'' with respect to the $\sigma$-algebra $\B_K$, in the sense that the component $f - \E(f|\B_K)$ of $f$ which is orthogonal to $\B_K$ is not Gowers-uniform, then we can refine $\B_K$ to a slightly more complex $\sigma$-algebra $\B_{K+1}$ such that the $L^2(\Z_N)$ norm (energy) of
$\E(f|\B_{K+1})$ is larger than $E(f|\B_K)$ by some quantitative amount.  Furthermore, $\nu$ remains uniformly distributed with respect to $\B_{K+1}$.\vs

\ni\textit{Proof of Proposition \ref{kvn-decomp} assuming Proposition \ref{relative}.}
Fix $\eps$, and let $K_0$ be the smallest integer greater than $2^{2^k}/\eps+1$; this quantity will be the upper bound for the number of iterations of an algorithm which we shall give shortly.  
We shall need a parameter $0 < \eta \ll \eps$ which we shall choose 
later (we assume $\eta < \eta_0(\eps, K_0)$, and then we shall assume $N > N_0(\eta,\eps)$ is sufficiently large.\vs

\ni To construct $\B$ and $\Omega$ we shall, for some $K \in [0,K_0]$, iteratively construct a sequence of basic Gowers anti-uniform
functions $\D F_1, \ldots, \D F_K$ on $\mathbb{Z}_N$ together with 
exceptional sets $\Omega_0 \subseteq\Omega_1 \subseteq\ldots \subseteq\Omega_K \subseteq\mathbb{Z}_N$ in the following manner.
\begin{itemize}

\item Step 0.  Initialise $K = 0$ and $\Omega_0 := \emptyset$.  (We will later increment the value of $K$).
\item Step 1.  Let $\B_K$ and $F_{K+1}$ be defined by \eqref{BK-def} and \eqref{FK1-def} respectively.  Thus for instance when $K=0$
we will have $\B_0 = \{\emptyset,\Z_N\}$ and $F_1 = f - \E(f)$.  Observe that in the $K=0$ case, the estimates \eqref{f-disc-ok}, \eqref{om-small}, \eqref{uniform-dist} are trivial (the latter bound following from \eqref{numean}).  As we shall see, these three 
estimates will be preserved throughout the algorithm.
\item Step 2.  If the estimate \eqref{non-terminate} \emph{fails}, or in other words that
$$ \| F_{K+1} \|_{U^{k-1}} \leq \eps^{1/2^k},$$
then we set $\Omega := \Omega_K$ and $\B := \B_K$, and successfully terminate the algorithm.
\item Step 3.  If instead \eqref{non-terminate} holds, then we define $\B_{K+1}$ by \eqref{BK1-def}.  (Here we of course need $K \leq K_0$, but
this will be guaranteed by Step 4 below).  We then invoke Proposition \ref{relative} to locate an exceptional 
set $\Omega_{K+1} \supset \Omega_K$ in $\B_{K+1}$ obeying the conditions\footnote{Of course, the constants in the $O()$ bounds are different at each stage of this iteration, but we are allowing these constants to depend on $K$, and $K$ will ultimately be bounded by $K_0$, which depends only on $\eps$ and $k$.} \eqref{om-small-2}, \eqref{uniform-dist-2}, and for which
we have the energy increment property \eqref{fj-increment}.  Also, we have \eqref{f-disc-ok-2}.
\item Step 4.  Increment $K$ to $K+1$; observe from construction that the estimates \eqref{f-disc-ok}, \eqref{om-small}, \eqref{uniform-dist} will be preserved when doing so.  If we now have $K > K_0$, then we terminate the algorithm with an error; otherwise, return to Step 1.

\end{itemize}

\ni \emph{Remarks.} The integer $K$ indexes the iteration number of the algorithm, thus we begin with the zeroth iteration when $K=0$, then the first
iteration when $K=1$, etc.  It is worth noting that apart from $O_{K,\eps}(\eta^{1/2})$ error terms, none of the 
bounds we will encounter while executing this algorithm will actually depend on $K$.  As we shall see, this algorithm will terminate well before $K$
reaches $K_0$ (in fact, for the application to the primes, the Hardy-Littlewood prime tuples conjecture implies that this algorithm will terminate at the first step $K=0$).
\vs

\ni Assuming Proposition \ref{relative}, we see that this algorithm will terminate after finitely many steps with one of two outcomes: either 
it will terminate successfully in Step 2 for some $K \leq K_0$, or else it will terminate with an error in Step 4 when $K$ 
exceeds $K_0$.  Assume for the moment that the former case occurs.  Then it is clear that at the successful conclusion of this algorithm,
we will have generated a $\sigma$-algebra $\B$ and an exceptional set $\Omega$ with the properties required for 
Proposition \ref{kvn-decomp}, with
error terms of $O_{K,\eps}(\eta^{1/2})$ instead of $o_\eps(1)$, if $N > N_0(\eta, K, \eps)$.
But by making $\eta$ decay sufficiently slowly to zero, we can replace the $O_{K,\eps}(\eta^{1/2})$ bounds by $o_\eps(1)$; note that 
the dependence of the error terms on $K$ will not be relevant since $K$ is bounded by $K_0$, which depends only on $\eps$.
\vs

\ni To conclude the proof of Proposition \ref{kvn-decomp} it will thus suffice to show that the above algorithm does not
terminate with an error.  Suppose for a contradiction that the algorithm ran until the $K_0^{\mbox{\scriptsize th}}$ iteration before terminating with an error in Step 4.
Then if we define the energies $E_K$ for $0 \leq K \leq K_0+1$ by the formula
$$ E_K := \| (1 - {\bf 1}_{\Omega_K}) \E(f|\B_K) \|_{L^2(\Z_N)}^2$$
then we see from \eqref{fj-increment} that
\begin{equation}\label{energy-increment}
E_{K+1} \geq E_K + 2^{-2^k + 1} \eps \hbox{ for all } 0 \leq K \leq K_0
\end{equation}
(for instance).  Also, by \eqref{fbk-bound} we have
$$ 0 \leq E_K \leq 1 + O_{K,\eps}(\eta^{1/2}) \hbox{ for all } 0 \leq K \leq K_0.$$
If $\eta < \eta_0(K,\eps)$ is sufficiently small, these last two statements contradict one another for $K = K_0$. Thus the above algorithm cannot reach the $K_0^{\mbox{\scriptsize th}}$ iteration, and instead terminates
successfully at Step 2.  This completes the proof of Proposition \ref{kvn-decomp}.\endproof\vs

\ni The only remaining task is to prove Proposition \ref{relative}.\vs

\ni\textit{Proof of Proposition \ref{relative}.}  Let $\nu$, $f$, $K$, $\eps$, $\eta$, $F_1,\ldots,F_K, F_{K+1}$, $\B_K$, $\Omega_K$, $\B_{K+1}$ 
be as in the proposition.  We begin by proving the bounds \eqref{fbk-bound}, \eqref{f-disc-ok-2}.  From \eqref{uniform-dist} we have
$$
\| (1 - {\bf 1}_{\Omega_{K}}) \E( \nu | \B_{K} )\|_{L^\infty} \leq 1 + O_{K,\eps}(\eta^{1/2});$$
since $f$ is non-negative and bounded pointwise by $\nu$, we obtain \eqref{fbk-bound}.  The bound \eqref{f-disc-ok-2} then follows
from \eqref{fbk-bound} and \eqref{FK1-def}, where we again use that $f$ is non-negative and bounded pointwise by $\nu$.
This shows in particular that $\D F_1, \ldots$, $\D F_{K_1+1}$ are basic Gowers anti-uniform functions (up to multiplicative 
errors of $1 + O_{K_1,\eps}(\eta^{1/2})$, which are negligible).\vs

\ni Applying Lemma \ref{antiuniform} (scaling out the multiplicative error of $1 + O_{K,\eps}(\eta^{1/2})$) and using
\eqref{f-disc-ok} and \eqref{f-disc-ok-2} we conclude that
\begin{equation}\label{FJ-bound}
\| \D F_j \|_{L^\infty(\Z_N)} \leq 2^{2^{k-1}-1} + O_{K,\eps}(\eta^{1/2}) \hbox{ for all } 0 \leq j \leq K+1, 
\end{equation}
since we are assuming $N$ to be large depending on $K$, $\eps$, $\eta$.\vs

\ni We now apply Proposition \ref{bohr-sigma} (absorbing the multiplicative errors of $1+O_{K,\eps}(\eta^{1/2})$) 
to conclude that we may find a set $\Omega$ in $\B_{K+1}$ such that
$$ \E( (\nu + 1) {\bf 1}_{\Omega} ) = O_{K,\eps}(\eta^{1/2})$$
and
$$
 \| (1 - {\bf 1}_{\Omega}) \E(\nu - 1 | \B_{K+1}) \|_{L^\infty} = O_{K,\eps}(\eta^{1/2}).$$
If we then set $\Omega_{K+1} := \Omega_K \cup \Omega$, then we can verify \eqref{om-small-2} and
\eqref{uniform-dist-2} from \eqref{om-small} and \eqref{uniform-dist}.\vs
 
\ni It remains to verify \eqref{fj-increment}, the energy increment property, that is to say the statement
\begin{equation}\label{restate} \| (1 - {\bf 1}_{\Omega_{K+1}}) \E(f|\B_{K+1}) \|_{L^2(\Z_N)}^2 \geq \| (1 - {\bf 1}_{\Omega_K}) \E(f|\B_K) \|_{L^2(\Z_N)}^2 + 2^{-2^k + 1} \eps.\end{equation}
To do this we exploit the hypothesis $\Vert F_{K+1} \Vert_{U^{k-1}} > \eps^{1/2^k}$, which was \eqref{non-terminate}.
By Lemma \ref{antiuniform} and the definition \eqref{FK1-def} we have
 \[
|\left\langle (1 - {\bf 1}_{\Omega_K}) (f-\E(f|\B_K)), \D F_{K+1} \right\rangle| = |\langle F_{K+1}, \D F_{K+1} \rangle| = \| F_{K+1} \|_{U^{k-1}}^{2^{k-1}} \geq \eps^{1/2}.
\]
On the other hand, from the bounds \eqref{f-disc-ok}, \eqref{om-small} and \eqref{FJ-bound} we have
\begin{eqnarray*}
& & 
 \big|\left\langle ({\bf 1}_{\Omega_{K+1}} - {\bf 1}_{\Omega_{K}}) (f-\E(f|\B_K)), \D F_{K+1} \right\rangle \big|  \\ & & \qquad\qquad\qquad\qquad \leq  \Vert \D F_{K+1} \Vert_{\infty} \E \big(({\bf 1}_{\Omega_{K+1}} - {\bf 1}_{\Omega_{K}}) |f - \E(f | \B_K)|\big) \\ &  & \qquad\qquad\qquad\qquad = O_{K,\eps}(1) \E \big(({\bf 1}_{\Omega_{K+1}} - {\bf 1}_{\Omega_{K}}) (\nu + 1)\big) \; = \; O_{K,\eps}(\eta^{1/2}),\end{eqnarray*}
while from \eqref{trivial} and  \eqref{fbk-bound} we have
\begin{eqnarray*}
 & & \big|\left\langle (1 - {\bf 1}_{\Omega_{K+1}}) (f-\E(f|\B_K)), \D F_{K+1} - \E(\D F_{K+1}|\B_{K+1}) \right\rangle\big| \\ & & \qquad\qquad\qquad\qquad \leq  \Vert \D F_{K+1} - \E(\D F_{K+1} | \B_{K+1}) \Vert_{\infty} \E \big((1 - {\bf 1}_{\Omega_{K+1}})|f - \E(f | \B_K)|\big) \\ & & \qquad\qquad\qquad\qquad \leq  O(\eps) \E\big((1 - {\bf 1}_{\Omega_{K+1}})(\nu + 1)\big) \; = \; O(\epsilon).
\end{eqnarray*}
By the triangle inequality we thus have
$$
 |\left\langle (1 - {\bf 1}_{\Omega_{K+1}}) (f-\E(f|\B_K)), \E(\D F_{K+1}|\B_{K+1}) \right\rangle| \geq \eps^{1/2} - O_{K,\eps}(\eta^{1/2}) - O(\eps).
$$
But since $(1 - {\bf 1}_{\Omega_{K+1}})$, $\E(\D F_{K+1}|\B_{K+1})$, and $\E(f|\B_K)$ are all measurable in $\B_{K+1}$, we can replace
$f$ by $\E(f|\B_{K+1})$, and so
$$
 |\left\langle (1 - {\bf 1}_{\Omega_{K+1}}) (\E(f|\B_{K+1})-\E(f|\B_K)), \E(\D F_{K+1}|\B_{K+1}) \right\rangle| \geq \eps^{1/2} - O_{K,\eps}(\eta^{1/2}) - O(\eps).
$$
By the Cauchy-Schwarz inequality and \eqref{FJ-bound} we obtain
\begin{equation}\label{pythagoras-targ}
 \| (1 - {\bf 1}_{\Omega_{K+1}}) (\E(f|\B_{K+1}) - \E(f|\B_K)) \|_{L^2(\Z_N)}
\geq 2^{-2^{k-1}+1} \eps^{1/2} - O_{K,\eps}(\eta^{1/2}) - O(\eps).
\end{equation}
Morally speaking, this implies \eqref{fj-increment} thanks to Pythagoras's theorem, but the presence of the exceptional sets $\Omega_K$ and
$\Omega_{K+1}$ means that we have to exercise caution, especially since we have no $L^2$ control on $\nu$.\vs

\ni Recalling that $\E(f | \B_K) \leq 1 + O_{\eps,K}(\eta^{1/2})$ outside $\Omega_K$ (cf. \eqref{fbk-bound}), we observe that if $\eta < \eta_0(\eps, K)$ is sufficiently small then
\[ \| ({\bf 1}_{\Omega_{K+1}} - {\bf 1}_{\Omega_K}) \E(f|\B_K) \|_{L^2(\Z_N)}^2 \leq 2 \Vert {\bf 1}_{\Omega_{K+1}} - {\bf 1}_{\Omega_K} \Vert_2^2 \leq 2 \Vert {\bf 1}_{\Omega_{K+1}} - {\bf 1}_{\Omega_K} \Vert_1 \leq 2\E {\bf 1}_{\Omega_{K+1}},\]
which, by \eqref{om-small}, is $O_{K,\eps}(\eta^{1/2})$.
By the triangle inequality (and \eqref{fbk-bound}) we thus see that to prove \eqref{fj-increment} it will suffice
to prove that
\[
 \| (1 - {\bf 1}_{\Omega_{K+1}}) \E(f|\B_{K+1}) \|_{L^2(\Z_N)}^2 \geq \| (1 - {\bf 1}_{\Omega_{K+1}}) \E(f|\B_K) \|_{L^2(\Z_N)}^2 + 
 2^{-2^k + 2} \eps - O_{K,\eps}(\eta^{1/2})
- O(\eps^{3/2}),\]
since we can absorb the error terms $- O_{K,\eps}(\eta^{1/2})
- O(\eps^{3/2})$ into the $2^{-2^k + 2} \eps$ term by choosing $\eps$ sufficiently small depending on $k$, and $\eta$ sufficiently small depending
on $K,\eps$.\vs

\ni We write the left-hand side as
$$ \| (1 - {\bf 1}_{\Omega_{K+1}}) \E(f|\B_K) + (1 - {\bf 1}_{\Omega_{K+1}}) (\E(f|\B_{K+1}) - \E(f|\B_K) \|_{L^2(\Z_N)}^2$$
which can be expanded using the cosine rule as
\begin{align*}
 \| (1 - {\bf 1}_{\Omega_{K+1}}) \E(f|\B_K) \|_{L^2(\Z_N)}^2 &+ \| (1 - {\bf 1}_{\Omega_{K+1}}) (\E(f|\B_{K+1}) - \E(f|\B_K)) \|_{L^2(\Z_N)}^2 \\
&+ 2 \left\langle (1 - {\bf 1}_{\Omega_{K+1}}) \E(f|\B_K), (1 - {\bf 1}_{\Omega_{K+1}}) (\E(f|\B_{K+1}) - \E(f|\B_K)) \right\rangle.
\end{align*}
Therefore by \eqref{pythagoras-targ} it will suffice to show the approximate orthogonality relationship
$$
\left\langle (1 - {\bf 1}_{\Omega_{K+1}}) \E(f|\B_K), 
(1 - {\bf 1}_{\Omega_{K+1}}) (\E(f|\B_{K+1}) - \E(f|\B_K)) \right\rangle
= O_{K,\eps}(\eta^{1/2}).
$$
Since $(1 - {\bf 1}_{\Omega_{K+1}})^2 = (1 - {\bf 1}_{\Omega_{K+1}})$, this can be rewritten as
\[ \left\langle (1 - {\bf 1}_{\Omega_{K+1}}) \E(f|\B_K), \E(f|\B_{K+1}) - \E(f|\B_K) \right\rangle.\]
Now note that $(1 - {\bf 1}_{\Omega_K}) \E(f|\B_K)$ is measurable with respect to $\B_K$, and hence
orthogonal to $\E(f|\B_{K+1}) - \E(f|\B_K)$, since $\B_K$ is a sub-$\sigma$-algebra of $\B_{K+1}$.  Thus the
above expression can be rewritten as
$$ \left\langle ({\bf 1}_{\Omega_{K+1}} - {\bf 1}_{\Omega_K}) \E(f|\B_K), \E(f|\B_{K+1}) - \E(f|\B_K) \right\rangle.$$
Again, since the left-hand side is measurable with respect to $\B_{K+1}$, we can rewrite this as
$$ \langle ({\bf 1}_{\Omega_{K+1}} - {\bf 1}_{\Omega_K}) \E(f|\B_K), f - \E(f|\B_K) \rangle.$$
Since $\E (f | \B_K)(x) \leq 2$ if $\eta < \eta_0(\eps, K)$ is sufficiently small and $x \notin \Omega_K$ (cf. \eqref{fbk-bound}), we may majorise this by
$$ 2\E \big( ({\bf 1}_{\Omega_{K+1}} - {\bf 1}_{\Omega_K}) |f - \E(f|\B_K)| \big).$$  Since we are working on the assumption that $0 \leq f(x) \leq \nu(x)$, we can bound this in turn by
$$ 2\E\big( ({\bf 1}_{\Omega_{K+1}} - {\bf 1}_{\Omega_K}) (\nu + \E(\nu|\B_K)) \big).$$
Since $\E(\nu | \B_K)(x) \leq 2$ for $x \notin \Omega_K$ (cf. \eqref{uniform-dist}) this is no more than
$$ 4\E\big( ({\bf 1}_{\Omega_{K+1}} - {\bf 1}_{\Omega_K}) (\nu + 1) \big),$$
which is $O_{K,\eps}(\eta^{1/2})$ as desired by \eqref{om-small}.  This concludes the proof of 
Proposition \ref{relative}, and hence Theorem \ref{main}.
\endproof

\section{A pseudorandom measure which majorises the primes}\label{sec8}

\ni Having concluded the proof of Theorem \ref{main}, we are now ready to apply it to the specific situation of locating arithmetic
progressions in the primes.  As in almost any additive problem involving the primes, we 
begin by considering the von Mangoldt function $\Lambda$ defined by $\Lambda(n) = \log p$ if $n = p^m$ and $0$ otherwise. 
Actually, for us the higher prime powers $p^2, p^3, \ldots$ will play no r\^ole whatsoever and will be discarded very shortly.\vs

\ni From the prime number theorem we know that the average value of $\Lambda(n)$ is $1 + o(1)$.
In order to prove Theorem \ref{mainthm} (or Theorem \ref{sz-primes}), it would suffice to exhibit a measure $\nu : \mathbb{Z}_N \rightarrow \mathbb{R}^{+}$ such that $\nu(n) \geq c(k)\Lambda(n)$ for some $c(k) > 0$ depending only on $k$, and which is $k$-pseudorandom.  Unfortunately, such a measure
cannot exist because the primes (and the von Mangoldt function) are concentrated on certain residue classes.  Specifically,
for any integer $q > 1$, $\Lambda$ is only non-zero on those $\phi(q)$ residue classes $a \md{q}$ for which $(a,q) = 1$, whereas
a pseudorandom measure can easily be shown to be uniformly distributed across all $q$ residue classes; here of course $\phi(q)$ is
the Euler totient function.  Since $\phi(q)/q$ can be made arbitrarily small, we therefore cannot hope to obtain a pseudorandom majorant
with the desired property $\nu(n) \geq c(k) \Lambda(n)$.\vs

\ni To get around this difficulty we employ a device which we call the $W$-trick\footnote{The reader will observe some similarity between this trick and the use of $\sigma$-algebras in the previous section to remove non-Gowers-uniformity from the system.  Here, of course, the precise obstruction to non-Gowers-uniformity in the primes is very explicit, whereas the exact structure of the $\sigma$-algebras constructed in the previous section are somewhat mysterious.  In the specific case of the primes, we expect (through such conjectures as the Hardy-Littlewood prime tuple conjecture) that the primes are essentially uniform once the obstructions from small primes are removed, and hence the algorithm of the previous section should in fact terminate immediately at the $K=0$ iteration.  However we emphasise that our argument does not give (or require) any progress on this very difficult prime tuple conjecture, as we allow $K$ to be non-zero.}, which effectively removes the arithmetic
obstructions to pseudorandomness arising from the very small primes. Let $\w = \w(N)$ be any function tending 
slowly\footnote{Actually, it will be clear at the end of the proof that we can in fact take $\w$ to be a sufficiently large
number independent of $N$, depending only on $k$, however it will be convenient for now to make $\w$ slowly growing in $N$ in order
to take advantage of the $o(1)$ notation.} to infinity with $N$, so that $1/\w(N) = o(1)$, and let $W = \prod_{p \leq \w(N)} p$ be the product of the primes up to $\w(N)$. Define the modified von Mangoldt function $\tilde \Lambda: \Z^+ \to \R^+$ by
\[ \widetilde{\Lambda}(n) := \left\{
\begin{array}{ll}
\frac{\phi(W)}{W} \log(Wn + 1) & \hbox{ when } Wn+1 \hbox{ is prime}\\
0 & \hbox{ otherwise.}
\end{array}\right.\]
Note that we have discarded the contribution of the prime powers since we ultimately wish to count arithmetic progressions in
the primes themselves.  This $W$-trick exploits the trivial observation that in order to obtain arithmetic progressions in the primes,
it suffices to do so in the modified primes $\{ n \in \Z: Wn+1 \hbox{ is prime} \}$ (at the cost of reducing the number of such progressions by a polynomial factor in $W$ at worst).  We also remark that one could replace $Wn+1$ here by $Wn+b$ for any 
integer $1 \leq b < W$ coprime to $W$ without affecting the arguments which follow.\vs

\ni Observe that if $\w(N)$ is sufficiently slowly growing ($\w(N) \ll \log \log N$ will suffice here) then by Dirichlet's theorem concerning
the distribution of the primes in arithmetic progressions\footnote{In fact, all we need is that $\sum_{N \leq n \leq 2N} \widetilde{\Lambda}(n) \gg N$. Thus one could avoid appealing to the theory of Dirichlet $L$-functions by replacing $n \equiv 1 \md{W}$ by $n \equiv b \md{W}$, for some $b$ coprime to $W$ chosen using the pigeonhole principle.} such as $\{ n: n \equiv 1 \md{W} \}$ we have $\sum_{n \leq N} \widetilde{\Lambda}(n) = N(1 + o(1))$.  With this modification, we can now majorise the primes by a pseudorandom
measure as follows:

\begin{proposition}\label{prime-majorant}
Write $\epsilon_k := 1/2^{k}(k+4)!$, and let $N$ be a sufficiently large prime number.  Then
there is a $k$-pseudorandom measure $\nu : \mathbb{Z}_N \rightarrow \R^+$ such that $\nu(n) \geq k^{-1}2^{-k-5}\widetilde{\Lambda}(n)$ for all 
$\epsilon_k N \leq n \leq 2 \epsilon_k N$.
\end{proposition}

\noindent\textit{Remark.} The purpose of $\epsilon_k$ is to assist in dealing with wraparound issues, which arise from the fact that we are working on $\mathbb{Z}_N$ and not on $[-N,N]$.  Standard sieve theory techniques (in particular the ``fundamental
lemma of sieve theory'') can come very close to providing such a majorant, but the error terms on the pseudorandomness are not
of the form $o(1)$ but rather something like $O(2^{-2^{Ck}})$ or so.  This unfortunately does not quite seem to be good
enough for our argument, which crucially relies on $o(1)$ type decay, and so we have to rely instead of recent arguments of
Goldston and Y{\i}ld{\i}r{\i}m.

\noindent\textit{Proof of Theorem \ref{mainthm} assuming Proposition \ref{prime-majorant}.} 
Let $N$ be a large prime number.  Define the function $f \in L^1(\Z_N)$ by setting $f(n) := 
k^{-1}2^{-k-5}\tilde{\Lambda}(n)$ for $\epsilon_k N \leq n \leq 2\epsilon_k N$ and $f(n) = 0$ otherwise. From Dirichlet's theorem we observe that
$$ \E(f) = \frac{k^{-1}2^{-k-5}}{N} \sum_{\epsilon_k N \leq n \leq 2 \epsilon_k N} \tilde \Lambda(n)
= k^{-1}2^{-k-5} \epsilon_k (1 + o(1)).$$
We now apply Proposition \ref{prime-majorant} and Theorem \ref{main} to conclude that
$$\mathbb{E} \big( f(x)f(x+r) \dots f(x + (k-1)r) \; \big| \; x,r \in \mathbb{Z}_N\big) \geq c(k,k^{-1}2^{-k-5} \epsilon_k) - o(1).$$
Observe that the degenerate case $r=0$ can only contribute at most $O(\frac{1}{N} \log^k N) = o(1)$ to the left-hand side and can
thus be discarded. Furthermore, every progression counted by the expression on the left is not just a progression in $\mathbb{Z}_N$, but a genuine arithmetic progression of integers since $\epsilon_k < 1/k$.  Since the right-hand side is positive (and bounded away from
zero) for sufficiently large $N$, the claim follows from the definition of $f$
and $\tilde \Lambda$.
\endproof\vspace{11pt}

\ni Thus to obtain arbitrarily long arithmetic progressions in the primes, it will suffice to prove Proposition \ref{prime-majorant}.
This will be the purpose of the remainder of this section (with certain number-theoretic computations being deferred to \S \ref{sec9} and the Appendix).\vs

\ni To obtain a majorant for $\tilde \Lambda(n)$, we begin with the well-known formula
\[ \Lambda(n) = \sum_{d|n} \mu(d) \log(n/d) = \sum_{d|n} \mu(d) \log(n/d)_+\]
for the von Mangoldt function, where $\mu$ is the M\"obius function, and $\log(x)_+ $ denotes the positive part of the logarithm, that is to say $\max(\log(x),0)$.  Here and in the sequel $d$ is always understood to be a positive integer.
Motivated by this, we define

\begin{definition}[Goldston--Y{\i}ld{\i}r{\i}m truncated divisor sum]\label{truncated-divisor-sum} Let $R$ be a parameter (in applications it will be a small power of $N$). Define
\[ \Lambda_R(n) := \sum_{\substack{d | n \\ d \leq R}} \mu(d)\log (R/d)
= \sum_{d|n} \mu(d) \log(R/d)_+.\]
\end{definition}

\ni These truncated divisor sums have been studied in several papers, most notably the works of Goldston and Y{\i}ld{\i}r{\i}m 
\cite{goldston-yildirim-old1,goldston-yildirim-old2,goldston-yildirim} concerning the problem of finding small gaps between primes. 
We shall use a modification of their arguments for obtaining asymptotics for these truncated primes to prove that the
measure $\nu$ defined below is pseudorandom.

\begin{definition}\label{mu-def} Let $R := N^{k^{-1}2^{-k-4}}$, and let $\epsilon_k := 1/2^{k}(k+4)!$.  We define the function
$\nu: \Z_N \to \R^+$ by
\[ \nu(n) \; := \; \left\{
\begin{array}{ll}
\frac{\phi(W)}{W} \frac{\Lambda_R(Wn + 1)^2}{\log R} & \hbox{ when } \epsilon_k N \leq n \leq 2\epsilon_k N\\
1 & \hbox{ otherwise}
\end{array}\right.
\]
for all $0 \leq n < N$, where we identify $\{0,\ldots,N-1\}$ with $\Z_N$ in the usual manner.
\end{definition}

\ni This $\nu$ will be our majorant for Proposition \ref{prime-majorant}.  We first verify that it is indeed a majorant.

\begin{lemma}\label{mu-majorises} $\nu(n) \geq 0$ for all $n \in \mathbb{Z}_N$, and furthermore we have $\nu(n) \geq k^{-1}2^{-k-5}\widetilde{\Lambda}(n)$ for all $\epsilon_k N \leq n \leq 2\epsilon_kN$ \emph{(}if $N$ is sufficiently large depending on $k$\emph{)}.
\end{lemma}

\proof The first claim is trivial. The second claim is also trivial unless $Wn+1$ is prime.  From definition of $R$, we see that $Wn+1 > R$ if
$N$ is sufficiently large.  Then the sum over $d | Wn + 1$, $d \leq R$ in \eqref{truncated-divisor-sum} in fact consists of just the one term $d = 1$. Therefore $\Lambda_R(Wn + 1) = \log R$, which means that $\nu(n) = \frac{\phi(W)}{W}\log R \geq k^{-1}2^{-k-5} \widetilde{\Lambda}(n)$
by construction of $R$ and $N$ (assuming $\w(N)$ sufficiently slowly growing in $N$).
\endproof\vspace{11pt}

\ni We will have to wait a while to show that $\nu$ is actually a measure (i.e. it verifies \eqref{numean}). The next proposition will be crucial in showing that $\nu$ has the linear forms property.

\begin{proposition}[Goldston-Y{\i}ld{\i}r{\i}m]\label{GY}  Let $m, t$ be positive integers.  For each $1 \leq i \leq m$, let
$\psi_i(\mathbf{x}) := \sum_{j=1}^t L_{ij} x_j + b_i$, be linear forms with integer coefficients $L_{ij}$ such that $|L_{ij}| \leq \sqrt{\w(N)}/2$ for all $i = 1,\ldots m$ and $j=1,\ldots,t$. We assume that the $t$-tuples $(L_{ij})_{j=1}^t$ are never identically zero, and that no two $t$-tuples
are rational multiples of each other.
Write $\theta_i := W\psi_i + 1$.  Suppose that $B$ is a product $\prod_{i = 1}^t I_i \subset \R^t$ of $t$ intervals $I_i$, each of which having length  at least $R^{10m}$. Then \textup{(}if the function $\w(N)$ is sufficiently slowly growing in $N$\textup{)}
\[ \E( \Lambda_R(\theta_1(\mathbf{x}))^2 \dots \Lambda_R(\theta_m(\mathbf{x}))^2 | \x \in B ) = 
(1 + o_{m,t}(1)) \left(\frac{W \log R}{\phi(W)}\right)^m.
\]
\end{proposition}
\noindent\textit{Remarks.} We have attributed this proposition to Goldston and Y{\i}ld{\i}r{\i}m, because it is a straightforward generalisation of \cite[Proposition 2]{goldston-yildirim}. The $W$-trick makes much of the analysis of the so-called singular series (which is essentially just $(W/\phi(W))^m$ here) easier in our case, but to compensate we have the slight extra difficulty of dealing with forms in several variables.\vs

\ni To keep this paper as self-contained as possible, we give a proof of Proposition \ref{GY}. In \S \ref{sec9} the reader will find a proof which depends on an estimation of a certain contour integral involving the Riemann $\zeta$-function. This is along the lines of \cite[Proposition 2]{goldston-yildirim} but somewhat different in detail. The aforementioned integral is precisely the same as one that Goldston and Y{\i}ld{\i}r{\i}m find an asymptotic for. We recall their argument in the Appendix.\vs

\ni Much the same remarks apply to the next proposition, which will be of extreme utility in demonstrating that $\nu$ has the correlation property (Definition \ref{correlation-condition}).

\begin{proposition}[Goldston-Y{\i}ld{\i}r{\i}m]\label{GY2}
Let $m \geq 1$ be an integer, and let $B$ be an interval of length  at least $R^{10m}$. Suppose that $h_1,\dots,h_m$ are distinct integers satisfying $|h_i| \leq N^2$ for all $1 \leq i \leq m$, and let $\Delta$ denote the integer
\[ \Delta := \prod_{1 \leq i < j \leq m} |h_i - h_j|.\]
Then (for $N$ sufficiently large depending on $m$, and assuming the function $\w(N)$ sufficiently slowly growing in $N$)
\begin{equation}\label{GY2-est}
\begin{split}
 \E( &\Lambda_R(W(x + h_1) + 1)^2 \dots \Lambda_R(W(x + h_m) + 1)^2| x \in B ) \\
&\leq (1 + o_m(1))
\left( \frac{W\log R}{\phi(W)}\right)^m \prod_{p | \Delta}(1 + O_m(p^{-1/2})).
\end{split}
\end{equation}
Here and in the sequel, $p$ is always understood to be prime.
\end{proposition}

\ni Assuming both Proposition \ref{GY} and Proposition \ref{GY2}, we can now
conclude the proof of Proposition \ref{prime-majorant}.  We begin by showing that $\nu$ is indeed a measure.

\begin{lemma}\label{nu-measure} The measure $\nu$ constructed in Definition \ref{mu-def} obeys the estimate
$\mathbb{E}(\nu) = 1 + o(1)$.
\end{lemma}
\proof Apply Proposition \ref{GY} with $m := t := 1$, $\psi_1(x_1) := x_1$ and $B := [\epsilon_k N, 2\epsilon_k N]$ (taking $N$
sufficiently large depending on $k$, of course). Comparing with Definition \ref{mu-def} we thus have
\[ \E (\nu(x) \; | \; x \in [\epsilon_k N, 2\epsilon_k N]) = 1 + o(1).\]
But from the same definition we clearly have
\[ \E (\nu(x) \; | \; x \in \Z_N \backslash [\epsilon_k N, 2\epsilon_k N] ) = 1;\]
Combining these two results confirms the lemma.\endproof\vs

\ni Now we verify the linear forms condition, which is proven in a similar spirit to the above lemma.

\begin{proposition}\label{mu-linear-forms}
The function $\nu$ satisfies the $(k \cdot 2^{k-1},3k-4,k)$-linear forms condition.
\end{proposition}
\proof 
Let $\psi_i(\mathbf{x}) = \sum_{j = 1}^t L_{ij} x_j + b_i$ be linear forms of the type which feature in Definition \ref{linear-forms-condition}. That is to say, we have $m \leq k \cdot 2^{k-1}$, $t \leq 3k - 4$, the $L_{ij}$ are rational numbers with numerator and denominator at most $k$ in absolute value, and none of the $t$-tuples $(L_{ij})_{j=1}^t$ is zero or is equal to a rational multiple of any other.
We wish to show that
\begin{equation}\label{to-prove} \mathbb{E}(\nu(\psi_1(\mathbf{x})) \dots \nu(\psi_m(\mathbf{x})) \; | \; \mathbf{x} \in \mathbb{Z}_N^t) = 1 + o(1).\end{equation}

\ni We may clear denominators and assume that all the $L_{ij}$ are integers, at the expense of increasing the bound on $L_{ij}$ to
$|L_{ij}| \leq (k+1)!$.  Since $\w(N)$ is growing to infinity in $N$, we may assume that 
$(k+1)! < \sqrt{\w(N)}/2$ by taking $N$ sufficiently large. This is required in order to apply Proposition \ref{GY} as we have stated it. \vs

\ni The two-piece definition of $\nu$ in Definition \ref{mu-def} means that we cannot apply Proposition \ref{GY} immediately, and we need the following localization argument. \vs

\ni We chop the range of summation in \eqref{to-prove} into $Q^t$ almost equal-sized boxes, where $Q = Q(N)$ is a slowly growing function of $N$ to be chosen later. Thus let 
\[ B_{u_1,\dots,u_t} = \{ \mathbf{x} \in \Z_N^t : x_j \in [\lfloor u_jN/Q\rfloor,\lfloor(u_j + 1)N/Q\rfloor), j = 1,\dots,t\},\] where the $u_j$ are to be considered $\md{Q}$. Observe that up to negligible multiplicative errors of $1 + o(1)$ (arising because the boxes do not quite have equal sizes) the left-hand side of \eqref{to-prove} can be rewritten as
$$ \E( \E( \nu(\psi_1(\x)) \ldots \nu(\psi_m(\x)) | \x \in B_{u_1,\ldots,u_t} ) | u_1, \ldots, u_t \in \Z_Q ).$$
Call a $t$-tuple $(u_1,\dots,u_t) \in \Z_Q^t$ \textit{nice} if for every $1 \leq i \leq m$,
the sets $\psi_i(B_{u_1,\ldots,u_t})$ are either completely contained in the interval $[\epsilon_k N, 2\epsilon_k N]$ or are completely
disjoint from this interval.  From Proposition \ref{GY} and Definition \ref{mu-def} we observe that
$$ \E( \nu(\psi_1(\x)) \ldots \nu(\psi_m(\x)) | \x \in B_{u_1,\ldots,u_t} ) = 1 + o_{m,t}(1)$$
whenever $(u_1,\ldots,u_t)$ is nice, since we can replace\footnote{There is a technical issue here due to the failure of the quotient map $\Z \to \Z_N$ to be a bijection.  More specifically, the functions $\psi_i(\x)$ only take values in the interval $[\epsilon_k N, 2\epsilon_k N]$ \emph{modulo} $N$, and so strictly speaking one needs to subtract a multiple of $N$ from $\psi_i$ in the formula below.  However, because of the relatively small dimensions of the box $B_{u_1,\ldots,u_t}$, the multiple of $N$ one needs to subtract is independent of $\x$, and so it can be absorbed into the constant term $b_i$ of the affine-linear form $\psi_i$ and thus be harmless.}
 each of the $\nu(\psi_i(\x))$ factors by either
$\frac{\phi(W)}{W \log R} \Lambda_R^2(\theta_i(\x))$ or $1$, and $N/Q$ will exceed $R^{10m}$ for $Q$ sufficiently slowly growing in $N$, by definition of $R$ and the upper bound on $m$.  
When $(u_1,\ldots,u_t)$ is not nice, then we can crudely
bound $\nu$ by $1 + \frac{\phi(W)}{W \log R} \Lambda_R^2(\theta_i(\x))$, multiply out, and apply Proposition \ref{GY} again to obtain
$$ \E( \nu(\psi_1(\x)) \ldots \nu(\psi_m(\x)) | \x \in B_{u_1,\ldots,u_t} ) = O_{m,t}(1) + o_{m,t}(1)$$
We shall shortly show that the proportion of non-nice $t$-tuples $(u_1,\ldots,u_t)$ in $\Z_Q^t$ is at 
most $O_{m,t}(1/Q)$, and thus the left-hand side of \eqref{to-prove} is $1 + o_{m,t}(1) + O_{m,t}(1/Q)$, and the claim follows
by choosing $Q$ sufficiently slowly growing in $N$.\vs

\ni It remains to verify the claim about the proportion of non-nice $t$-tuples.
Suppose $(u_1,\ldots,u_t)$ is not nice. Then there exists $1 \leq i \leq m$ and $\x, \x' \in B_{u_1,\ldots,u_t}$ such that
$\psi_i(\x)$ lies in the interval $[\epsilon_k N, 2\epsilon_k N]$, but $\psi_i(\x')$ does not.  But from definition of $B_{u_1,\ldots,u_t}$(and the boundedness of the $L_{ij}$)  we have
$$ \psi_i(\x), \psi_i(\x') = \sum_{j=1}^t L_{ij}\lfloor Nu_j/Q\rfloor + b_i + O_{m,t}(N/Q). $$
Thus we must have
$$ a \epsilon_k N = \sum_{j=1}^t L_{ij}\lfloor Nu_j/Q\rfloor + b_i + O_{m,t}(N/Q)$$
for either $a=1$ or $a=2$.  Dividing by $N/Q$, we obtain
$$ \sum_{j=1}^t L_{ij} u_j = a \epsilon_k Q + b_i Q/N + O_{m,t}(1) \quad \md{Q}.$$
Since $(L_{ij})_{j=1}^t$ is non-zero, the number of $t$-tuples $(u_1,\ldots,u_t)$ which satisfy this equation is at most
$O_{m,t}(Q^{t-1})$.  Letting $a$ and $i$ vary we thus see that the proportion of non-nice $t$-tuples is at most $O_{m,t}(1/Q)$
as desired (the $m$ and $t$ dependence is irrelevant since both are functions of $k$).
\endproof\vspace{11pt}

\ni In a short while we will use Proposition \ref{GY2} to show that $\nu$ satisfies the correlation condition (Definition \ref{correlation-condition}). Prior to that, however, we must look at the average size of the ``arithmetic'' factor $\prod_{p | \Delta} (1 + O_m(p^{-1/2}))$ appearing in that proposition.

\begin{lemma}\label{additive-weights}  Let $m \geq 1$ be a parameter.
There is a weight function $\tau = \tau_m: \Z \to \R^+$ such that $\tau(n) \geq 1$ for all $n \neq 0$, and such that
for all distinct $h_1,\ldots,h_j \in [\epsilon_k N, 2\epsilon_k N]$ we have
\[ \prod_{p | \Delta}(1 + O_m(p^{-1/2})) \leq \sum_{1 \leq i < j \leq m} \tau(h_i - h_j),\]
where $\Delta$ is defined in Proposition \ref{GY2},
and such that $\E( \tau^q(n) | 0 < |n| \leq N ) = O_{m,q}(1)$ for all $0< q < \infty$.
\end{lemma}

\proof We observe that
$$ \prod_{p|\Delta} (1 + O_m(p^{-1/2})) \leq \prod_{1 \leq i < j \leq m} \bigg(\prod_{p|h_i-h_j} (1 + p^{-1/2}) \bigg)^{O_m(1)}.$$
By the arithmetic mean-geometric mean inequality (absorbing all constants into the $O_m(1)$ factor) we can thus take 
$\tau_m(n) := O_m(1) \prod_{p | n} (1 + p^{-1/2})^{O_m(1)}$
for all $n \neq 0$.  (The value of $\tau$ at 0 is irrelevant for this lemma since we are taking all the $h_i$ to be distinct).
To prove the claim, it thus suffices to show that
$$ \E\bigg( \prod_{p|n} (1 + p^{-1/2})^{O_m(q)} \; \bigg| \; 0 < |n| \leq N \bigg) = O_{m,q}(1) \hbox{ for all } 0 < q < \infty.$$
Since $(1 + p^{-1/2})^{O_m(q)}$ is bounded
by $1 + p^{-1/4}$ for all but $O_{m,q}(1)$ many primes $p$, we have
$$ \E\bigg( \prod_{p|n} (1 + p^{-1/2})^{O_m(q)} \;  \bigg| \; 0 < |n| \leq N \bigg) \leq O_{m,q}(1) 
\E\bigg( \prod_{p|n} (1 + p^{-1/4}) \; \bigg| \; 0 < n \leq N \bigg).$$
But $\prod_{p|n} (1 + p^{-1/4}) \leq \sum_{d|n} d^{-1/4}$, and hence
\begin{eqnarray*} \E\bigg( \prod_{p|n} (1 + p^{-1/2})^{O_m(q)} \; \bigg| \; 0 < |n| \leq N \bigg) & \leq  & O_{m,q}(1)
\frac{1}{2N} \sum_{1 \leq |n| \leq N} \sum_{d|n} d^{-1/4} \\ &
\leq & O_{m,q}(1) \frac{1}{2N} \sum_{d=1}^N \frac{N}{d} d^{-1/4},\end{eqnarray*} which is $O_{m,q}(1)$ as desired.\endproof\vspace{11pt}

\ni We are now ready to verify the correlation condition.

\begin{proposition}\label{mu-cor-con} The measure $\nu$ satisfies the $2^{k-1}$-correlation condition.\end{proposition}
\proof Let us begin by recalling what it is we wish to prove. For any $1 \leq m \leq 2^{k-1}$ and
$h_1,\dots,h_m \in \mathbb{Z}_N$ we must show a bound
\begin{equation}\label{to-prove2} \E\big( \nu(x+h_1) \nu(x+h_2) \ldots \nu(x+h_m) \; \big| \; x \in \mathbb{Z}_N\big)
\leq \sum_{1 \leq i < j \leq m} \tau(h_i-h_j),\end{equation} where the weight function $\tau= \tau_m$ is bounded in $L^q$ for all $q$.\vs

\ni Fix $m$, $h_1, \ldots, h_m$.
We shall take the weight function constructed in Lemma \ref{additive-weights} (identifying $\Z_N$ with the integers between $-N/2$ and $+N/2$), 
and set \[ \tau(0) := \exp(Cm \log N/\log \log N)\]
for some large absolute constant $C$.  From the previous lemma we see that $\E(\tau^q) = O_{m,q}(1)$ for all $q$, since the
addition of the weight $\tau(0)$ at 0 only contributes $o_{m,q}(1)$ at most.\vs

\ni We first dispose of the easy case when at least two of the $h_i$ are equal.  In this case we bound the left-hand side
of \eqref{to-prove} crudely by $\| \nu \|_{L^\infty}^m$.  But from Definitions \ref{truncated-divisor-sum}, \ref{mu-def}
and by standard estimates for the maximal order of the divisor function $d(n)$ we have the crude bound
$\| \nu \|_{L^\infty} \ll \exp(C\log N/\log \log N)$, and the claim follows thanks to our choice of $\tau(0)$.\vs

\ni Suppose then that the $h_i$ are distinct. Since, in \eqref{to-prove2}, our aim is only to get an upper bound, there is no need to subdivide $\mathbb{Z}_N$ into intervals as we did in the proof of Proposition \ref{mu-linear-forms}.
Write 
$$g(n) := \frac{\phi(W)}{W} \frac{\Lambda_R^2(Wn+1)}{\log R} {\bf 1}_{[\epsilon_k N, 2\epsilon_k N]}(n).$$
Then by construction of $\nu$ (Definition \ref{mu-def}), we have
\begin{eqnarray*} &&\mathbb{E}\big(\nu(x + h_1) \dots \nu(x+h_m) \; \big| \; x \in \mathbb{Z}_N\big) \\ && \qquad\qquad\qquad\leq  \mathbb{E}\big((1 + g(x + h_1)) \dots (1 + g(x+h_m)) \; \big| \; x \in \mathbb{Z}_N\big).\end{eqnarray*}
The right-hand side may be rewritten as
\[ \sum_{A \subseteq \{1,\ldots,m\}}  \mathbb{E}\bigg( \prod_{i \in A} g(x+h_i) \; \bigg| \; x \in \mathbb{Z}_N\bigg)\]
(cf. the proof of Lemma \ref{halfway}).  Observe that for $i,j \in A$ we may assume $|h_i - h_j| \leq \epsilon_k N$, since
the expectation vanishes otherwise.  By Proposition \ref{GY2} and Lemma \ref{additive-weights}, we therefore have
\[ \mathbb{E}\bigg( \prod_{i \in A} g(x+h_i) \; \bigg| \; x \in \mathbb{Z}_N\bigg) \leq \sum_{1 \leq i < j \leq m} \tau(h_i - h_j) + o_m(1).\] 
Summing over all $A$, and adjusting the weights $\tau$ by a bounded factor (depending only on $m$ and hence on $k$), we obtain the result.
\endproof\vspace{11pt}

\noindent\textit{Proof of Proposition \ref{prime-majorant}.} This is immediate from Lemma \ref{mu-majorises}, Lemma \ref{nu-measure}, Proposition \ref{mu-linear-forms}, Proposition \ref{mu-cor-con} and the definition of $k$-pseudorandom measure, which is Definition \ref{mu-pseudo-def}.\endproof

\section{Correlation estimates for $\Lambda_R$}\label{sec9}

\ni To conclude the proof of Theorem \ref{mainthm} it remains to verify
Propositions \ref{GY} and \ref{GY2}.  That will be achieved in this section, assuming an estimate (Lemma \ref{zeta-bound}) for a certain class of contour integrals involving the $\zeta$-function. The proof of that estimate is given in the preprint \cite{goldston-yildirim}, and will be repeated in the Appendix for sake of completeness. 
The techniques of this section are also rather close to those in \cite{goldston-yildirim}.  We are greatly indebted to Dan Goldston for sharing this preprint with us.\vspace{11pt}

\noindent\textit{The linear forms condition for $\Lambda_R$.} We begin by proving Proposition \ref{GY}.
Recall that for each $1 \leq i \leq m$ we have a linear form
$\psi_i(\mathbf{x}) = \sum_{j=1}^t L_{ij} x_j + b_i$ in $t$ variables $x_1,\dots,x_t$. The coefficients $L_{ij}$ satisfy $|L_{ij}| \leq \sqrt{\w(N)}/2$, where $\w(N)$ is the function, tending to infinity with $N$, which we used to set up the $W$-trick. We assume that none of the $t$-tuples $(L_{ij})_{j=1}^t$ are zero or are rational multiples of any other. Define $\theta_i := W\psi_i + 1$. \vs

\ni Let $B := \prod_{j = 1}^t I_j$ be a product of intervals $I_j$, each of length at least $R^{10m}$.
We wish to prove the estimate
\[ \E\big( \Lambda_R(\theta_1(\mathbf{x}))^2 \dots \Lambda_R(\theta_m(\mathbf{x}))^2 \; \big| \; \x \in B \big) = (1 + o_{m,t}(1)) 
\left( \frac{W \log R}{\phi(W)} \right)^m.
\]
The first step is to eliminate the role of the box $B$.  We can 
use Definition \ref{truncated-divisor-sum} to expand the left-hand side as
\[ \E \bigg(\prod_{i=1}^m \sum_{\substack{d_i,d'_i \leq R \\ d_i, d'_i | \theta_i(\x)}}
\mu(d_i) \mu(d'_i) \log \frac{R}{d_i} \log \frac{R}{d'_i} \; \bigg| \; \x \in B \bigg)\]
which we can rearrange as
\begin{equation}\label{ds}
 \sum_{d_1,\ldots,d_m,d'_1,\ldots,d'_m \leq R}
\bigg(\prod_{i=1}^m \mu(d_i) \mu(d'_i) \log \frac{R}{d_i} \log \frac{R}{d'_i}\bigg)
\E\bigg( \prod_{i=1}^m {\bf 1}_{d_i, d'_i | \theta_i(\x)} \; \bigg| \; \x \in B \bigg).
\end{equation}
Because of the presence of the M\"obius functions we may assume that all the $d_i$, $d'_i$ are square-free.
Write $D := [d_1,\ldots,d_m,d'_1,\ldots,d'_m]$ to be the least common multiple of the $d_i$ and $d'_i$, thus
$D \leq R^{2m}$.
Observe that the expression $\prod_{i=1}^m {\bf 1}_{d_i, d'_i | \theta_i(\x)}$ 
is periodic with period $D$ in each of the components of $\x$, and can thus can be safely defined on $\Z_D^t$.  
Since $B$ is a product of intervals of length at least $R^{10m}$, we thus see that
$$ \E\bigg( \prod_{i=1}^m {\bf 1}_{d_i, d'_i | \theta_i(\x)} \; \bigg| \; \x \in B \bigg)
= \E\bigg( \prod_{i=1}^m {\bf 1}_{d_i, d'_i | \theta_i(\x)} \; \bigg| \; \x \in \Z_D^t \bigg) + O_{m,t}(R^{-8m}).$$
The contribution of the error term $O_m(R^{-8m})$ to \eqref{ds} can be crudely estimated by $O_{m,t}(R^{-6m} \log^{2m} R)$, which
is easily acceptable.    Our task is thus to show that
\begin{eqnarray}\nonumber
 && \sum_{d_1,\ldots,d_m,d'_1,\ldots,d'_m \leq R}
\bigg(\prod_{i=1}^m \mu(d_i) \mu(d'_i) \log \frac{R}{d_i} \log \frac{R}{d'_i}\bigg)
\E\bigg( \prod_{i=1}^m {\bf 1}_{d_i, d'_i | \theta_i(\x)} \; \bigg| \; \x \in \Z_D^t \bigg)\\
&& \qquad\qquad \qquad\qquad= (1 + o_{m,t}(1))
\left( \frac{W \log R}{\phi(W)} \right)^m.\label{ddr}
\end{eqnarray}
To prove \eqref{ddr}, we shall perform a number of standard manipulations (as in \cite{goldston-yildirim}) to rewrite
the left-hand side as a contour integral of an Euler product, which in turn can be rewritten in
terms of the Riemann $\zeta$-function and some other simple factors.
We begin by using the Chinese remainder theorem (and the square-free nature of
$d_i, d'_i$) to rewrite
$$ \E\bigg( \prod_{i=1}^m {\bf 1}_{d_i, d'_i | \theta_i(\x)} \; \bigg| \; \x \in \Z_D^t \bigg)
= \prod_{p|D} \E\bigg( \prod_{i: p | d_i d'_i} {\bf 1}_{\theta_i(\x) \equiv 0 \mdsub{p}} \; \bigg| \; \x \in \Z_p^t \bigg).$$
Note that the restriction that $p$ divides $D$ can be dropped since the multiplicand is 1 otherwise.
In particular, if we write $X_{d_1,\ldots,d_m}(p) := \{ 1 \leq i \leq m: p | d_i \}$ and
\begin{equation}\label{omega-def}
\omega_X(p) := \E\bigg( \prod_{i \in X} {\bf 1}_{\theta_i(\x) \equiv 0 \mdsub{p}} \; \bigg| \; \x \in \Z_p^t \bigg)
\end{equation} for each
subset $X \subseteq\{1,\ldots,m\}$, then we have
$$ \E\bigg( \prod_{i=1}^m {\bf 1}_{d_i, d'_i | \theta_i(\x)} \; \bigg| \; \x \in \Z_D^t \bigg)
= \prod_p \omega_{X_{d_1,\ldots,d_m}(p) \cup X_{d'_1,\ldots,d'_m}(p)}(p).$$
We can thus write the left-hand side of \eqref{ddr} as
\[ \sum_{d_1,\ldots,d_m,d'_1,\ldots,d'_m \in \Z^+}
\bigg(\prod_{i=1}^m \mu(d_i) \mu(d'_i) (\log \frac{R}{d_i})_+ (\log \frac{R}{d'_i})_+\bigg)
\prod_p \omega_{X_{d_1,\ldots,d_m}(p) \cup X_{d'_1,\ldots,d'_m}(p)}(p).\]
To proceed further, we need to express the logarithms in terms of multiplicative functions of the $d_i$, $d'_i$.  To
this end, we introduce the vertical line contour $\Gamma_1$ parameterised by
\begin{equation}\label{gamma1-def}
\Gamma_1(t) := \frac{1}{\log R} + it; \quad -\infty < t < +\infty
\end{equation}
and observe the contour integration identity
$$ \frac{1}{2\pi i} \int_{\Gamma_1} \frac{x^z}{z^2} \, dz = (\log x)_+$$
valid for any real $x > 0$.
The choice of $\frac{1}{\log R}$ for the real part of $\Gamma_1$ is not currently relevant, but will be convenient later
when we estimate the contour integrals that emerge (in particular, $R^z$ is bounded on $\Gamma_1$, while $1/z^2$ is not too large).
Using this identity, we can rewrite the left-hand side of \eqref{ddr} as
\begin{equation}\label{T-exp} (2\pi i)^{-2m}\int_{\Gamma_1} \dots \int_{\Gamma_1} 
F(z,z') \prod_{j = 1}^m \frac{R^{z_j + z'_j}}{z_j^2z^{\prime 2}_j} \, dz_jdz'_j
\end{equation} 
where there are $2m$ contour integrations in the variables $z_1, \ldots, z_m, z'_1, \ldots, z'_m$ on $\Gamma_1$,
$z := (z_1,\ldots,z_m)$ and $z' := (z'_1,\ldots,z'_m)$, and
\begin{equation}\label{fzz-def}
 F(z,z') :=
  \sum_{d_1,\ldots,d_m,d'_1,\ldots,d'_m \in \Z^+}
\bigg(\prod_{j=1}^m \frac{\mu(d_j) \mu(d'_j)}{d_j^{z_j} d_j^{\prime z'_j}}\bigg)
\prod_p \omega_{X_{d_1,\ldots,d_m}(p) \cup X_{d'_1,\ldots,d'_m}(p)}(p).
\end{equation}
We have changed the indices from $i$ to $j$ to avoid conflict with the square
root of $-1$.
Observe that the summand in \eqref{fzz-def} is a multiplicative function of $D = [d_1,\dots,d_m,d'_1,\dots,d'_m]$ and thus we have  (formally,
at least) the Euler product representation $F(z,z') = \prod_p E_p(z,z')$, where
\begin{equation}\label{ep-def}
 E_p(z,z') := \sum_{X,X' \subseteq\{1,\ldots,m\}} \frac{(-1)^{|X|+|X'|} \omega_{X \cup X'}(p)}{p^{\sum_{j \in X} z_j +
\sum_{j \in X'} z'_j}}.
\end{equation}
From \eqref{omega-def} we have $\omega_\emptyset(p) = 1$ and $\omega_X(p) \leq 1$, and so $E_p(z,z') = 1 + O_\sigma(1/p^\sigma)$
when $\Re(z_j), \Re(z'_j) > \sigma$ (we obtain more precise estimates below).  
Thus this Euler product is absolutely convergent to $F(z,z')$ in the domain $\{ \Re(z_j), \Re(z'_j) > 1\}$ at least.\vs

\ni To proceed further we need to exploit the hypothesis that the linear parts
of $\psi_1, \ldots, \psi_m$ are non-zero and not rational multiples of each other.
This shall be done via the following elementary estimates on $\omega_X(p)$.

\begin{lemma}[Local factor estimate]\label{lem9.2}  If $p \leq \w(N)$, then $\omega_X(p) = 0$ for all non-empty $X$; in particular, $E_p = 1$ when
$p \leq \w(N)$.  If instead $p > \w(N)$, then $\omega_X(p) = p^{-1}$ when $|X|=1$ and $\omega_X(p) \leq p^{-2}$
when $|X| \geq 2$.
\end{lemma}

\proof The first statement is clear, since the maps $\theta_j: \Z_p^t \to \Z_p$ are identically 1 when $p \leq \w(N)$.  The second
statement (when $p > \w(N)$ and $|X|=1$) is similar since in this case $\theta_j$ uniformly covers $\Z_p$.  Now
suppose $p >\w(N)$ and $|X|=2$.  We claim that none of the $s$ pure linear forms $W(\psi_i - b_i)$ is a multiple of any other $\md{p}$. Indeed, if this were so then we should have $L_{ij}L^{-1}_{i'j} \equiv \lambda \md{p}$ for some $\lambda$, and for all $j = 1,\dots,t$. But if $a/q$ and $a'/q'$ are two rational numbers in lowest terms, with $|a|, |a'|, q,q' < \sqrt{\w(N)}/2$, then clearly $a/q \not\equiv a'/q' \md{p}$ unless $a = a'$, $q = q'$. It follows that the two pure linear forms $\psi_i - b_i$ and $\psi_{i'} - b_{i'}$ are rational multiples of one another, contrary to assumption. Thus the set of $\mathbf{x} \in (\mathbb{Z}/p\mathbb{Z})^t$ for which $\theta_i(\mathbf{x}) \equiv 0 \md{p}$ for all $i \in X$ is contained in the intersection of two skew affine subspaces of $(\mathbb{Z}/p\mathbb{Z})^t$, and as such has cardinality at most $p^{t - 2}$.\endproof\vspace{11pt}

\ni This lemma implies, comparing with \eqref{ep-def}, that
\[ E_p(z,z') = 1 - {\bf 1}_{p > \w(N)} \sum_{j = 1}^m (p^{-1 - z_j} + p^{-1 - z'_j} - p^{-1 - z_j - z'_j}) \]
\begin{equation}\label{ep-est}
 \qquad\qquad\qquad\qquad + 
{\bf 1}_{p > \w(N)} \sum_{\substack{X,X' \subseteq \{1,\dots,m\} \\ |X \cup X'| \geq 2}} \frac{ O(1/p^2)}{p^{\sum_{j \in X} z_j + \sum_{j \in X'} z'_j}} ,
\end{equation}
where the $O(1/p^2)$ numerator does not depend on $z$, $z'$.
To take advantage of this expansion, we factorise $E_p = E_p^{(1)} E_p^{(2)} E_p^{(3)}$, where
\begin{align*}
E_p^{(1)}(z,z') &:=
 \frac{E_p(z,z')}{
 \prod_{j=1}^m(1 - {\bf 1}_{p > \w(N)} p^{-1 - z_j}) (1 - {\bf 1}_{p > \w(N)}p^{-1 - z'_j}) (1 - {\bf 1}_{p > \w(N)}p^{-1 - z_j - z'_j})^{-1}}\\
E_p^{(2)}(z,z') &:= \prod_{j=1}^m (1 - {\bf 1}_{p \leq \w(N)}p^{-1 - z_j})^{-1}(1 - {\bf 1}_{p \leq \w(N)} p^{-1 - z'_j})^{-1}(1 - {\bf 1}_{p \leq \w(N)} p^{-1 - z_j - z'_j}) \\
E_p^{(3)}(z,z') &:= \prod_{j=1}^m(1 - p^{-1 - z_j})(1 - p^{-1 - z'_j})(1 - p^{-1 - z_j - z'_j})^{-1}.
\end{align*}
Writing $G_j := \prod_p E_p^{(j)}$ for $j=1,2,3$, one thus has $F = G_1G_2G_3$ (at least for $\Re(z_j), \Re(z'_j)$ sufficiently large).  
If we introduce the Riemann $\zeta$-function $\zeta(s) := \prod_p (1 - \frac{1}{p^s})^{-1}$ then we have
\begin{equation}\label{g3-form}
 G_3(z,z') = \prod_{j=1}^m \frac{\zeta(1 + z_j + z'_j)}{\zeta(1 + z_j)\zeta(1 + z'_j)}
\end{equation}
so in particular $G_3$ can be continued meromorphically to all of $\C^{2m}$.  As for the other two factors, we have
the following estimates which allow us to continue these factors a little bit to the left of the imaginary axes.

\begin{definition}  For any $\sigma > 0$, let $\D^m_\sigma \subseteq\C^{2m}$ denote the domain
\[ \D^m_\sigma := \{z_j,z'_j : -\sigma < \Re(z_j),\Re(z'_j) < 100, j = 1,\dots,m\}.\]
If $G = G(z,z')$ is an analytic function of $2m$ complex variables on $\D^m_\sigma$,
we define the $C^k(\D^m_\sigma)$ norm of $G$ for any integer $k \geq 0$ as
$$ \| G \|_{C^k(\D^m_\sigma)} := \sup_{a_1,\dots,a_m,a'_1,\dots,a'_m}
\big\Vert \big(\frac{\partial}{\partial z_1}\big)^{a_1} \ldots \big(\frac{\partial}{\partial z_m}\big)^{a_m} \big(\frac{\partial}{\partial z'_1}\big)^{a'_1} \ldots \big(\frac{\partial}{\partial z'_m}\big)^{a'_m}
G \big\Vert_{L^\infty(\D^m_\sigma)}$$
where $a_1,\ldots,a_m,a'_1,\dots,a'_m$ range over all non-negative integers with total sum at most $k$.
\end{definition}

\begin{lemma}\label{midnight}  The Euler products $\prod_p E_p^{(j)}$ for $j=1,2$ are absolutely convergent in the domain
$\D^m_{1/6m}$.  In particular, $G_1$, $G_2$ can be continued analytically to this 
domain.  Furthermore, we have the estimates
\begin{align*}
\| G_1 \|_{C^m(\D^m_{1/6m})} &\leq O_m(1) \\
\| G_2 \|_{C^m(\D^m_{1/6m})} &\leq O_{m,\w(N)}(1) \\
G_1(0,0) &= 1 + o_m(1) \\
G_2(0,0) &= (W/\phi(W))^m.
\end{align*}
\end{lemma}

\ni\emph{Remark.} The choice $\sigma = 1/6m$ is of course not best possible, but in fact any small positive quantity depending on $m$ would suffice for our argument here.  The dependence of $O_{m,\w(N)}(1)$ on $\w(N)$ is not important, but one can easily obtain (for instance) growth bounds of the form $\w(N)^{O_m(\w(N))}$.\vs

\ni\proof First consider $j=1$.  From \eqref{ep-est} and Taylor expansion
 we have the crude bound $E_p^{(1)}(z,z') = 1 + O_m(p^{-2+4/6m})$ in $\D^m_{1/6m}$, which gives the desired
convergence and also the $C^m(\D^m_{1/6m})$ bound on $G_1$; the estmate for $G_1(0,0)$ also follows since
the Euler factors $E_p^{(1)}(z,z')$ are identically 1 when $p \leq \w(N)$.  The bound for $G_2$ are easy since
this is just a finite Euler product involving at most $\w(N)$ terms; the formula for $G_2(0,0)$ follows
from direct calculation since $\frac{\phi(W)}{W} = \prod_{p < \w(n)} (1 - \frac{1}{p})$.
\endproof \vs

\ni To estimate \eqref{T-exp}, we now invoke the following contour integration lemma.

\begin{lemma}\label{zeta-bound}\cite{goldston-yildirim}
Let $R$ be a positive real number. Let $G = G(z,z')$ be an analytic function of $2m$ complex variables on the domain $\D^m_\sigma$ for some $\sigma > 0$, and suppose that 
\begin{equation}\label{der-bound}\|G\|_{C^m(\D^m_\sigma)} = \exp(O_{m,\sigma}(\log^{1/3} R)).\end{equation}
 Then
\[ \frac{1}{(2\pi i)^{2m}}\int_{\Gamma_1} \dots \int_{\Gamma_1} G(z,z')\prod_{j=1}^m \frac{\zeta(1+z_j + z'_j)}{\zeta(1 + z_j)\zeta(1 + z'_j)} \frac{R^{z_j + z'_j}}{z_j^2z^{\prime 2}_j} \, dz_jdz'_j  \]
\[ = G(0,\dots,0) \log^m R + \sum_{j=1}^m O_{m,\sigma}(\| G \|_{C^j(\D^m_\sigma)} \log^{m-j} R) + O_{m,\sigma}(e^{-\delta \sqrt{\log R}})\]
for some $\delta = \delta(m) > 0$. 
\end{lemma}

\proof While this lemma is essentially in \cite{goldston-yildirim}, we shall give a complete proof in the Appendix
for sake of completeness. \endproof\vspace{11pt}

\ni We apply this lemma with $G := G_1 G_2$ and $\sigma := 1/6m$.  From Lemma \ref{midnight} and the Leibnitz rule we have the bounds
$$ \| G \|_{C^j(\D^m_{1/6m})} \leq O_{j,m,\w(N)}(1) \hbox{ for all } 0 \leq j \leq m,$$
and in particular we obtain \eqref{der-bound} by choosing $\w(N)$ to grow sufficiently slowly in $N$.  Also we have
$G(0,0) = (1 + o_m(1)) (\frac{W}{\phi(W)})^m$ from that lemma.  We conclude 
(again taking $\w(N)$ sufficiently slowly growing in $N$) 
that the quantity in \eqref{T-exp}
is $(1 + o_m(1)) (\frac{W \log R}{\phi(W)})^m$, as desired. 
This concludes the proof of Proposition \ref{GY}. \endproof\vspace{11pt}

\noindent\textit{Higher order correlations for $\Lambda_R$.}
We now prove Proposition \ref{GY2}, using arguments very similar to those used to prove
Proposition \ref{GY}.  The main differences
here are that the number of variables $t$ is just equal to 1, but on
the other hand all the linear forms are equal to each other, 
$\psi_i(x_1) = x_1$.  In particular, these linear forms are now rational
multiples of each other and so Lemma  \ref{lem9.2} no longer applies.
However, the arguments before that Lemma are still valid; thus
we can still write the left-hand side of \eqref{GY2-est} as
an expression of the form \eqref{T-exp} plus an acceptable error,
where $F$ is again defined by \eqref{fzz-def} and $E_p$ is defined by
\eqref{ep-def}; the difference now is that $\omega_X(p)$ is the quantity
\[ \omega_X(p) := \E\bigg( \prod_{i \in X} {\bf 1}_{W(x + h_i) + 1 \equiv 0 \mdsub{p}} \; \bigg| \; x \in \Z_p \bigg).\]
Again we have $\omega_\emptyset(p) = 1$ for all $p$.  The analogue of Lemma
\ref{lem9.2} is as follows.

\begin{lemma}\label{omega2-lem}  
If $p \leq \w(N)$, then $\omega_X(p) = 0$ for all non-empty $X$; in particular, $E_p = 1$ when
$p \leq \w(N)$.  If instead $p > \w(N)$, then $\omega_X(p) = p^{-1}$ when 
$|X|=1$ and $\omega_X(p) \leq p^{-1}$ when $|X| \geq 2$. Furthermore, if $|X| \geq 2$ then $\omega_X(p) = 0$ unless $p$ divides $\Delta := \prod_{1 \leq i < j \leq s} |h_i-h_j|$.
\end{lemma}

\begin{proof} When $p \leq \w(N)$ then $W(x+h_i)+1 \equiv 1 \md{p}$ and the claim
follows.  When $p > \w(N)$ and $|X| \geq 1$, 
$\omega_X(p)$ is equal to $1/p$ when the residue classes 
$\{ h_i \md{p}: i \in X \}$ are all equal, and zero otherwise, and the
claim again follows.
\end{proof}

\ni In light of this lemma, the analogue of \eqref{ep-est} is now
\begin{equation}\label{ep-est2}
  E_p(z,z') = 1 - {\bf 1}_{p > \w(N)} 
\sum_{j = 1}^m (p^{-1 - z_j} + p^{-1 - z'_j} - p^{-1 - z_j - z'_j}) + 
{\bf 1}_{p > \w(N), p | \Delta} \lambda_p(z,z')
\end{equation}
where $\lambda_p(z,z')$ is an expression of the form
$$ \lambda_p(z,z') =
\sum_{\substack{X,X' \subseteq \{1,\dots,m\} \\ |X \cup X'| \geq 2}} 
\frac{ O(1/p)}{p^{\sum_{j \in X} z_j + \sum_{j \in X'} z'_j}}$$
and the $O(1/p)$ quantities do not depend on $z,z'$.
We can thus factorise \[E_p = E_p^{(0)} E_p^{(1)} E_p^{(2)} E_p^{(3)},\]
where
\begin{align*}
E_p^{(0)} &= 1 + {\bf 1}_{p > \w(N), p | \Delta} \lambda_p(z,z') \\
E_p^{(1)} &= \frac{E_p}{E_p^{(0)}
\prod_{j=1}^m (1 - {\bf 1}_{p > \w(N)}p^{-1-z_j})(1 - {\bf 1}_{p > \w(N)}p^{-1-z'_j})(1 - {\bf 1}_{p > \w(N)}p^{-1-z_j-z'_j})^{-1}} \\
E_p^{(2)} &= \prod_{j=1}^m (1 - {\bf 1}_{p \leq \w(N)}p^{-1-z_j})^{-1}(1 - {\bf 1}_{p \leq \w(N)}p^{-1-z'_j})^{-1}(1 - {\bf 1}_{p \leq w(N)}p^{-1-z_j-z'_j}) \\
E_p^{(3)} &= \prod_{j=1}^m (1 - p^{-1-z_j})(1 - p^{-1-z'_j})(1 - p^{-1-z_j-z'_j})^{-1}.
\end{align*}
Write $G_j = \prod_p E_p^{(j)}$. Then, as before, $F = G_0G_1G_2G_3$ and $G_3$
is given by \eqref{g3-form} as before.  As for $G_0, G_1, G_2$, we have
the following analogue of Lemma \ref{midnight}.

\begin{lemma}\label{noon}  Let $0 < \sigma < 1/6m$.  Then the Euler products $\prod_p E_p^{(l)}$ for $l=0,1,2$ are absolutely convergent in the domain
$\D^m_{\sigma}$.  In particular, $G_0$, $G_1$, $G_2$ can be continued analytically to this 
domain.  Furthermore, we have the estimates
\begin{align}
\| G_0 \|_{C^r(\D^m_{\sigma})} &\leq O_m \left( \frac{\log R}{\log \log R} \right)^r
\prod_{p | \Delta} (1 + O_m(p^{2m\sigma-1})) \quad 
\hbox{ for } 0 \leq r \leq m \label{G1-bound}\\
\| G_0 \|_{C^m(\D^m_{1/6m})} &\leq \exp(O_m(\log^{1/3} R)) \label{G1-bound2}\\
\| G_1 \|_{C^m(\D^m_{1/6m})} &\leq O_m(1) \nonumber\\
\| G_2 \|_{C^m(\D^m_{1/6m})} &\leq O_{m,\w(N)}(1) \nonumber\\
G_0(0,0) &= \prod_{p | \Delta} (1 + O_m(p^{-1/2}) ) \label{G1-bound3}\\
G_1(0,0) &= 1 + o_m(1) \nonumber\\
G_2(0,0) &= (W/\phi(W))^m.\nonumber
\end{align}
\end{lemma}

\proof  The estimates for $G_1$ and $G_2$ proceed exactly as in
Lemma \ref{midnight} (the additional factors of $\lambda_p(z,z')$ which
appear on both the numerator and denominator of $E_p^{(1)}$ cancel to
first order, and thus do not present any new difficulties); it is the
estimates for $G_0$ which are the most interesting.\vs

\ni We begin by proving \eqref{G1-bound}.  Fix $l$.
First observe that $G_0 = \prod_{p | \Delta} E_p^{(0)}$.  Now the number of
primes dividing $\Delta$ is at most $O(\log \Delta / \log \log \Delta)$.
Using the crude bound
\begin{equation}\label{badbound} 
\Delta = \prod_{1 \leq i < j \leq m}|h_i - h_j| 
\leq N^{m^2} \leq R^{O_m(1)},
\end{equation} 
we thus see that the number of factors in the Euler product for $G_0$
is $O_m(\frac{\log R}{\log\log R})$. Upon differentiating $r$ times for any $0 \leq r \leq m$ using the Leibnitz rule, one gets a 
sum of $O_m((\log R/\log \log R)^r)$ terms, each of which
consists of $O_m(\log R/\log\log R)$ factors, each of which is equal to some derivative of
$1 + \lambda_p(z,z')$ of order between 0 and $r$.
On $\D^m_\sigma$, each factor is bounded by $1 + O_m(p^{2m\sigma-1})$ (in fact, the terms containing a non-zero number of derivatives will be much smaller
since the constant term 1 is eliminated).  This gives \eqref{G1-bound}.\vs

\ni Now we prove \eqref{G1-bound2}.  In light of \eqref{G1-bound}, it suffices
to show that
\[
\prod_{p | \Delta} (1  + O_m(p^{2m\sigma - 1}))
\leq \exp( O_m(\log^{1/3} R)).\]
Taking logarithms and using the hypothesis $\sigma < 1/6m$ (and 
\eqref{badbound}), we reduce to showing
\[ \sum_{p | \Delta} p^{-2/3} \leq O(\log^{1/3} \Delta).\]
But there are at most $O(\log \Delta/ \log \log \Delta)$ primes dividing
$\Delta$, hence the left-hand side can be crudely bounded by
\[ \sum_{1 \leq n \leq O(\log \Delta/\log\log\Delta)} n^{-2/3} = O(\log^{1/3} \Delta)\]
as desired.\vs

\ni The bound \eqref{G1-bound3} now follows from the crude
estimate $E_p^{(0)}(z,z') = 1 + O_m( p^{-1/2})$.
\endproof\vspace{11pt}

\ni We now apply Lemma \ref{zeta-bound} with $\sigma := 1/6m$ and
$G := G_0 G_1 G_2$.  Again by the Leibnitz rule we have the bound \eqref{der-bound}, and furthermore
\[ \| G \|_{C^r(\D^m_\sigma)} 
\leq O_m(1) O_{m,\w(N)}(1) \left( \frac{\log R}{\log \log R} \right)^r \prod_{p | \Delta} \left(1 + O_m(p^{-1/2})\right).\]
for all $0 \leq r \leq m$.  From Lemma \ref{noon} and Lemma \ref{zeta-bound}
we can then estimate \eqref{T-exp} as
\begin{align*}
 \leq (1 + o_m(1)) & \left( \frac{W}{\phi(W)} \right)^m \log^m R \prod_{p | \Delta}\left( 1 + O_m(p^{-1/2})\right) \\ & + O_{m,w(N)}\left(  \frac{\log^m R}{\log\log R} \right) \prod_{p | \Delta} \left(1 + O_m(p^{-1/2})\right) + O_m(e^{-\delta\sqrt{\log R}}).\end{align*}
The claim \eqref{GY2-est} then follows by choosing $\w(N)$ (and hence $W$)
sufficiently slowly growing in $N$ (and hence in $R$).
Proposition \ref{GY2} follows.\endproof\vspace{11pt}

\ni\emph{Remark.}  It should be clear that the above argument not only gives
an upper bound for the left-hand side of \eqref{GY2-est}, but in fact
gives an asymptotic, by working out $G_0(0,0)$ more carefully; 
this is worked out in detail (in the $W=1$ case) in \cite{goldston-yildirim}.

\section{Further remarks}

\ni In this section we discuss some extensions and refinements of our main result. First of all, notice that our proof actually shows that that there is some constant $\gamma(k)$ such that the number of $k$-term progressions of primes, all less than $N$, is at least $(\gamma(k) + o(1))N^2/\log^k N$. This is because the error term in \eqref{recurrence} does not actually need to be $o(1)$, but merely
less than $\frac{1}{2} c(k,\delta) + o(1)$ (for instance).  Working backwards through the proof, this eventually reveals that the quantity
$\w(N)$ does not actually need to be growing in $N$, but can instead be a fixed number depending only on $k$ (although this number will
be very large because our final bounds $o(1)$ decayed to zero extremely slowly).  Thus $W$ can be made independent of $N$, and so
the loss incurred by the $W$-trick when passing from primes to primes equal to 1 mod $W$ is bounded uniformly in $N$.
Nevertheless the bound we obtain on $\gamma(k)$ is extremely poor, in part because of the growth of constants in the best known
bounds $c(k,\delta)$ on Szemer\'edi's theorem in \cite{gowers}, but also because we have not attempted to optimise the decay rate of the
$o(1)$ factors and hence will need to take $\w(N)$ to be extremely large.  In the other direction, standard sieve theory arguments
show that the number of $k$-term progressions of primes all less than $N$ are at most $O_k(N^2 / \log^k N)$, and so the lower bounds are only
off by a constant depending on $k$.\vs

\ni As we remarked earlier, our method also extends to prove Theorem \ref{sz-primes}, namely that any subset of the primes with positive relative upper density contains a $k$-term arithmetic progression.   The only significant change\footnote{Also, since we are only assuming positivity of the upper density and not the lower density, we only have good density control for an infinite sequence $N_1, N_2, \ldots \to \infty$ of integers, which may not be prime.  However one can easily use Bertrands postulate (for instance) to make the $N_j$ prime, giving up a factor of $O(1)$ at most.} to the proof is that one must use the pigeonhole principle to replace the residue class $n \equiv 1 \md{W} $ by a more general 
residue class $n \equiv b \md{W}$ for some $b$ coprime
to $W$, since the set $A$ in Theorem \ref{sz-primes} does not need to obey a Dirichlet-type theorem in these residue classes.  However it is easy to verify 
that this does not significantly affect the rest of the argument, and we 
leave the details to the reader.\vs

\ni Applying Theorem \ref{sz-primes} to the set of primes $p \equiv 1 \md{4}$, we obtain the previously unknown fact that there are arbitrarily long progressions consisting of numbers which are the sum of two squares. For this problem, more satisfactory results were known for small $k$ than was the case for the primes. Let $S$ be the set of sums of two squares. It is a simple matter to show
that there are infinitely many $4$-term arithmetic progressions in $S$. Indeed, Heath-Brown \cite{heath-brown2} observed that the numbers $(n-1)^2 + (n-8)^2$, $(n - 7)^2 + (n + 4)^2$, $(n + 7)^2 + (n - 4)^2$ and $(n + 1)^2 + (n + 8)^2$ always form such a progression; in fact, he was able to prove much more, in particular finding an asymptotic for the number of 4-term progressions in $S$, all of whose members are at most $N$
(weighted by $r(n)$, the number of representations of $n$ as the sum of two squares).\vs

\ni It is reasonably clear that our method will produce long arithmetic progressions for many sets of primes for which one can give a lower bound which agrees with some upper bound coming from a sieve, up to a multiplicative constant. Invoking Chen's famous theorem\cite{Chen} to the effect that there are $\gg N/\log^2 N$ primes $p \leq N$ for which $p + 2$ is a prime or a product of two primes, it ought to be a simple matter to adapt our arguments to show that there are arbitrarily long arithmetic progressions $p_1,\dots,p_k$ of primes, such that each $p_i + 2$ is either prime or the product of two primes; indeed there should be $N / \log^{2k} N$ such progressions with entries less than $N$.  Whilst we do not plan\footnote{Very briefly, the idea is to replace the function $\Lambda_R(Wn+1)$ in the definition of the pseudorandom measure $\nu$ with a variant such as $\Lambda_R(Wn+b) \Lambda_R(Wn+b+2)$ for some $1 \leq b < W$ for which $b, b+2$ are both coprime to $W$; one can use Chen's theorem and the pigeonhole principle to locate a $b$ for which this majorant will capture a large number of Chen primes.  We leave the details to the reader.} to write a detailed proof of this fact, we will in \cite{green-tao} give a proof of the case $k = 3$ using harmonic analysis.\vs

\ni The methods in this paper suggest a more general ``transference principle'', in that if a type of pattern (such as an arithmetic progression)
is forced to arise infinitely often within sets of positive density, then it should also be forced to arise infinitely often inside the prime numbers,
or more generally inside any subset of a pseudorandom set (such as the ``almost primes'') of positive relative density.  Thus, for instance,
one is led to conjecture a Bergelson-Leibman type result (cf. \cite{bergelson-leibman}) for primes. That is, one could hope to show that if $F_i : \mathbb{N} \rightarrow \mathbb{N}$ are polynomials with $F(0) = 0$, then there are infinitely many configurations $(a + F_1(d),\dots,a + F_k(d))$ in which all $k$ elements are prime. This however seems to require some modification\footnote{Note added in press: such a result has been obtained by the second author and T. Ziegler, to appear in \emph{Acta Math.}} to our current argument, in large part because
of the need to truncate the step parameter $d$ to be at most a small power of $N$.  In a similar spirit, the work of Furstenberg and Katznelson \cite{fk} on multidimensional analogues of Szemer\'edi's theorem, combined with this transference principle, now suggests that one should be able to show\footnote{Note added in press: such a result has been obtained by the second author, \emph{J. d'Analyse Math\'ematique} \textbf{99} (2006), 109--176.} that the Gaussian primes in $\Z[i]$ contain infinitely many constellations of any prescribed shape, and similarly for other number fields.  Furthermore, the later work of Furstenberg and Katznelson \cite{fk2} on density Hales-Jewett theorems suggests that one could also show 
that for any finite field $F$, the monic irreducible polynomials in $F[t]$ contain affine subspaces over $F$ of arbitrarily high dimension.
Again, these results would require non-trivial modifications to our argument for a number of reasons,
not least of which is the fact that the characteristic factors for these more advanced generalizations of Szemer\'edi's theorem are much 
less well understood.

\appendix
\section{Proof of Lemma \ref{zeta-bound}} \label{app}

\ni In this appendix we prove Lemma \ref{zeta-bound}. This Lemma was essentially proven in \cite{goldston-yildirim}, but for the sake of self-containedness
we provide a complete proof here (following very closely the approach
in \cite{goldston-yildirim}).\vs

\ni Throughout this section, $R \geq 2$, $m \geq 1$, and $\sigma > 0$ will be fixed.  We shall use $\delta > 0$
to denote various small constants, which may vary from line to line (the previous interpretation of $\delta$ as the average value of a function $f$ will
now be irrelevant). 
We begin by recalling the classical zero-free region for the Riemann $\zeta$ 
function.

\begin{lemma}[Zero-free region]\label{basic-zeta-facts} Define the \emph{classical zero free region} $\mathcal{Z}$ to be the closed region
\[ \mathcal{Z} := \{ s \in \C: 10 \geq \Re s \geq 1 - \frac{\beta}{\log(|\Im s| + 2)} \}\]
for some small $0 < \beta < 1$.  Then if $\beta$ is sufficiently small,
$\zeta$ is non-zero and meromorphic in $\mathcal{Z}$ with a simple pole at 1 and no other singularities.
Furthermore we have
the bounds
$$\zeta(s) - \frac{1}{s-1} = O(\log(|\Im s| + 2)); \quad \frac{1}{\zeta(s)} = O(\log(|\Im s| + 2))$$
for all $s \in \mathcal{Z}$.
\end{lemma}

\proof  
See Titchmarsh \cite[Chapter 3]{titchmarsh}.
\endproof\vspace{11pt}

\ni Fix $\beta$ in the above lemma; we may take $\beta$ to be small enough
that $\mathcal{Z}$ is contained in the region where $1-\sigma < \Re(s) < 101$.
We will allow all our constants in the $O()$ notation
to depend on $\beta$ and $\sigma$, and omit explicit mention of these dependencies from our subscripts.\vs

\ni In addition to the contour $\Gamma_1$ defined in \eqref{gamma1-def}, we will need the two further contours $\Gamma_0$
and $\Gamma_2$, defined by
\begin{equation}\label{gamma-1-param} 
\begin{split}
\Gamma_0(t) &:= - \frac{\beta}{\log(|t| + 2)} + it, \qquad -\infty < t < \infty\\
\Gamma_2(t) &:= 1 + it, \quad -\infty < t < \infty.
\end{split}
\end{equation}
Thus $\Gamma_0$ is the left boundary of $\mathcal{Z}-1$ (which therefore lies to the left of the origin), 
while $\Gamma_1$ and $\Gamma_2$ are vertical lines to the right of the origin.  The usefulness of $\Gamma_2$ for us
lies in the simple observation that $\zeta(1 + z + z')$ has no poles when $z \in \mathcal{Z}-1$ and 
$z' \in \Gamma_2$, but we will not otherwise attempt to estimate any integrals on $\Gamma_2$.\vs

\ni We observe the following elementary integral estimates.

\begin{lemma}\label{preliminary-bounds-1} Let $A,B$ be fixed constants with $A > 1$. Then we have the bounds.
\begin{align}
\int_{\Gamma_0} \log^B(|z| + 2)\left| \frac{R^zdz}{z^A}\right| &\leq O_{A,B}(e^{-\delta \sqrt{\log R}});
\label{I1-bound}\\
\int_{\Gamma_1} \log^B(|z| + 2)\left| \frac{R^zdz}{z^2}\right| &\leq O_B(\log R).\label{I2-bound}
\end{align}
Here $\delta = \delta(A,B,\beta) > 0$ is a constant independent of $R$.
\end{lemma}
\proof We first bound the left-hand side of \eqref{I1-bound}.  Substitute in the parametrisation \eqref{gamma-1-param}. Since $\Gamma_0'(t) = O(1)$ and $|z| \gg |t| + \beta$ we have, for any $T \geq 2$,
\begin{eqnarray*} 
\int_{\Gamma_0} \log^B(|z| + 2)\left| \frac{R^zdz}{z^A}\right| & 
\leq & O_B(\int^{\infty}_0 R^{-\beta/(\log(|t| + 2))} \frac{\log^B(|t| + 2)}{(|t| + \beta)^A} dt) \\ & \leq & O_B( \log^BT \int^T_0  R^{-\beta/\log(t+2)} dt + \int^{\infty}_T \frac{\log^Bt}{t^A} dt) \\ & \leq & O_{A,B}(T\log^BT\exp(-\beta \log R/\log T) + T^{1-A}\log^BT).\end{eqnarray*}
Choosing $T = \exp(\sqrt{\beta \log R/2})$ one obtains the claimed bound.
The bound \eqref{I2-bound} is much simpler, and can be obtained by noting that $R^z$ is bounded on $\Gamma_1$,
and substituting in \eqref{gamma1-def} splitting the integrand up into the ranges $|t| \leq 1/\log R$ and $|t| > 1/\log R$.\endproof\vspace{11pt}

\ni The next lemma is closely related to the case $m = 1$ of Lemma \ref{zeta-bound}.

\begin{lemma}\label{preliminary-bounds-2}
Let $f(z,z')$ be analytic in $\mathcal{D}_{\sigma}^1$ and suppose that \[|f(z,z')| \leq \exp(O_m(\log^{1/3}R))\] uniformly in this domain. Then the integral
\[ I := \frac{1}{(2\pi i)^2} \int_{\Gamma_1}\int_{\Gamma_1} f(z,z') \frac{\zeta(1 + z + z')}{\zeta(1 + z)\zeta(1 + z')} \frac{R^{z + z'}}{z^2 z^{\prime 2}} dz dz'  \]
obeys the estimate
$$ I = f(0,0)\log R + 
\frac{\partial f}{\partial z'}(0,0) + \frac{1}{2\pi i} \int_{\Gamma_1} f(z,-z) \frac{dz}{\zeta(1 + z)\zeta(1 - z)z^4}
+ O_{m}(e^{-\delta \sqrt{\log R}})$$
for some $\delta = \delta(\sigma,\beta) > 0$ independent of $R$.
\end{lemma}

\proof  We observe from Lemma \ref{basic-zeta-facts} that we have enough decay of the integrand in the domain $\D_\sigma^1$ to interchange the order of integration, and to shift contours in either one of the variables $z,z'$ while keeping the 
other fixed, without any difficulties when $\Im(z), \Im(z') \to \infty$; the only issue is to keep track of when the contour
passes through a pole of the integrand.  
In particular we can shift the $z'$ contour from $\Gamma_1$ to $\Gamma_2$, since we do not
encounter any of the poles of the integrand while doing so. 
Let us look at the integrand for each fixed $z' \in \Gamma_2$, viewing it as an analytic function of $z$.  We now attempt to shift the $z$ contour of integration to $\Gamma_0$.
In so doing the contour passes just one pole, a simple one at $z = 0$. The residue there is $\frac{1}{2\pi i} \int_{\Gamma_2} f(0,z') \frac{R^{z'}}{z^{\prime 2}} dz'$, and so we have $I = I_1 + I_2$, where
\begin{align*}
I_1 &:= \frac{1}{2\pi i} \int_{\Gamma_2} f(0,z') \frac{R^{z'}}{z^{\prime 2}}dz' \\
I_2 &:= \frac{1}{(2\pi i)^2} \int_{\Gamma_2} \int_{\Gamma_0} f(z,z') \frac{\zeta(1 + z + z') R^{z + z'}}{\zeta(1 + z)\zeta(1 + z') z^2 z^{\prime 2}}dz dz'.
\end{align*}

\ni To evaluate $I_1$, we shift the $z'$ contour of integration to $\Gamma_0$. Again there is just one pole, a double one at $z' = 0$. The residue there is $f(0,0)\log R + \frac{\partial f}{\partial z'}(0,0)$, and so 
\begin{eqnarray*}
I_1 & = & f(0,0)\log R + \frac{\partial f}{\partial z'}(0,0) + \frac{1}{2\pi i}\int_{\Gamma_0} f(0,z') \frac{R^{z'}}{z^{\prime 2}} dz' \\ & = & f(0,0)\log R + \frac{\partial f}{\partial z'}(0,0) + O_m(e^{-\delta \sqrt{\log R}}),\end{eqnarray*}
for some $\delta > 0$,
the latter step being a consequence of our bound on $f$ and \eqref{I1-bound} (in the case $B = 0$).\vs

\ni To estimate $I_2$, we first swap the order of integration and, for each fixed $z$, view the integrand as an analytic function of $z'$. We move the $z'$ contour from $\Gamma_2$ to $\Gamma_0$, this again being allowed since we have sufficient decay in vertical strips as $|\Im z'| \rightarrow \infty$. In so doing we pass exactly two simple poles, at $z' = -z$ and $z' = 0$. The residue at the first is exactly 
\[ \frac{1}{2\pi i} \int_{\Gamma_0} f(z,-z) \frac{dz}{\zeta(1 + z)\zeta(1 - z)z^4},\] which is one of the terms appearing in our formula for $I$.\vs

\ni The residue at $z' = 0$ is
\[ \int_{\Gamma_0} f(z,0) \frac{R^{z}}{z^2} dz,\] which is $O(e^{-\delta \sqrt{\log R}})$ for some $\delta > 0$ by \eqref{I1-bound}. The value of $I_2$ is the sum of these two quantities and the integral over the new contour $\Gamma_0$, which is 
\begin{equation}\label{eq11.1} \int_{\Gamma_0}\int_{\Gamma_0} f(z,z') \frac{\zeta(1 + z + z') R^{z + z'}}{\zeta(1 + z)\zeta(1 + z') z^2 z^{\prime 2}}dz dz'.\end{equation}
In this integrand we have $|f| = \exp(O_m(\log^{1/3} R))$ and, by Lemma \ref{basic-zeta-facts}, $1/|\zeta(1 + z)| \ll \log (|\Im z| + 2)$ and $1/|\zeta(1 + z')| \ll \log (|\Im z'| + 2)$. Assume that $\beta < 1/10$, as we obviously may. We claim that
\begin{equation}\label{troublesome} |\zeta(1 + z + z')| \ll  (1 + |z| + |z'|)^{1/4} \ll (1 + |z|)^{1/4}(1 + |z'|)^{1/4}\end{equation}
for all $z,z' \in \Gamma_0$. Once this is proven it follows from \eqref{I1-bound}, applied with $A = 7/4$ and $A = 2$, that the integral \eqref{eq11.1} is bounded by $O_m(e^{-\delta\sqrt{\log R}})$ for some $\delta > 0$. Now if $1/2 \leq \sigma \leq 1$ and $|t| \geq 1/100$ we have the convexity bound $|\zeta(\sigma + it)| \ll_{\epsilon} |t|^{1 - \sigma + \epsilon}$ (cf. \cite[Chapter V]{titchmarsh}), and so \eqref{troublesome} is indeed true provided that $|\Im(z + z')| \geq 1/100$. However since $z,z' \in \Gamma_0$ one may see that if $|\Im(z)|,|\Im(z')| \leq t$ then $|z + z'| \gg 1/\log(t + 2)$. It follows from Lemma \ref{basic-zeta-facts} that \eqref{troublesome} holds when $|\Im(z + z')| \leq 1/100$ as well.\vs

\ni Thus we now have estimates for $I_1$ and $I_2$ up to errors of $O_m(e^{-\delta\sqrt{\log R}})$. Putting all of this together completes the proof of the lemma.\endproof\vspace{11pt}

\noindent\textit{Proof of Lemma \ref{zeta-bound}.}  Let $G = G(z,z')$ be
an analytic function of $2m$ complex variables on the domain $\D^m_\sigma$
obeying the derivative bounds \eqref{der-bound}.  We will allow all our implicit constants in the $O()$ notation to
depend on $m$, $\beta$, $\sigma$.
We are interested in the integral
\[ I(G,m) := \frac{1}{(2\pi i)^{2m}}\int_{\Gamma_1} \dots \int_{\Gamma_1} G(z,z')\prod_{j=1}^s \frac{\zeta(1+z_j + z'_j)}{\zeta(1 + z_j)\zeta(1 + z'_j)} \frac{R^{z_j + z'_j}}{z_j^2z^{\prime 2}_j} \, dz_jdz'_j, \]
and wish to prove the estimate
\[ I(G,m) := G(0,\dots,0) (\log R)^{m} + \sum_{j=1}^mO(\|G\|_{C^j(\mathcal{D}^s_\sigma)}(\log R)^{m-j}) + O(e^{-\delta\sqrt{\log R}}).\]
The proof is by induction on $m$.
The case $m = 1$ is a swift deduction from Lemma \ref{preliminary-bounds-2}, the only issue being an estimation of the term 
\[ \frac{1}{2\pi i} \int_{\Gamma_0} G(z_1,-z_1) \frac{dz_1}{\zeta(1 + z_1)\zeta(1 - z_1)z_1^4}.\]
It is not hard to check (using Lemma \ref{basic-zeta-facts}) that 
\begin{equation}\label{zb1} \int_{\Gamma_0} \left| \frac{dz_1}{\zeta(1 + z_1)\zeta(1 - z_1) z_1^4} \right| = O(1),\end{equation}
and so this term is $O(\sup_{z \in \mathcal{D}^1_\sigma} |G(z)|) = O(\|G\|_{C^1(\mathcal{D}_{\sigma}^1)})$.\vs

\ni Suppose then that we have established the result for $m \geq 1$ and wish to deduce it for $m+1$. 
Applying Lemma \ref{preliminary-bounds-2} in the variables $z_{m+1},z'_{m+1}$, we get $I(G,m+1) = $
\begin{align*}
 & \frac{\log R}{(2\pi i)^{2m}}\int_{\Gamma_1}\dots\int_{\Gamma_1} G(z_1,\dots,z_m,0,z'_1,\dots,z'_m,0) \prod_{j=1}^m \frac{\zeta(1 + z_j + z'_j)}{\zeta(1 + z_j)\zeta(1 + z'_j)}\frac{R^{z_j + z'_j}}{z_j^2z^{\prime 2}_j} dz_jdz'_j \\
& + \frac{1}{(2\pi i)^{2m}}\int_{\Gamma_1}\dots\int_{\Gamma_1} H(z_1,\dots,z_m,z'_1,\dots,z'_m)\prod_{j=1}^m\frac{\zeta(1 + z_j + z'_j)}{\zeta(1 + z_j)\zeta(1 + z'_j)}\frac{R^{z_j + z'_j}}{z_j^2z^{\prime 2}_j}dz_jdz'_j \\
& + O(e^{-\delta\sqrt{\log R}}) \\
=&  I(G(z_1,\dots,z_m,0,z'_1,\dots,z'_m,0),m) \log R
+ I(H,m) + O(e^{-\delta \sqrt{\log R}})
\end{align*}
where $\delta > 0$ and $H: \D^m_\sigma \to \C$ is the function
\[ H(z_1,\dots,z_m,z'_1,\dots,z'_m) := \frac{\partial G}{\partial z'_{m+1}}(z_1,\dots,z_m,0,z'_1,\dots,z'_m,0) \qquad\qquad\] \[ \qquad + \frac{1}{2\pi i} \int_{\Gamma_0} G(z_1,\dots,z_m,z_{m+1},z'_1,\dots,z'_m,-z_{m+1}) \frac{dz_{m+1}}{\zeta(1 + z_{m+1})\zeta(1 - z_{m+1})z_{m+1}^4}.\]
The error term $O(e^{-\delta\sqrt{\log R}})$ which we claim here arises by applying \eqref{der-bound} and several
applications of \eqref{I2-bound}.\vs

\ni Now both of the functions $G(z_1,\dots,z_m,0,z'_1,\dots,z'_m,0)$ and $H(z_1,\dots,z_m,z'_1,\dots,z'_m)$ are analytic on
$\D^m_\sigma$ and (appealing to \eqref{zb1}) we have
$\|H\|_{C^j(\mathcal{D}^m_\sigma)} = O_m(\|G\|_{C^{j+1}(\mathcal{D}^{m+1}_\sigma)})$ for
$0 \leq j \leq m$. Using the inductive hypothesis,
we therefore obtain $I(G,m+1) = $
\begin{eqnarray*}
&  & G(0,\dots,0) (\log R)^{m + 1} + \sum_{j=1}^m O_m(\|G(\cdot,0,\cdot,0)\|_{C^j(\mathcal{D}^m_\sigma)}(\log R)^{m + 1 -j}) \\ & &  + H(0,\dots,0) (\log R)^{m} + \sum_{j=1}^m O_m(\|H\|_{C^j(\mathcal{D}^m_\sigma)}(\log R)^{m  -j}) + O(e^{-\delta\sqrt{\log R}})
 \\ & = & G(0,\dots,0)(\log R)^{m+1} + \sum_{j=1}^m O_m(\|G\|_{C^j(\mathcal{D}^{m+1}_\sigma)}(\log R)^{m + 1 -j}) \\ & & + H(0,\dots,0) (\log R)^{m} + \sum_{j=1}^m O_m(\|G\|_{C^{j+1}(\mathcal{D}^{m+1}_\sigma)}(\log R)^{m  -j}) + O(e^{-\delta\sqrt{\log R}})\\
 & = & G(0,\dots,0)(\log R)^{m+1} + \sum_{j=1}^{m+1} O_m(\|G\|_{C^j(\mathcal{D}^{m+1}_\sigma)}(\log R)^{m + 1 -j}) + O(e^{-\delta\sqrt{\log R}}),\end{eqnarray*}which is what we wanted to prove.\endproof

\providecommand{\bysame}{\leavevmode\hbox to3em{\hrulefill}\thinspace}
\providecommand{\MR}{\relax\ifhmode\unskip\space\fi MR }
\providecommand{\MRhref}[2]{%
  \href{http://www.ams.org/mathscinet-getitem?mr=#1}{#2}
}
\providecommand{\href}[2]{#2}


\begin{thebibliography}{10}

\bibitem{assani}
I. Assani, \emph{Pointwise convergence of ergodic averages along cubes}, preprint.

\bibitem{balog1} A. Balog, \emph{Linear equations in primes,} Mathematika \textbf{39} (1992) 367--378.

\bibitem{balog2} \bysame, \emph{Six primes and an almost prime in four linear equations,} Can. J. Math. \textbf{50} (1998), 465--486.

\bibitem{bergelson-leibman} V. Bergelson and A. Leibman, \emph{Polynomial extensions of van der Waerden's and Szemer\'edi's theorems,} J. Amer. Math. Soc. \textbf{9} (1996), 725--753.

\bibitem{bourg2} J. Bourgain, \emph{A Szemer\'edi-type theorem for sets of positive density in $\mathbb{R}^k$,} Israel J. Math \textbf{54} (1986), no. 3, 307--316.

\bibitem{bourgain-triples}
\bysame, \emph{On triples in arithmetic progression}, GAFA \textbf{9} (1999), 968--984.

\bibitem{Chen} J.-R. Chen, \emph{On the representation of a large even integer as the sum of a prime and a product of at most two primes,} Sci. Sinica \textbf{16} (1973), 157--176.

\bibitem{chowla}
S. Chowla, \emph{There exists an infinity of 3---combinations of primes in A. P.},
 Proc. Lahore Philos. Soc. \textbf{6}, (1944). no. 2, 15--16.

\bibitem{erdos}
P. Erd\H{o}s, P. Tur\'an, \emph{On some sequences of integers}, J. London Math. Soc. \textbf{11} (1936), 261--264.

\bibitem{furst}
H. Furstenberg, \emph{Ergodic behavior of diagonal measures and a theorem of Szemer\'edi on arithmetic progressions}, J. Analyse Math. \textbf{31} (1977), 204--256.

\bibitem{fk}
H. Furstenberg, Y. Katznelson, \emph{An ergodic Szemer\'edi theorem for commuting transformations}. J. Analyse Math. \textbf{34} (1978), 275--291.

\bibitem{fk2}
H. Furstenberg, Y. Katznelson, \emph{A density version of the Hales-Jewett theorem}, J. d'Analyse Math. \textbf{57} (1991), 64--119.

\bibitem{FKO}
H. Furstenberg, Y. Katznelson and D. Ornstein, \emph{The ergodic-theoretical proof of Szemer\'edi's theorem,} Bull. Amer. Math. Soc. \textbf{7} (1982), 527--552.

\bibitem{furst-weiss}
H. Furstenberg, B. Weiss, \emph{A mean ergodic theorem for $1/N \sum_{n=1}^N f(T^n x) g(T^{n^2} x)$}, Convergence in ergodic theory and probability (Columbus OH 1993), 193--227, Ohio State Univ. Math. Res. Inst. Publ., 5. de Gruyter, Berlin, 1996.

\bibitem{goldston-yildirim-old1} D. Goldston and C.Y. Y{\i}ld{\i}r{\i}m \emph{Higher correlations of divisor sums related to primes, I: Triple correlations,} Integers \textbf{3} (2003) A5, 66pp.

\bibitem{goldston-yildirim-old2} \bysame, \emph{Higher correlations of divisor sums related to primes, III: $k$-correlations,} preprint (available at AIM preprints)

\bibitem{goldston-yildirim} \bysame, \emph{Small gaps between primes, I,} preprint available at \\
\texttt{http://www.arxiv.org/abs/math/0504336.}

\bibitem{gowers-4}
T. Gowers, \emph{A new proof of Szemer\'edi's theorem for arithmetic 
progressions of length four}, GAFA \textbf{8} (1998), 529--551.
 
\bibitem{gowers}
\bysame, \emph{A new proof of Szemer\'edi's theorem}, GAFA \textbf{11} (2001), 465-588.

\bibitem{gowers-reg}
\bysame, \emph{Hypergraph regularity and the multidimensional Szemer\'edi theorem,} preprint.

\bibitem{green}
B.J. Green, \emph{Roth's theorem in the primes,} Ann. Math. \textbf{161} (2005), no. 3, 1609--1636.

\bibitem{green-reg}
\bysame, \emph{A Szemer\'edi-type regularity lemma in abelian groups,} GAFA \textbf{15} (2005), no. 2, 340--376.

\bibitem{green-tao}
B.J. Green and T. Tao, \emph{Restriction theory of Selberg's sieve, with applications,} J. Th\'eorie des Nombres de Bordeaux \textbf{18} (2006), 147--182.

\bibitem{hardy-littlewood} G.H. Hardy and J.E. Littlewood \emph{Some problems of ``partitio numerorum''; III: On the expression of a number as a sum of primes,} Acta Math. \textbf{44} (1923), 1--70

\bibitem{heath-brown1}
D.R. Heath-Brown, \emph{Three primes and an almost prime in arithmetic progression,} J. London Math. Soc. (2) \textbf{23} (1981), 396--414.

\bibitem{heath-brown2}
\bysame, \emph{Linear relations amongst sums of two squares,} Number theory and algebraic geometry --- to Peter Swinnerton-Dyer on his 75th birthday, CUP (2003). 

\bibitem{host-kra1} B. Host and B. Kra, \emph{Convergence of Conze-Lesigne averages,} Ergodic Theory and Dynamical Systems \textbf{21} (2001), no. 2, 493--509.

\bibitem{host-kra2}
\bysame, \emph{Non-conventional ergodic averages and nilmanifolds,} Ann. Math. \textbf{161} (2005), no. 1, 397--488.

\bibitem{host-kra3} \bysame, \emph{Convergence of polynomial ergodic averages,} Israel. Jour. Math. \textbf{149} (2005), 1--19.

\bibitem{klr}
Y. Kohayakawa, T. Luczsak, V. R\"odl, \emph{Arithmetic progressions of length three in subsets of a random set}, 
Acta Arith. \textbf{75} (1996), no. 2, 133--163.

\bibitem{laba-lacey} I. {\L}aba and M. Lacey, \emph{On sets of integers not containing long arithmetic progressions,} unpublished. Available at
\texttt{http://www.arxiv.org/pdf/math.CO/0108155}.

\bibitem{moran-pritchard-thyssen} A. Moran, P. Pritchard and A. Thyssen, \emph{Twenty-two primes in arithmetic progression,} Math. Comp. \textbf{64} (1995), no. 211, 1337--1339.

\bibitem{ramare} O. Ramar\'e, \emph{On Snirel'man's constant,} Ann. Scu. Norm. Pisa \textbf{21} (1995), 645--706.

\bibitem{ramare-ruzsa} O. Ramar\'e and I.Z. Ruzsa, \emph{Additive properties of dense subsets of sifted sequences,} J. Th. Nombres de Bordeaux \textbf{13} (2001) 559--581.

\bibitem{rankin} R. Rankin, \emph{Sets of integers containing not more than a given number of terms in arithmetical 
progression.} Proc. Roy. Soc. Edinburgh Sect. A, \textbf{65} 1960/1961 332--344 (1960/61).

\bibitem{roth}
K.F. Roth, \emph{On certain sets of integers}, J. London Math. Soc. \textbf{28} (1953), 245-252.

\bibitem{szemeredi-4}
E. Szemer\'edi, \emph{On sets of integers containing no four elements in arithmetic progression},
Acta Math. Acad. Sci. Hungar. \textbf{20} (1969), 89--104.

\bibitem{szemeredi}
\bysame, \emph{On sets of integers containing no $k$ elements in arithmetic progression},
Acta Arith. \textbf{27} (1975), 299--345.

\bibitem{szem-reg}
\bysame, \emph{Regular partitions of graphs,} in ``Proc. Colloque Inter. CNRS'' (J.-C. Bermond, J.-C. Fournier, M. Las Vergnas, D. Sotteau, eds.) (1978), 399--401.

\bibitem{tao:ergodic}
T. Tao, \emph{A quantitative ergodic theory proof of Szemer\'edi's theorem}, to appear in Electronic J. Combinatorics.

\bibitem{titchmarsh}
E.C. Titchmarsh, \emph{The theory of the Riemann zeta function,} Oxford University Press, 2nd ed, 1986.

\bibitem{van-der-corput}
J.G. van der Corput, \emph{\"Uber Summen von Primzahlen und Primzahlquadraten,} Math. Ann. \textbf{116} (1939), 1--50.

\bibitem{varnavides} P. Varnavides, \emph{On certain sets of positive density,} J. London Math. Soc. \textbf{34} (1959) 358--360.

\bibitem{ziegler} T. Ziegler, \emph{Universal characteristic factors and Furstenberg averages}, J. Amer. Math. Soc.  \textbf{20}  (2007),  no. 1, 53--97.

\bibitem{ziegler2}
\bysame, \emph{A non-conventional ergodic theorem for a nilsystem}, A non-conventional ergodic theorem for a nilsystem.  Ergodic Theory Dynam. Systems  \textbf{25}  (2005),  no. 4, 1357--1370.

\end{thebibliography}
     \end{document}